\documentclass[a4paper,11pt, reqno]{amsart}
\usepackage{graphicx}
\usepackage{standalone}
\usepackage{amsmath,amsthm,amssymb}
\usepackage{mathtools}
\usepackage[hidelinks]{hyperref}
\usepackage{cleveref}
\usepackage[a4paper, total={6in, 8in}]{geometry}
\usepackage{todonotes}
\usepackage{subcaption}
\usepackage{enumitem}
\usepackage[foot]{amsaddr}

\usetikzlibrary{decorations.pathreplacing}
\usetikzlibrary{calc}

\hypersetup{
    colorlinks=false,
    linkcolor=blue,
    filecolor=blue,      
    urlcolor=blue,
    citecolor=blue,
    pdftitle={lrp_project},
    pdfpagemode=FullScreen
}
    
\numberwithin{equation}{section}
\theoremstyle{plain}

\newtheorem{theorem}{Theorem}[section]
\newtheorem{proposition}[theorem]{Proposition}
\newtheorem{corollary}[theorem]{Corollary}
\newtheorem{lemma}[theorem]{Lemma}

\crefname{theorem}{Theorem}{Theorems}
\Crefname{theorem}{Theorem}{Theorems}

\crefname{Proposition}{Proposition}{propositions}
\Crefname{proposition}{Proposition}{Propositions}

\crefname{corollary}{Corollary}{Corollaries}
\Crefname{corollary}{Corollary}{Corollaries}

\crefname{lemma}{Lemma}{Lemmas}
\Crefname{lemma}{Lemma}{Lemmas}

\theoremstyle{remark}
\newtheorem{remark}[theorem]{Remark}

\newcommand{\LRP}{long-range percolation}
\newcommand{\BNNP}{Bernoulli nearest-neighbour percolation}

\renewcommand{\epsilon}{\varepsilon}
\newcommand{\aut}{\mathrm{Aut}}
\newcommand{\stab}{\mathrm{Stab}}
\newcommand{\jump}[1]{\overset{r_{#1}}{\sim}}
\newcommand{\normone}[1]{\lVert #1 \rVert_{1}}

\newcommand{\cpei}{c_{\mathrm{pei}}}
\newcommand{\ciso}{c_{\mathrm{iso}}}

\newcommand{\charf}{\mathbf{1}}

\newcommand{\PP}{\mathcal{P}}

\newcommand{\jmax}{J^{\star}}

\newcommand{\mm}{\mathbf{m}}
\newcommand{\rr}{\mathbf{r}}
\newcommand{\ff}{\mathbf{f}}

\newcommand{\cvol}{c_G}
\newcommand{\Cvol}{C_G}

\newcommand{\Z}{\mathbf{Z}}
\newcommand{\R}{\mathbf{R}}
\newcommand{\N}{\mathbf{N}}

\renewcommand{\P}{\mathbf{P}}
\newcommand{\E}{\mathbf{E}}

\newcommand{\FF}{\mathcal{F}}
\newcommand{\HH}{\mathcal{H}}
\newcommand{\BB}{\mathcal{B}}
\newcommand{\MM}{\mathcal{M}}
\newcommand{\RR}{\mathcal{R}}
\renewcommand{\AA}{\mathcal{A}}
\renewcommand{\SS}{\mathcal{S}}

\newcommand{\rmax}{m}

\newcommand{\cluster}{K}

\newcommand{\expdeg}[1]{\E_{\beta}\left[\deg(#1)\right]}
\newcommand{\expdegrad}[2]{\E_{\beta}\left[\deg(#1,#2)\right]}
\newcommand{\expleng}[1]{\E_{\beta}\left[\mathrm{length}(#1)\right]}

\newcommand{\ablocks}{\mathbf{A}}
\newcommand{\bblocks}{\mathbf{B}}

\title[Finite Cluster radius of supercritical LRP]{Supercritical long-range percolation on graphs of polynomial growth: the truncated one-arm exponent}
\date{\today}
\author[Y. M. Alonso]{Yago Moreno Alonso$^1$}
\address{$^1$Delft Institute of Applied Mathematics, Delft University of Technology.}
\author[J. Komj\'athy]{J\'ulia Komj\'athy$^2$}
\address{$^2$Delft Institute of Applied Mathematics, Delft University of Technology.}

\email{Y.M.A.MorenoAlonso@tudelft.nl, j.komjathy@tudelft.nl}

\begin{document}

\begin{abstract}
    We consider supercritical long-range percolation on transitive graphs of polynomial growth. In this model, any two vertices $x$ and $y$ of the underlying graph $G$ connect by a direct edge with probability $1-\exp(-\beta J(x,y))$, where $J(x,y)$ is a function that is invariant under the automorphism group of $G$, and we assume that $J$ decays polynomially with the graph distance between $x$ and $y$. 
    We give up-to-constant bounds on the decay of the radius of finite cluster for $\beta > \beta_c$.
    In the same setting, we also give upper and lower bounds on the tail volume of finite clusters. 
    The upper and lower bounds are of matching order, conjecturally on sharp volume bounds for spheres in transitive graphs of polynomial growth.
    As a corollary, we obtain a lower bound on the anchored isoperimetric dimension of the infinite component.
\end{abstract}

\maketitle

{\footnotesize
    \hspace{1em}
    
    \hspace{1em}
    MSC2020 Class: 82B43, 20F65.
}

\section{Introduction}
Let $G = (V,E)$ be a transitive graph of polynomial growth with a fixed origin $o$.
We write $B(r)$ for the ball of radius $r$ centred at the origin with respect to the graph metric $d_G$.
It is a consequence of the theorems of Gromov, Trofimov, Bass and Guivarc'h \cite{gromov_groups_1981,trofimov_graphs_1984,bass_degree_1972,guivarch_groupes_1970} that there exist $d \in \N$ and $c_G,C_G > 0$ such that
\begin{equation}
    \label{eq:ball_volume_bounds}
    c_G 
    r^d
    \leq
    \# B(r)
    \leq 
    C_G 
    r^d
\end{equation}
for all $r \in \N$. 
Let $J : V \times V \to \R_+$ be \textbf{transitive}, meaning that for any graph automorphism $\gamma$ of $G$ we have $J(\gamma(x), \gamma(y)) = J(x,y)$.
For $\beta \geq 0$, \textbf{\LRP{}} on $G$ with kernel $J$ is the random graph with vertex set $V$ where we include an edge between any pair of distinct vertices $x,y \in V$ independently at random with probability $1 - \exp(-\beta J(x,y))$.
We are mostly interested in the case where $J(x,y) = \Theta(d_G(x,y)^{-d \alpha})$ with $\alpha > 1$, meaning that there exist $c_J,C_J,R_J > 0$ such that 
\begin{equation}
\label{eq:polynomial_kernel}
    c_J 
    d_G(x,y)^{-d \alpha}
    \leq 
    J(x,y) 
    \leq 
    C_J d_G(x,y)^{-d\alpha}
\end{equation}
for all $x,y \in V$ with $d_G(x,y) \geq R_J$.
 Such a choice ensures that the kernel is \textbf{integrable}, meaning that $\sum_{y \in G \setminus \{x\}} J(x,y) < \infty$ for all $x \in G$, see \Cref{lem:finite_expdeg} below.
We write $\P_{\beta}$ for the law of the resulting random graph and $\E_{\beta}$ for the expectation operator. We refer to the connected components of this random graph as \textbf{clusters}, we write $\cluster$ for the cluster of the origin, and we write $\# \cluster$ for its cardinality. 
This model undergoes a phase transition at the \textbf{critical value}
\begin{equation}
    \label{eq:critical_parameter}
    \beta_c 
    = 
    \inf
    \{
        \beta \geq 0 : \P_{\beta}(\# K = \infty) > 0 
    \}
\end{equation}
which satisfies $0 < \beta_c < \infty$ if $d \geq 2$ and $\alpha > 1$ or $d = 1$ and $\alpha \in (1,2]$ \cite{schulman_long_1983,newman_one-dimensional_1986}.
In this paper we are interested in the supercritical regime when $\beta > \beta_c$.
In particular, we consider the probability of the \textbf{truncated one-arm} event
\begin{equation}
    \P_{\beta}
    \left(
        o \leftrightarrow B(r)^c, 
        \# \cluster < \infty
    \right)
\end{equation}
that there exists an open path from the origin to the complement of the ball $B(r)$ and the cluster of the origin $K$ is finite. 
In \BNNP{} $\{o \leftrightarrow B(r)^c, \# K < \infty\}$ implies $\{r \leq \# K < \infty \}$, namely the truncated one-arm event implies that there are at least $r+1$ open edges in the cluster of the origin.
This eventually leads to the exponential decay of the truncated one-arm event, due to \cite{kesten_critical_1980,chayes_correlation_1989,contreras_supercritical_2024}.
In \LRP{} however, edge lengths are unbounded, and a single long edge from $o$ to $B(r)^c$ already occurs with probability polynomial in $r$.

\begin{lemma}
    \label[lemma]{lem:one_arm_lower_bound}
    Let $G$ be a transitive graph of polynomial growth with $d \geq 1$, and suppose that $J : V \times V \to \R_+$ is a transitive kernel with $J(x,y) = \Theta(d_G(x,y)^{-d \alpha})$ with $\alpha > 1$.
    Let $\beta > \beta_c$.
    Then there exists $c = c(\beta,G,J) > 0$ such that for all $r \in \N$
    \begin{equation}
        c
        r^{d(1-\alpha)}
        \leq 
        \P_{\beta}
        \left(
            o 
            \leftrightarrow 
            B(r)^c,
            \# \cluster
            < \infty
        \right).
    \end{equation}
\end{lemma}

It could be that potential one-arm paths with multiple shorter edges push the probability of the truncated one-arm event up to a different order, or that they cause lower-order correction terms. 
In this paper we show that the order $r^{d(1-\alpha)}$ given by this simple lower bound is in fact the correct order for the probability of the truncated one-arm event. 

\subsection{Main results}
Our first result is the decay of the radius of a finite cluster for $d \geq 2$ and $\beta > \beta_c$ sufficiently large.

\begin{theorem}
\label[theorem]{thm:arm_exponent_1}
    Let $G$ be a transitive graph of polynomial growth with $d \geq 2$, and suppose that $J : V \times V \to \R_+$ is a transitive kernel with $J(x,y) = \Theta(d_G(x,y)^{-d \alpha})$ with $\alpha > 1$.
    Let $\beta > \beta_c$ be sufficiently large.
    Then there exist $C = C(\beta,G,J) > 0$ such that for all $r \in \N$
    \begin{equation}\label{eq:theorem-1-equation}
        \P_{\beta}
        \left(
            o 
            \leftrightarrow 
            B(r)^c,
            \# \cluster
            < \infty
        \right) 
        \leq 
        C
        r^{d(1-\alpha)}.
    \end{equation}
\end{theorem}

At the heart of the proof is a pigeonhole argument which roughly says that if $\# K \leq k$ for $k \in \N$, then at least one of the edges in the one-arm path must have length at least $r/k$.
The main term for the decay of the truncated one-arm comes from this long edge: an upper bound with polylogarithmic correction terms is proved in \Cref{thm:lazy_one_arm}, and a more sophisticated pigeonhole argument gives the sharp result above.
In order for these arguments to work, we must have good bounds on the distribution of finite clusters. 
Indeed, for $k \in \N$ we have by a union bound that
\begin{align}
\label{eq:one_arm_decomposition}
    \P_{\beta}
    \left(
        o
        \leftrightarrow 
        B(r)^c, 
        \# \cluster
        < 
        \infty
    \right)
    \leq
    \P_{\beta}
    \left(
        o 
        \leftrightarrow 
        B(r)^c, 
        \# \cluster
        \leq 
        k
    \right)
    +
    \P_{\beta}
    \left(
        k 
        \leq 
        \# \cluster
        < 
        \infty
    \right),
\end{align}
and to prove \Cref{thm:arm_exponent_1} we must show that $\P_{\beta} \left( k \leq \# \cluster < \infty \right) = o(r^{d(1-\alpha)})$ for some choice of $k=k(r)$.
In this paper we prove the subexponential decay of the distribution of finite clusters for $d \geq 2$, $\alpha > 1 + 1/d$, and $\beta > \beta_c$ sufficiently large, which is new for transitive graphs of polynomial growth.

\begin{theorem}[Cluster-size decay]
    \label[theorem]{thm:cluster_size_strong_decay}
    Let $G$ be a transitive graph of polynomial growth with $d \geq 2$, and suppose that $J : V \times V \to \R_+$ is a transitive kernel with $J(x,y) =O(d_G(x,y))^{-d \alpha})$ with $\alpha > 1 + 1/d$.
    Let $\beta > \beta_c$ be sufficiently large.
    Then there exists $A = A(G,J) > 0$ such that 
    for all $k \in \N$
    \begin{equation}
        \label{eq:cluster_size_strong_decay}
        \P_{\beta}
        \left(
            k 
            < 
            \# K 
            < 
            \infty
        \right)
        \le
        \exp
        \left(
            -
            \beta k ^{(d-1)/d}
            /
            A
        \right).
    \end{equation}      
\end{theorem}

When $d \geq 1$, $\alpha \in(1, 2)$, and $\beta > \beta_c$, the decay of the distribution of finite clusters is one of the main results of our companion paper \cite{moreno-alonso_supercritical_2025}, stated below as \Cref{thm:cluster_size_weak_decay}.
\Cref{thm:cluster_size_strong_decay} is an extension and simplification of \cite[Theorem 1.2]{jorritsma_cluster-size_2024}, which proves the corresponding result for \LRP{} on $\Z^d$.
The proof of \cite{jorritsma_cluster-size_2024} proceeds by a finite-size analysis of the second largest component in a finite box, before extending this analysis to the whole infinite graph. 
This finite-size analysis and some other key arguments in the proof rely on geometric arguments specific to $\Z^d$ and its embedding into $\R^d$, in a way which does not obviously generalise to transitive graphs of polynomial growth.
We avoid the finite-size analysis and provide a simpler proof of the decay of the distribution of finite clusters by improving on some of the combinatorial bounds in \cite{jorritsma_cluster-size_2024}.
To replace the geometric arguments specific to $\Z^d$, we use the $d$-dimensional isoperimetric inequality and a Varopoulos-type isoperimetric inequality based on a mass-transport argument (\Cref{lem:want_to_show}) as well as Tim\'ar's \cite{timar_cutsets_2007,timar_boundary-connectivity_2013} notion of coarse connectivity.
Our proof also generalises \cite{jorritsma_cluster-size_2024} to a broad class of transitive kernels $J$.
Finally, the proof in \cite{jorritsma_cluster-size_2024} uses that in $\Z^d$ the sphere of radius $r$, denoted by $S(r)$, satisfies $\# S(r) \asymp r^{d-1}$.
In the setting of transitive graphs of polynomial growth however, it remains an open problem to prove the analogous result (see Section \ref{subsubsec:spheres_in_graphs}).
When dealing with the upper bound in \Cref{thm:cluster_size_strong_decay}, we are able to use the observation that ``spheres partition balls'' to avoid the use of sphere volume estimates and maintain sharp upper bounds.
For a lower bound on the distribution of finite clusters however, it seems unavoidable to use sphere volume estimates.
Conjecturally on sharp volume bounds for spheres, a renormalisation argument gives the following matching lower bound to \Cref{thm:cluster_size_strong_decay}.

\begin{theorem}
    \label[theorem]{thm:cluster_lower_bound}
    Let $G$ be a transitive graph of polynomial growth with $d \geq 2$, and suppose that $J : V \times V \to \R_+$ is a transitive kernel with $J(x,y) = \Theta(d_G(x,y)^{-d \alpha})$ with $\alpha > 1$.
    Let $\beta > \beta_c$ be sufficiently large.
    Let $\delta > 0$ and $c>0$ be such that $\# S(r) \leq c r^{d - \delta}$ for all $r \in \N$.
    Then there exists $A = A(\beta,G,J) > 0$ such that for all $k \in \N$
    \begin{equation}
    \label{eq:cluster-lower-bound-alt}
        \exp
        \left(
            - 
            A  
            k^{
                \max(
                    2-\alpha, 
                    (d-1)/d
                ) 
                +
                (1-\delta)/d
            } 
            (\log r)^{
                \charf
                \left(
                    \alpha=1+1/d
                \right)
            }
        \right)
        \leq 
        \P_{\beta}
        \left(
            k 
            < 
            \# K 
            < 
            \infty
        \right).
    \end{equation}
\end{theorem} 

    Our next result proves that upper bounds on the distribution of finite clusters imply an anchored isoperimetric inequality for the infinite cluster of long-range percolation.
    We say that a locally finite connected graph $G = (V,E)$ satisfies an \textbf{anchored $d$-dimensional isoperimetric inequality} if for some (and hence every) vertex $v \in V$ there exists $c = c(v) > 0$ such that 
\begin{equation}
    \label{eq:anchored_isop_ineq}
    \# 
    \partial_{E} W 
    \geq 
    c
    \left(
        \# 
        W
    \right)^{(d-1)/d}
\end{equation}
for every finite connected set of vertices $W \subset V$ that contains $v$, where $\partial_E W$ denotes the edge boundary of $W$ in $G$.
The \textbf{anchored isoperimetric dimension} of $G$ is defined to be the supremal value of $d$ for which $G$ satisfies an anchored $d$-dimensional isoperimetric inequality. Note that $G$ may be a random graph. 

\begin{theorem}
    \label[theorem]{thm:cluster_tail_implies_iso}
    Let $G$ be a transitive graph, and suppose that $J: V \times V \to \R_+$ is a transitive integrable kernel.
    Let $\beta>\beta_c$. 
    If there exist $\zeta$ and $c > 0$ such that $\P_{\beta} \left(k < \# \cluster <\infty \right) \leq \exp \left(-c k^{\zeta}\right)$ for every $k \in \N$, then every infinite cluster has anchored isoperimetric dimension at least $1/(1-\zeta)$ almost surely. 
\end{theorem}

This extends arguments of Pete and Hutchcroft \cite[Theorem 4.1]{pete_note_2008} \cite[Theorem 1 (iii)-(iv)]{hutchcroft_transience_2023} to the setting of \LRP{}. 
Note that the statement holds for arbitrary transitive graphs.
It follows from \cite{gandolfi_uniqueness_1992} that supercritical \LRP{} on a transitive graph of polynomial
growth has a unique infinite cluster, which we denote by $\cluster_{\infty}$.
Together with \Cref{thm:cluster_size_strong_decay,thm:cluster_size_weak_decay} below, we obtain the following immediate corollary.
\begin{corollary}
    \label[corollary]{cor:anchored_isop_dimension}
    Let $G$ be a transitive graph of polynomial growth with $d \geq 1$, and suppose that $J: V \times V \to \R_+$ is a transitive kernel with $J(x,y) = O(d_G(x,y)^{-d\alpha})$ with $\alpha>1$. 
    Let $\beta>\beta_c$. 
    If $d \geq 2$, $\alpha> 1+1/d$, and $\beta$ is sufficiently large, or if $\alpha\in (1,2)$, then $\cluster_{\infty}$ has anchored isoperimetric dimension at least $\max(d, 1/(\alpha-1))$ almost surely. 
 \end{corollary}
 
\Cref{thm:cluster_size_strong_decay,thm:cluster_size_weak_decay} prove the required upper bound on the tail of finite clusters in \Cref{thm:cluster_tail_implies_iso} with exponent $\zeta=(d-1)/d$, $\zeta=2-\alpha$, respectively. 
 So, the anchored isoperimetric dimension of the infinite component can be set to $\max(d, 1/(\alpha-1))$ in both cases. If $\alpha<1+1/d$, then LRP experiences a dimension increase. In dimension $1$, the model is only supercritical if $\alpha \in (1,2]$, and then  the anchored dimension is $1/(\alpha-1)>1$.  Here we leave the case $\alpha=2$ open.

The strategy to proving \Cref{thm:arm_exponent_1,thm:cluster_size_strong_decay} can also be used to bound the probability that a finite connected set $S$ connects to infinity.
Given a set $S \subseteq V$, we write $\{S \leftrightarrow \infty\}$ for the event that $S$ is connected to infinity by an open path, and we write $\{S \not\leftrightarrow \infty\}$ for the complement of this event.

\begin{theorem}
    \label[theorem]{thm:s_connects_infty}
    Let $G$ be a transitive graph of polynomial growth with $d \geq 2$, and suppose that $J: V \times V \to \R_+$ is a transitive kernel with $J(x,y) = O(d_G(x,y)^{-d\alpha})$ with $\alpha > 1 + 1/d$. 
    Let $\beta > \beta_c$ be sufficiently large.
    Then there exists $A=A(G,J) > 0$ such that for every finite connected set $S \subset V$
    \begin{equation}\label{eq:S-notconnecting}
        \P_{\beta}
        \left(
            S 
            \not\leftrightarrow
            \infty
        \right) 
        \leq 
        \exp
        \left(
            - 
            \beta(\# 
            S)^{(d-1)/d}
            /
            A
        \right).
    \end{equation}
\end{theorem}

We prove the above theorem in a more general setting where $S$ is any finite set with not too many components: this is made precise in \Cref{thm:s_few_components}.
In the setting of \BNNP{}, \cite{hutchcroft_transience_2023} proves that statements analogous to \Cref{thm:cluster_size_strong_decay,thm:s_connects_infty} and \Cref{cor:anchored_isop_dimension} are equivalent characterisations of the isoperimetric dimension of percolation clusters. 
 
When $d = 1$, the isoperimetric arguments in the proof of the truncated one-arm event in \Cref{thm:arm_exponent_1} break down and the proof does not apply. 
A simple argument gives the following weaker result with polylogarithmic error terms,  valid for all $d \geq 1$, $\alpha > 1$, and $\beta > \beta_c$,  provided that we know a strong enough decay on the distribution of finite clusters.

\begin{theorem}
    \label[theorem]{thm:lazy_one_arm}
    Let $G$ be a transitive graph of polynomial growth with $d \geq 1$, and suppose that $J : V \times V \to \R_+$ is a transitive kernel with $J(x,y) = \Theta(d_G(x,y)^{-d \alpha})$ with $\alpha > 1$. 
    Let $\beta > \beta_c$, and suppose that there exist $\zeta > 0$ and $c > 0$ such that $\P_{\beta} (k < \# \cluster < \infty) \leq \exp(-ck^{\zeta})$.
    Then there exist $C = C(\beta, J, G) > 0$ such that for all $r \in \N$
    \begin{equation}
        \P_{\beta}
        \left(
            o 
            \leftrightarrow 
            B(r)^c,
            \# \cluster
            < \infty
        \right) 
        \leq
        C
        \mathrm{polylog}(r)
        r^{d(1-\alpha)}.
    \end{equation}
\end{theorem}

\subsection{Discussion}
\subsubsection{The driving exponent in the decay of finite clusters}
\label{subsubsec:driving_exponent}
In our companion paper \cite{moreno-alonso_supercritical_2025} we prove the following result.

\begin{theorem}[Cluster-size decay in the weak-decay regime \cite{moreno-alonso_supercritical_2025}]
\label[theorem]{thm:cluster_size_weak_decay}
    Let $G$ be a transitive graph of polynomial growth with $d \geq 1$, and suppose that $J : V \times V \to \R_+$ is a transitive kernel with $J(x,y) = \Theta(d_G(x,y)^{-d \alpha})$ with $\alpha \in (1, 2)$. 
    Let $\beta > \beta_c$.
    Then there exists $A = A(\beta,G, J) > 0$ such that for all $k \in \N$
    \begin{equation}
        \label{eq:cluster_size_decay}
        \P_{\beta}
        (
            k
            \leq
            \# \cluster
            < 
            \infty
        )
        \leq
        \exp
        \left(
            - k^{2 - \alpha} 
            / 
            A
        \right).
    \end{equation}
\end{theorem}

Whereas the proof of \Cref{thm:cluster_size_strong_decay} is combinatorial, the proof of \Cref{thm:cluster_size_weak_decay} proceeds by a renormalisation argument specific to the regime where $\alpha \in (1,2)$.
In particular, \Cref{thm:cluster_size_weak_decay} holds for all $\beta > \beta_c$.
When $d=1$ it is a consequence of \cite{schulman_long_1983,newman_one-dimensional_1986} that there is a non-trivial phase transition if and only if $ \alpha \in (1, 2]$, and in this setting \Cref{thm:cluster_size_weak_decay} gives that the decay of the distribution of finite clusters is driven by the long-range exponent $2-\alpha$ for all $\beta > \beta_c$.
When $d \geq 2$ however, there is a transition in the exponent driving the decay of the distribution of finite clusters at $\alpha = 1 + 1/d$: the driving exponent goes from the long-range exponent $2-\alpha$ in \Cref{thm:cluster_size_strong_decay} to the surface-tension exponent $(d-1)/d$ in \Cref{thm:cluster_size_weak_decay}.
The \textbf{expected total edge-length}
\begin{equation}
    \expleng{o}
    =
    \sum_{y \in G \setminus \{o\}}
    d_G(o,y)
    \P_{\beta}
    \left(
        o 
        \sim
        y
    \right)
\end{equation}
is finite if and only if $\alpha > 1 + 1/d$ (see \Cref{lem:finite_exp_leng}), and hence the value $\alpha = 1+1/d$ can be thought of as a transition to ``short-range'' behaviour of the model, as indeed the surface-tension exponent $(d-1)/d$ is the same as that for \BNNP{}, due to \cite{kunz_essential_1978,kesten_probability_1990,contreras_supercritical_2024}.
Further, \Cref{cor:anchored_isop_dimension} gives that for $\alpha<1+1/d$ the infinite component $\cluster_{\infty}$ has anchored isoperimetric dimension at least $1/(\alpha-1)>d$, so the parameter $\alpha$ may be thought of as continuously varying the dimension of the model.
The transition in the driving exponent as $\alpha$ increases above $1+1/d$ may be attributed to the increasing importance of isoperimetry and short edges. 
See a related discussion in \Cref{rem:crucial_alpha_assumption}.

\begin{figure}
    \centering
    \begin{tikzpicture}[scale=2,>=latex,every node/.style={font=\small}]
    % Approximate widths in cm for text in regions with text
    \def\wOne{7}    % region (1, 1+1/d)
    \def\wTwo{8.3}    % region (1+1/d, 2)
    \def\wThree{8.3}  % green region width
    
    % Width for (0,1) region (small)
    \def\wZeroOne{1.0} 
    
    % Convert cm to tikz units (1cm ~ 2 units at scale=2)
    \pgfmathsetmacro{\widthZeroOne}{\wZeroOne/2}  
    \pgfmathsetmacro{\widthOne}{\wOne/2}          
    \pgfmathsetmacro{\widthTwo}{\wTwo/2}          
    \pgfmathsetmacro{\widthThree}{\wThree/2}      
    
    % Padding inside each region
    \def\pad{0.15}
    
    % Set alpha boundaries:
    \coordinate (a0) at (0,0); % 0
    \coordinate (a1) at ($(a0) + (\widthZeroOne,0)$); % 1 (small region, no text)
    \coordinate (a2) at ($(a1) + (\widthOne + 2*\pad,0)$); % 1 + 1/d
    \coordinate (a3) at ($(a2) + (\widthTwo + 2*\pad,0)$); % 2
    \coordinate (B)  at ($(a3) + (\widthThree + 2*\pad,0)$); % infinity
    
    % Shade heights
    \def\shadeHeight{1.0}
    
    % Shade regions (no shade on (0,1))
    \fill[blue!10]  (a1) rectangle ($(a2)+(0,\shadeHeight)$);
    \fill[red!10]   (a2) rectangle ($(a3)+(0,\shadeHeight)$);
    \fill[green!10] (a3) rectangle ($(B)+(0,\shadeHeight)$);
    
    % Axis
    \draw[->] (a0) -- (B) node[right] {$\alpha$};
    
    % Tick marks and labels
    \draw (a0) -- ++(0,-0.08) node[below] {$0$};
    \draw (a1) -- ++(0,-0.08) node[below] {$1$};
    \draw (a2) -- ++(0,-0.08) node[below] {$1+1/d$};
    \draw (a3) -- ++(0,-0.08) node[below] {$2$};
    
    % Dashed lines
    \draw[dashed] (a1) -- +(0,\shadeHeight);
    \draw[dashed] (a2) -- +(0,\shadeHeight);
    \draw[dashed] (a3) -- +(0,\shadeHeight);
    
    % --- Texts for regions with text, centered inside each ---
    % Region (1, 1+1/d)
    \node[align=center,above,yshift=0.5cm,text width=\wOne cm] 
    at ($(a1)!0.5!(a2)$) 
    {$\displaystyle 
    \mathbb{P}_{\beta}(k \leq K < \infty) 
    \asymp \exp\bigl(- k^{2-\alpha}\bigr) 
    \quad \forall\, \beta > \beta_c$};
    
    % Region (1+1/d, 2) two stacked lines
    \node[align=center,above,yshift=0.75cm,text width=\wTwo cm] 
    at ($(a2)!0.5!(a3)$) 
    {$\displaystyle 
    \mathbb{P}_{\beta}(k \leq K < \infty) 
    \leq \exp\bigl(- k^{2-\alpha}\bigr) 
    \quad \forall\, \beta > \beta_c$};
    
    \node[align=center,above,yshift=0.15cm,text width=\wTwo cm] 
    at ($(a2)!0.5!(a3)$) 
    {$\displaystyle 
    \mathbb{P}_{\beta}(k \leq K < \infty) 
    \asymp \exp\bigl(- k^{\frac{d-1}{d}}\bigr) 
    \quad \forall\, \beta > \beta_c \text{ suff. large}$};
    
    % Region (2, infinity) — single line expression, no scaling
    \node[align=center,above,yshift=0.6cm,text width=\wThree cm] 
    at ($(a3)!0.5!(B)$) 
    {$\displaystyle 
    \mathbb{P}_{\beta}(k \leq K < \infty) 
    \asymp \exp\bigl(- k^{\frac{d-1}{d}}\bigr) 
    \quad \forall\, \beta > \beta_c \text{ suff. large}$};
\end{tikzpicture}
    \caption{The results of \Cref{thm:cluster_size_strong_decay,thm:cluster_lower_bound,thm:cluster_size_weak_decay} as a function of the long-range parameter $\alpha$, conjectural on sharp bounds for the volume of spheres in transitive graphs of polynomial growth. We prove both upper and lower bounds for all $\beta>\beta_c$ when $\alpha\in(1,2)$ in \cite{moreno-alonso_supercritical_2025}. 
    Here, $d \geq 2$.}
    \label{fig:alpha_regions}
\end{figure}
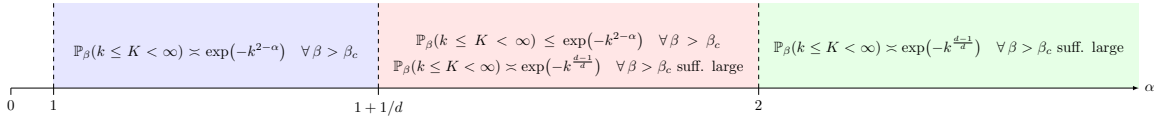

\subsubsection{Related results}
In the setting of supercritical \BNNP{} in $\Z^d$, \cite{kesten_critical_1980} proved the exponential decay of the truncated one-arm event when $d = 2$, and \cite{chayes_correlation_1989} proved the matching result for $d \geq 3$. 
The surface-tension order decay of the distribution of finite clusters is proved by \cite{kunz_essential_1978} for $p$ sufficiently close to $1$, and \cite{kesten_probability_1990} extend this to all $p > p_c$ (conditional on the then-conjectural Grimmett-Marstrand theorem).
\cite{hermon_supercritical_2021} prove the decay of the truncated one-arm event and the distribution of finite clusters for supercritical \BNNP{} on transitive non-amenable graphs.
Most relevant to this paper, \cite{contreras_supercritical_2024} prove the decay of the truncated one-arm event and the distribution of finite clusters for supercritical \BNNP{} on transitive graphs of polynomial growth, in particular giving a new proof of the Grimmett-Marstrand theorem.
In the study of scaling relations for percolation (see \cite[Chapter 9]{grimmett_percolation_1999}), the exponent in the decay of the one-arm event at criticality is sometimes called the one-arm exponent.
The one-arm exponent continues to be an active area of investigation in percolation, and we refer to the recent papers \cite{engelenburg_one-arm_2025-1,engelenburg_one-arm_2025} for an overview of the literature on the one-arm exponent.

For \LRP{}, the one-arm exponent for the mean-field critical regime is established in \cite{hulshof_one-arm_2015}.
In the subcritical regime $\beta < \beta_c$, \cite{duminil-copin_new_2016} prove the exponential decay of the one-arm event for finite-range percolation on arbitrary transitive graphs.
For $\beta \leq \beta_c$, \cite{hutchcroft_power-law_2021} proves power-law upper bounds for the distribution of finite clusters for \LRP{} on transitive unimodular graphs.
These bounds are improved (but for $\Z^d$) in \cite{hutchcroft_sharp_2022}, and lower bounds are proved in \cite{baumler_isoperimetric_2024}.
Several properties of supercritical \LRP{} on $\Z^d$, including the distribution of finite clusters, are considered in \cite{jorritsma_cluster-size_2024,jorritsma_cluster-size_2025,jorritsma_large_2025}.
Critical \LRP{} has recently been considered in \cite{hutchcroft_critical_2025,hutchcroft_critical_2025-2,hutchcroft_critical_2025-3}.

\subsubsection{Supercritical sharpness}
Proving sharp bounds on the distribution of finite clusters and the truncated one-arm event for the entire supercritical regime is sometimes called \textit{supercritical sharpness}.
Supercritical sharpness is proved for \BNNP{} on transitive non-amenable graphs \cite{hermon_supercritical_2021}, transitive graphs of polynomial growth \cite{contreras_supercritical_2024}, and Voronoi percolation \cite{dembin_supercritical_2025}.
While for $\alpha \in (1, 1+1/d)$ we have sharp bounds on the distribution of finite clusters for all $\beta > \beta_c$ (\Cref{thm:cluster_size_weak_decay}), it is inherent to the combinatorial techniques of this paper that we can only prove \Cref{thm:arm_exponent_1,thm:cluster_size_strong_decay,thm:cluster_lower_bound} for $\beta > \beta_c$ sufficiently large.
In the setting of transitive graphs of polynomial growth, it is not obvious how to extend these results to all $\beta > \beta_c$.
 
\subsubsection{Sphere volumes in transitive graphs of polynomial growth}
\label{subsubsec:spheres_in_graphs}
The volume of spheres appears naturally when working with \LRP{}. 
Perhaps the most natural instance is when verifying that a transitive kernel $J$ satisfying $J(x,y) = O(d_G(x,y)^{-d \alpha})$ with $\alpha > 1$ is integrable, namely that $\sum_{y \in G \setminus \{x\}} J(x,y) < \infty$ for all $x \in G$.
The integrability of the kernel implies that the open subgraph sampled by $\P_{\beta}$ almost surely has finite expected degree, since for all $x \in G$ we have
\begin{equation}        
\label{eq:expdeg_first_bound}
    \expdeg{x}
    =
    \sum_{y \in G \setminus \{x\}}
    \left(
        1 
        - 
        \exp
        \left(
            -
            \beta
            J(x,y)
        \right) 
    \right)
    \leq
    \beta
    \sum_{y \in G \setminus \{x\}}
    J(x,y)
    < 
    \infty.
\end{equation}
We can write
\begin{equation}
    \label{eq:kernel_integrability_spheres}
    \sum_{y \in G \setminus \{x\}} 
    J(x,y)
    \le 
    \sum_{r = 1}^{R_J - 1}
    \sum_{y \in S(x,r)}
    J(x,y)
    +
    C_J
    \sum_{r = R_J}^{\infty}
    \# 
    S(x,r)
    r^{-d \alpha}.
\end{equation}
When $G = \Z^d$, it is a standard fact that $\# S(r) \asymp r^{d-1}$ and it is immediate that $J$ is integrable for $\alpha > 1$.
In the setting of transitive graphs of polynomial growth however, it remains an open problem to find sharp bounds for the volume of spheres.
For any locally finite transitive graph it holds that $\# B(r) \leq 2r \# S(r)$, and in the setting transitive graphs of polynomial growth the lower bound $\# S(r) \gtrsim r^{d-1}$ follows immediately from the volume bounds for balls.
It follows from \cite{colding_liouville_1998,tessera_volume_2007} that in the general setting of doubling metric spaces, there exists $\delta > 0$ depending on the \textit{doubling constant} such that $\# S(r) \lesssim r^{d-\delta}$ (see \cite[Proposition 6.10]{spanos_spread-out_2024} for a precise statement in the setting of transitive graphs).
When $G$ is the Cayley graph of an $s$-step virtually nilpotent group, it follows from \cite{breuillard_rate_2013,gianella_asymptotics_2017,bodart_intermediate_2025} that $S(r) \lesssim r^{d-\delta_s}$ with $\delta_s = 1$ if $s = 1,2$ and $\delta_s = 1/s$ if $s \geq 3$.
It is conjectured \cite[Conjecture 10]{breuillard_rate_2013} that, at least in the setting of Cayley graphs of nilpotent groups, the volume of spheres satisfies $\#S(r) \lesssim r^{d-1}$, so that $\# S(r) \asymp r^{d-1}$. 

Although the estimates in \cite{colding_liouville_1998,tessera_volume_2007} can be used in \eqref{eq:kernel_integrability_spheres} to conclude that the kernel is integrable for $\alpha > 2 - \delta$, this does not give the true threshold for integrability, which is $\alpha > 1$.
More generally, the volume of spheres appears repeatedly in our arguments in the form of the expected degree $\expdeg{x}$ and the expected total edge-length $\expleng{x}$, and if we were to use sphere volume estimates throughout, the error term $1-\delta$ would find its way into the statements of \Cref{thm:arm_exponent_1,thm:cluster_size_strong_decay}.
In Section \ref{subsec:integrability_expected_degree} we use the elementary observation that ``spheres partition balls'' to establish the integrability of the kernel for $\alpha > 1$, and in particular that $\expdeg{x} < \infty$ and $\expleng{x} < \infty$ for $\alpha > 1$ and $\alpha > 1 + 1/d$ respectively, without the use of sphere volume bounds.

For the lower bound for the distribution of finite clusters in \Cref{thm:cluster_lower_bound}, it seems inevitable to use estimates for the upper bound on the volume of spheres (see \Cref{rem:sphere_inevitability}).
This is the only place where we require the results of \cite{colding_liouville_1998,tessera_volume_2007,breuillard_rate_2013,gianella_asymptotics_2017,bodart_intermediate_2025}.

\subsection{About the proofs}
\subsubsection{The truncated one-arm from cluster-size decay}
Given \Cref{thm:cluster_size_strong_decay,thm:cluster_size_weak_decay} on the decay of the distribution of finite clusters, we can choose $k(r)$ as a polylogarithmic function of $r$ in \eqref{eq:one_arm_decomposition} to obtain the required decay. 

\begin{lemma}
    \label[lemma]{lem:one_arm_decomp}
    Let $G$ be a transitive graph of polynomial growth with $d \geq 1$, let $J : G \times G \to \R_+$ be transitive and suppose that $J(x,y) = \Theta(d_G(x,y)^{-d\alpha})$ with $\alpha > 1$. 
    Let $\beta > \beta_c$.
    If $d \geq 2$, $\alpha > 1 + 1/d$, and $\beta$ is sufficiently large, or if $\alpha \in (1,2)$, then there exists $k(r) = \mathrm{polylog}(r) > 0$ such that for all $r \in \N$
    \begin{equation}
        \label{lem:one_arm_decomposition_lemma}
        \P_{\beta}
        \left(
            o 
            \leftrightarrow 
            B(r)^c,
            \# K < \infty
        \right)
        = 
        o
        \left(
            r^{d(1-\alpha)} 
        \right) 
        + 
        \P_{\beta}
        \left(
            o 
            \leftrightarrow B(r)^c, 
            \# \cluster
            \leq 
            k(r)
        \right).
    \end{equation}
\end{lemma}

\begin{proof}
    Let $k_1(r) = \left(2 A_1 d(1-\alpha) \log(r) \right)^{(d-1)/d}$ and $k_2(r) = \left( 2 A_2 d(1-\alpha) \log(r)\right)^{d(2-\alpha)}$, where $A_1$ and $A_2$ are as in \Cref{thm:cluster_size_strong_decay,thm:cluster_size_weak_decay} respectively.
    If $d \geq 2$, $\alpha > 1 + 1/d$, and $\beta$ is sufficiently large, then by \Cref{thm:cluster_size_strong_decay} we have
    $
        \P_{\beta}
        \left(
            k_1(r)
            \leq 
            \#
            \cluster
            < 
            \infty
        \right)
        \leq
        r^{
            -2d(1-\alpha)
        }.
    $ 
    If $\alpha\in (1,2)$, then by \Cref{thm:cluster_size_weak_decay} we have
    $
        \P_{\beta}
        \left(
            k_2(r)
            \leq 
            \#
            \cluster
            < 
            \infty
        \right)
        \leq
        r^{
            -2d(1-\alpha)
        }.
    $
    Both upper bounds are $o(r^{d(1-\alpha)})$, and together with the union bound in \eqref{eq:one_arm_decomposition} the statement follows.
\end{proof}

The goal of the remainder of the paper is to bound $\P_{\beta}\left( 0 \leftrightarrow B(r)^c, \# \cluster \leq k(r) \right)$, and to prove the bounds on the distribution of finite clusters in \Cref{thm:cluster_size_strong_decay,thm:cluster_lower_bound}.

\subsubsection{Outline of the proof}
For an upper bound on the probability of an event involving $\{\# \cluster < \infty \}$, such as the cluster-size and truncated one-arm events $\{k \leq \# \cluster < \infty\}$ and $\{o \leftrightarrow B(r)^c, \# \cluster < \infty\}$, we need to specify that certain edges in the cluster of the origin are open or closed.
In \BNNP{}, this is typically done using the theory of \textit{lattice animals} \cite[Chapter 4.2]{grimmett_percolation_1999}.
A lattice animal is a finite connected subgraph of $G$ containing the origin, and in the setting of \BNNP{} the finite cluster of the origin is a random lattice animal.
If we write $\P_{p}$ for the law of Bernoulli nearest-neighbour (bond) percolation on $G$ with parameter $0 \leq p \leq 1$, then
\begin{equation}
    \label{eq:bnnp_lattice_animals}
    \P_{p}
    \left(
        \# \cluster
        < \infty
    \right)
    =
    \sum_{\substack{
        A \subset G, \# A < \infty
        \\
        A \ni o \text{ connected}
    }}
    p^{\# E(A)}
    (1-p)^{\# \partial A}.
\end{equation}  

\begin{figure}
  \centering
  \begin{subfigure}{0.32\textwidth}
    \centering
    \begin{tikzpicture}[scale=0.6,
    vertex/.style={circle, fill=black, inner sep=0.8pt},
    gridline/.style={gray!30, very thin},
    longedge/.style={line width=1.2pt, gray!60, bend left=15},
    edgeA/.style={line width=1.2pt, red!70!black, opacity=0.85},
    vertA/.style={circle, fill=red!60, inner sep=2pt},
    edgeB/.style={line width=1.2pt, orange!80!black, opacity=0.85},
    vertB/.style={circle, fill=orange!70, inner sep=2pt},
    edgeC/.style={line width=1.2pt, green!70!black, opacity=0.85},
    vertC/.style={circle, fill=green!60, inner sep=2pt},
    edgeD/.style={line width=1.2pt, purple!70!black, opacity=0.85},
    vertD/.style={circle, fill=purple!60, inner sep=2pt},
    edgeE/.style={line width=1.2pt, cyan!70!black, opacity=0.85},
    vertE/.style={circle, fill=cyan!60, inner sep=2pt},
    edgeF/.style={line width=1.2pt, pink!80!black, opacity=0.9},
    vertF/.style={circle, fill=pink!80!red, inner sep=2pt},
    edgeG/.style={line width=1.2pt, teal!80!black, opacity=0.9},
    vertG/.style={circle, fill=teal!70, inner sep=2pt},
    edgeH/.style={line width=1.2pt, magenta!80!black, opacity=0.9},
    vertH/.style={circle, fill=magenta!70, inner sep=2pt}
]

% ----------------------------------------------------------------
% Draw 13x13 grid
\foreach \x in {-6,...,6} \draw[gridline] (\x,-6.5) -- (\x,6.5);
\foreach \y in {-6,...,6} \draw[gridline] (-6.5,\y) -- (6.5,\y);
\foreach \x in {-6,...,6}{
  \foreach \y in {-6,...,6}{
    \node[vertex] at (\x,\y) {};
  }
}

% ================================================================
% Animal A
% ================================================================
\draw[edgeA]
  (-1,-1)--(-1,0)--(-1,1)--(0,1)--(1,1)--(1,2)--(0,2)--(0,1)--(0,0)--(-1,0);
\foreach \p in {
  (-1,-1),(-1,0),(-1,1),(0,0),(0,1),(0,2),(1,1),(1,2)
}{
  \node[vertA] at \p {};
}

% ================================================================
% Animal B
% ================================================================
\draw[edgeB]
  (-5,2)--(-5,5)--(-2,5)--(-2,2)--(-5,2)
  (-2,2)--(-3,2)
  (-2,2)--(-2,3)
  (-3,3)--(-2,3)
  (-3,3)--(-3,2);
\foreach \p in {
  (-5,2),(-5,3),(-5,4),(-5,5),
  (-4,5),(-3,5),(-2,5),
  (-2,4),(-2,3),(-2,2),
  (-3,2),(-3,3),(-4,2)
}{
  \node[vertB] at \p {};
}

% ================================================================
% Animal C
% ================================================================
\draw[edgeC]
  (3,5)--(5,5)--(5,3)--(4,3)--(4,4)--(3,4)--cycle;
\foreach \p in {(3,5),(4,5),(5,5),(5,4),(5,3),(4,3),(4,4),(3,4)}{
  \node[vertC] at \p {};
}

% ================================================================
% Animal D
% ================================================================
\draw[edgeD]
  (-5,-1)--(-5,-2)--(-3,-2)--(-3,-3)--(-3,-4)--(-3,-5)--(-4,-5)--(-5,-5)--(-5,-4)--(-5,-3)--(-5,-2)
  (-4,-4)--(-5,-4)
  (-4,-4)--(-3,-4)
  (-4,-4)--(-4,-5);
\foreach \p in {
  (-5,-1),
  (-5,-2),(-4,-2),(-3,-2),
  (-3,-3),(-3,-4),(-3,-5),
  (-4,-4),(-4,-5),
  (-5,-5),(-5,-4),(-5,-3)
}{
  \node[vertD] at \p {};
}

% ================================================================
% Animal E
% ================================================================
\draw[edgeE]
  (2,-1)--(1,-1)--(1,-2)--(1,-3)--(0,-3)--(0,-4)--(1,-4)--(1,-5)--(2,-5)--(3,-5)--(4,-5)--(5,-5)--(5,-4)--(5,-3)--(5,-2)--(5,-1)--(5,0)--(2,0)--(2,-1)
  (2,-4)--(1,-4)
  (2,-4)--(2,-5)
  (4,-4)--(4,-5)
  (4,-4)--(5,-4);
\foreach \p in {
  (2,0),(3,0),(4,0),(5,0),
  (5,-1),(5,-2),(5,-3),(5,-4),(5,-5),
  (4,-5),(4,-4),(3,-5),(2,-5),(2,-4),
  (1,-5),(1,-4),(1,-3),(1,-2),(1,-1),
  (0,-3),(0,-4),
  (2,-1)
}{
  \node[vertE] at \p {};
}

% ================================================================
% Animal F
% ================================================================
\draw[edgeF] (3,-2)--(3,-3);
\foreach \p in {(3,-2),(3,-3)}{
  \node[vertF] at \p {};
}

% ================================================================
% Animal G
% ================================================================
\draw[edgeG] (-4,0)--(-3,0);
\foreach \p in {(-4,0),(-3,0)}{
  \node[vertG] at \p {};
}

% ================================================================
% Animal H
% ================================================================
\draw[edgeH] (0,4)--(1,4)--(1,5);
\foreach \p in {(0,4),(1,4),(1,5)}{
  \node[vertH] at \p {};
}

% ================================================================
% Long curved edges
% ================================================================
\draw[longedge, bend left=20] (-1,1) to (3,4);      % A → C
\draw[longedge, bend right=25] (-3,5) to (-4,-2);   % B → D
\draw[longedge, bend right=20] (-3,0) to (1,-4);    % G → E
\draw[longedge, bend left=25] (1,5) to (4,0);       % H → E
\draw[longedge, bend left=20] (-2,5) to (0,2);      % B → A
\draw[longedge, bend right=20] (5,0) to (5,4);      % E → C
\draw[longedge, bend right=20] (-2,2) to (-1,1);      
\draw[longedge, bend left=20] (1,1) to (2,0);   
\draw[longedge, bend left=20] (0,-4) to (-3,-5);
\draw[longedge, bend left=20] (-5,-1) to (-5,3);

% ================================================================
% Long edges within blocks
% ================================================================
\draw[longedge, bend left=15] (-1,-1) to (1,2);      % inside A
\draw[longedge, bend left=20] (-5,5) to (-3,2);      % inside B
\draw[longedge, bend right=20] (3,5) to (5,3);       % inside C
\draw[longedge, bend left=25] (-5,-1) to (-3,-4);    % inside D
\draw[longedge, bend right=25] (1,-1) to (4,-4);     % inside E
\draw[longedge, bend left=30] (3,-3) to (2,-4);      % F → E

\end{tikzpicture}
    \caption{}
    \label{fig:lrp_configuration}
  \end{subfigure}
  \hfill
  \begin{subfigure}{0.32\textwidth}
    \centering
    \begin{tikzpicture}[scale=0.6,
    vertex/.style={circle, fill=black, inner sep=0.8pt},
    gridline/.style={gray!30, very thin},
    vertA/.style={circle, fill=red!60, inner sep=2pt},
    vertB/.style={circle, fill=orange!70, inner sep=2pt},
    vertC/.style={circle, fill=green!60, inner sep=2pt},
    vertD/.style={circle, fill=purple!60, inner sep=2pt},
    vertE/.style={circle, fill=cyan!60, inner sep=2pt},
    vertG/.style={circle, fill=teal!70, inner sep=2pt},
    vertH/.style={circle, fill=magenta!70, inner sep=2pt}
]

% ----------------------------------------------------------------
% Draw 13x13 grid
\foreach \x in {-6,...,6} \draw[gridline] (\x,-6.5) -- (\x,6.5);
\foreach \y in {-6,...,6} \draw[gridline] (-6.5,\y) -- (6.5,\y);
\foreach \x in {-6,...,6}{
  \foreach \y in {-6,...,6}{
    \node[vertex] at (\x,\y) {};
  }
}

% ================================================================
% Animal A
% ================================================================
\path[fill=red!30, draw=red!30, line width=6pt, opacity=0.35, line cap=round, line join=round]
  (-1,-1)--(-1,0)--(-1,1)--(0,1)--(1,1)--(1,2)--(0,2)--(0,1)--(0,0)--(-1,0)--cycle;

\foreach \p in {(-1,-1),(-1,0),(-1,1),(0,0),(0,1),(0,2),(1,1),(1,2)}{
  \node[vertA] at \p {};
}

% ================================================================
% Animal B
% ================================================================
\path[fill=orange!30, draw=orange!30, line width=6pt, opacity=0.35, line cap=round, line join=round]
  (-5,2)--(-5,5)--(-2,5)--(-2,2)--cycle;

\foreach \x in {-5,-4,-3,-2}{
  \foreach \y in {2,3,4,5}{
    \node[vertB] at (\x,\y) {};
  }
}

% ================================================================
% Animal C
% ================================================================
\path[fill=green!30, draw=green!30, line width=6pt, opacity=0.35, line cap=round, line join=round]
  (3,5)--(5,5)--(5,3)--(4,3)--(4,4)--(3,4)--cycle;

\foreach \p in {(3,4),(3,5),(4,3),(4,4),(4,5),(5,3),(5,4),(5,5)}{
  \node[vertC] at \p {};
}

% ================================================================
% Animal D
% ================================================================
\path[fill=purple!30, draw=purple!30, line width=6pt, opacity=0.35, line cap=round, line join=round]
  (-5,-1)--(-5,-2)--(-3,-2)--(-3,-5)--(-5,-5)--cycle;

\foreach \x in {-5,-4,-3}{
  \foreach \y in {-5,-4,-3,-2}{
    \node[vertD] at (\x,\y) {};
  }
}
\node[vertD] at (-5,-1) {};

% ================================================================
% Animal E
% ================================================================
\path[fill=cyan!30, draw=cyan!30, line width=6pt, opacity=0.35, line cap=round, line join=round]
  (0,-4)--(1,-4)--(1,-5)--(5,-5)--(5,0)--(2,0)--(2,-1)--(1,-1)--(1,-3)--(0,-3)--cycle;

\foreach \p in {
  (0,-4),(0,-3),
  (1,-5),(1,-4),(1,-3),(1,-2),(1,-1),
  (2,-5),(2,-4),(2,-3),(2,-2),(2,-1),(2,0),
  (3,-5),(3,-4),(3,-3),(3,-2),(3,-1),(3,0),
  (4,-5),(4,-4),(4,-3),(4,-2),(4,-1),(4,0),
  (5,-5),(5,-4),(5,-3),(5,-2),(5,-1),(5,0)
}{
  \node[vertE] at \p {};
}

% ================================================================
% Animal G
% ================================================================
\draw[teal!30, line width=6pt, opacity=0.35, line cap=round, line join=round] (-4,0)--(-3,0);
\foreach \p in {(-4,0),(-3,0)}{
  \node[vertG] at \p {};
}

% ================================================================
% Animal H
% ================================================================
\draw[magenta!30, line width=6pt, opacity=0.35, line cap=round, line join=round] (0,4)--(1,4)--(1,5);
\foreach \p in {(0,4),(1,4),(1,5)}{
  \node[vertH] at \p {};
}

\end{tikzpicture}
    \caption{}
    \label{fig:associated_decomposition}
  \end{subfigure}
  \hfill
  \begin{subfigure}{0.32\textwidth}
    \centering
    \begin{tikzpicture}[scale=0.6,
    vertex/.style={circle, fill=black, inner sep=0.8pt},
    gridline/.style={gray!30, very thin},
    vertA/.style={circle, fill=gray!70, inner sep=2pt},
    vertB/.style={circle, fill=gray!70, inner sep=2pt},
    vertC/.style={circle, fill=gray!70, inner sep=2pt},
    vertD/.style={circle, fill=gray!70, inner sep=2pt},
    vertE/.style={circle, fill=cyan!80!blue, inner sep=2pt}, % highlighted
    vertG/.style={circle, fill=gray!70, inner sep=2pt},
    vertH/.style={circle, fill=gray!70, inner sep=2pt}
]

% ----------------------------------------------------------------
% Draw 13x13 grid
\foreach \x in {-6,...,6} \draw[gridline] (\x,-6.5) -- (\x,6.5);
\foreach \y in {-6,...,6} \draw[gridline] (-6.5,\y) -- (6.5,\y);
\foreach \x in {-6,...,6}{
  \foreach \y in {-6,...,6}{
    \node[vertex] at (\x,\y) {};
  }
}

% ================================================================
% Animal A (grayed out)
% ================================================================
\path[fill=gray!60, draw=gray!60, line width=6pt, opacity=0.6, line cap=round, line join=round]
  (-1,-1)--(-1,0)--(-1,1)--(0,1)--(1,1)--(1,2)--(0,2)--(0,1)--(0,0)--(-1,0)--cycle;

\foreach \p in {(-1,-1),(-1,0),(-1,1),(0,0),(0,1),(0,2),(1,1),(1,2)}{
  \node[vertA] at \p {};
}

% ================================================================
% Animal B (grayed out)
% ================================================================
\path[fill=gray!60, draw=gray!60, line width=6pt, opacity=0.6, line cap=round, line join=round]
  (-5,2)--(-5,5)--(-2,5)--(-2,2)--cycle;

\foreach \x in {-5,-4,-3,-2}{
  \foreach \y in {2,3,4,5}{
    \node[vertB] at (\x,\y) {};
  }
}

% ================================================================
% Animal C (grayed out)
% ================================================================
\path[fill=gray!60, draw=gray!60, line width=6pt, opacity=0.6, line cap=round, line join=round]
  (3,5)--(5,5)--(5,3)--(4,3)--(4,4)--(3,4)--cycle;

\foreach \p in {(3,4),(3,5),(4,3),(4,4),(4,5),(5,3),(5,4),(5,5)}{
  \node[vertC] at \p {};
}

% ================================================================
% Animal D (grayed out)
% ================================================================
\path[fill=gray!60, draw=gray!60, line width=6pt, opacity=0.6, line cap=round, line join=round]
  (-5,-1)--(-5,-2)--(-3,-2)--(-3,-5)--(-5,-5)--cycle;

\foreach \x in {-5,-4,-3}{
  \foreach \y in {-5,-4,-3,-2}{
    \node[vertD] at (\x,\y) {};
  }
}
\node[vertD] at (-5,-1) {};

% ================================================================
% Animal E (highlighted)
% ================================================================
\path[fill=cyan!40, draw=cyan!70!blue, line width=6pt, opacity=0.8, line cap=round, line join=round]
  (0,-4)--(1,-4)--(1,-5)--(5,-5)--(5,0)--(2,0)--(2,-1)--(1,-1)--(1,-3)--(0,-3)--cycle;

\foreach \p in {
  (0,-4),(0,-3),
  (1,-5),(1,-4),(1,-3),(1,-2),(1,-1),
  (2,-5),(2,-4),(2,-3),(2,-2),(2,-1),(2,0),
  (3,-5),(3,-4),(3,-3),(3,-2),(3,-1),(3,0),
  (4,-5),(4,-4),(4,-3),(4,-2),(4,-1),(4,0),
  (5,-5),(5,-4),(5,-3),(5,-2),(5,-1),(5,0)
}{
  \node[vertE] at \p {};
}

% ================================================================
% Animal G (grayed out)
% ================================================================
\draw[gray!60, line width=6pt, opacity=0.6, line cap=round, line join=round] (-4,0)--(-3,0);
\foreach \p in {(-4,0),(-3,0)}{
  \node[vertG] at \p {};
}

% ================================================================
% Animal H (grayed out)
% ================================================================
\draw[gray!60, line width=6pt, opacity=0.6, line cap=round, line join=round] (0,4)--(1,4)--(1,5);
\foreach \p in {(0,4),(1,4),(1,5)}{
  \node[vertH] at \p {};
}

% ================================================================
% Highlight vertex (2, -1)
% ================================================================
\node[circle, fill=cyan!40, draw=black, thick, minimum size=10pt] at (1,-1) {};

% ================================================================
% Red path
% ================================================================
\draw[line width=3pt, red!90!orange, opacity=0.8, line cap=round, line join=round, ->]
  (1,-1)--(1,0)--(2,1)--(2,2)--(1,3)--(0,3)--(-1,2)--(-2,1)--(-2,0)--(-3,-1)--(-4,-1)--(-5,0);

% Mark vertices along the path
\foreach \p in {(1,-1),(1,0),(2,1),(2,2),(1,3),(0,3),(-1,2),(-2,1),(-2,0),(-3,-1),(-4,-1)}{
  \node[circle, fill=red!90!orange, inner sep=2pt, opacity=0.8] at \p {};
}

\end{tikzpicture}
    \caption{}
    \label{fig:isolation}
  \end{subfigure}

  \caption{
    Figures \ref{fig:lrp_configuration} and \ref{fig:associated_decomposition} show a sketch of a \LRP{} configuration and the associated block decomposition.
    In Figure \ref{fig:isolation} the red path is a coarsely connected path which starts from a boundary vertex and whose subsequent vertices are all guaranteed to \textit{not} be in $K$, so that all edges from the boundary vertex to path vertices must be closed.
    This will be used to bound the isolation of blocks.
  }
  \label{fig:overall}
\end{figure}
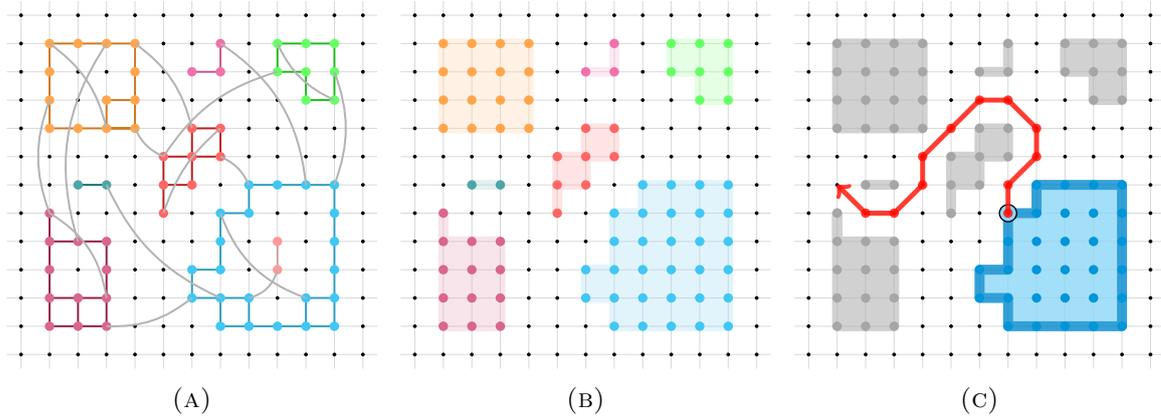

By counting lattice animals, Kunz and Souillard \cite{kunz_essential_1978} proved a bound on the distribution of finite clusters for $p$ sufficiently close to $1$, later extended to all $p > p_c$ by Kesten and Zhang \cite{kesten_probability_1990}, conditionally on the then-conjectural Grimmett-Marstrand theorem. 
In \LRP{} however, the finite cluster of the origin $K$ is not a lattice animal: it is a finite connected subgraph of the \textit{complete} graph on $G$ which contains the origin. 
Counting the possible configurations of $K$ is combinatorially more complex, and we approach this using the \textit{block decomposition} of $K$, due to \cite{jorritsma_cluster-size_2024}.
A block is a set of nearest-neighbour vertices in $G$ which contains all vertices for which it is a cut-set.
In  \Cref{lem:lattice_animal} we show that $K$ is contained in a unique set of mutually disjoint blocks such that the distance between two blocks is at least two.
This is the block decomposition of $K$.
The block graph of $K$ is the resulting graph when blocks are identified as vertices and edges are given by the percolation configuration. 
This procedure is illustrated in Figure \ref{fig:overall}.
The block graph may be thought of as a connected graph of lattice animals, in the sense that each block contains at least one cluster of nearest-neighbour vertices of $K$ corresponding to a lattice animal.
Specifying the block graph provides a coarse way of specifying the percolation configuration of $K$.

In \Cref{prop:block_one_arm,prop:block_cluster} we express the one-arm and cluster-size events in terms of events concerning the block graph. 
For the cluster-size event, we require that $K$ decomposes into finitely many blocks containing in total at least $k$ vertices, and that the block graph is connected via  open (long) edges between the blocks. 
The one-arm event implies that within the block decomposition of $\cluster$ there is a sequence of blocks connected by (long) edges forming a one-arm path.
For the one-arm event, the main decay term comes from a pigeonhole argument, while for the cluster-size event it comes from the $d$-dimensional isoperimetric inequality.
We use a Peierls argument, the $d$-dimensional isoperimetric inequality, and a Varopoulos-type isoperimetric inequality \cite{lyons_probability_2016} to bound the connection events in terms of the number of blocks and their total boundary size.  
In the setting of transitive graphs of polynomial growth, the Peierls argument is made possible by Tim\'ar's \cite{timar_cutsets_2007, timar_boundary-connectivity_2013} notion of coarse connectivity for minimal cutsets.
It is these last three methods where our proof for the cluster-size decay significantly differs from \cite{jorritsma_cluster-size_2024}, which relies on geometric arguments specific to the embedding of $\Z^d$ into $\R^d$.

A further difficulty when compared to \BNNP{} arises when identifying closed edges. 
Because we work with a general kernel $J$, only assumed to satisfy an asymptotic decay, we must identify closed edges of arbitrarily long length in a way that only relies on the limited knowledge that $K$ is contained in the blocks. 
We use Tim\'ar's \cite{timar_cutsets_2007,timar_boundary-connectivity_2013} notion of coarse connectivity combined with the key properties of the block graph construction in Lemma \ref{lem:one_arm_isolation} to resolve this issue, and  we identify a fixed probabilistic cost for isolating the block graph of $K$ from the rest of the graph in terms of its total boundary length.

Bounding connection and isolation events in terms of the block decomposition and its boundary size gives an expression resembling that in \eqref{eq:bnnp_lattice_animals}, see \Cref{lem:independence,prop:cluster_size_independence} below.
By choosing $\beta > \beta_c$ sufficiently large, we can sum over the possible block decompositions of $\cluster$ for any unbounded number of blocks and we obtain the desired decay for both events. 
This improves on the combinatorics in \cite{jorritsma_cluster-size_2024}, where the analogous combinatorial sum blows up so that a finite-size analysis with a bounded number of blocks is required.
\subsection{Organisation}
In Section \ref{sec:group_theory_preliminaries} we collect notation for the paper, some elementary calculations and properties concerning the \LRP{} kernel, the $d$-dimensional and the Varopoulos-type isoperimetric inequalities, and we introduce the notion of coarse connectivity in the setting of transitive graphs of polynomial growth.
In Section \ref{sec:warm_up} we prove the lower bound \Cref{lem:one_arm_lower_bound} and the preliminary result \Cref{thm:lazy_one_arm}.
In Section \ref{sec:cluster_size_lower_bound} we prove the lower bound on the size of finite clusters in \Cref{thm:cluster_size_weak_decay}.
In Section \ref{sec:anchored_isoperimetric_ineq} we prove \Cref{thm:cluster_tail_implies_iso} concerning the anchored isoperimetric inequality.
In Section \ref{sec:block_decomp} we introduce the block decomposition and we prove the isolation cost of the block graph and the Peierls argument, Lemmas \ref{lem:one_arm_isolation} and \ref{lem:peierls_argument} respectively.
In Section \ref{sec:one_arm} we prove the upper bound \Cref{thm:arm_exponent_1}.
In Section \ref{sec:cluster_decay_above} we prove  \Cref{thm:cluster_size_strong_decay,thm:s_connects_infty}.
\section{Preliminaries}
\label{sec:group_theory_preliminaries}
\subsection{Notation}
Throughout this paper, $G$ is an infinite, locally finite, transitive graph of polynomial growth with $d \geq 1$ equipped with the graph metric $d_G$ and a fixed origin $o$. 
We write $V$ for the vertex-set of $G$, $B(x,r)$ for the ball of radius $r$ centred at $x \in G$, and $S(x,r)$ for the sphere of radius $r$ centred at $x$.
We write $B(r)$ and $S(r)$ for the ball and sphere of radius $r$ centred at the origin.
For a set $A \subset V$, we write $\# A$ for the cardinality of $A$. 
We write $\R_+$ for the set $\{x \in \R : x \geq 0\}$.
The kernel $J : V \times V \to \R_+$ will generally be assumed to satisfy $J(x,y) = \Theta(d_G(x,y)^{-d \alpha})$ with $\alpha > 1$, as defined in \eqref{eq:polynomial_kernel}.
Without loss of generality we fix the constants $c_J,C_J$, and $R_J$ in \eqref{eq:polynomial_kernel}, and we write $J(x,y) = \Omega(d_G(x,y)^{-d \alpha})$ and $J(x,y) = O(d_G(x,y)^{-d \alpha})$ to mean that the kernel satisfies either the lower bound or the upper bound, respectively.
For two disjoint sets $A,B \subset V$, we write
\begin{equation}
    \label{eq:jab_abbreviation}
    J(A,B) 
    = 
    \sum_{x \in A}
    \sum_{y \in B}
    J(x,y).
\end{equation}
We write $\expdeg{x} = \E_{\beta} [\# \{y \in V : x \sim y\}]$ for the expected degree of $x \in V$ and $\expdegrad{x}{r} = \E_{\beta} [\# \{y \in V : x \sim y, d_G(x,y) = r\}]$ for the expected degree with edges of length $r \in \N$.
When we need to specify the kernel with respect to which we are considering \LRP{}, we write $\P_{\beta,J}$ and $\E_{\beta,J}$.
We write $\cluster(x)$ for the cluster of a vertex $v\in V$ and $\cluster$ for the cluster of the origin $o$.

\subsection{Kernel bounds}
By the Taylor expansion of the exponential function centred at $0$ we have $x - x^2/2 \leq 1 - \exp(-x) \leq  x$.
It follows that for all $x,y \in G$
\begin{equation}    
    \label{eq:kernel_bounds}
    1
    - 
    \exp
    \left(
        - 
        \beta
        J(x,y)
    \right)
    \leq 
    \beta 
    J(x,y),
\end{equation}  
and for all $x,y \in G$ with $d_G(x,y)$ sufficiently large such that $J(x,y) \leq 1$,
\begin{equation}
    \label{eq:kernel_lower_bound}
    \beta 
    J(x,y)/2
    \leq 
    1
    - 
    \exp
    \left(
        - 
        \beta
        J(x,y)
    \right).
\end{equation}
We assume without loss of generality that the bound in \eqref{eq:kernel_lower_bound} holds for all $x,y \in G$ with $d_G(x,y) \geq R_J$.

\subsection{Integrability, expected degree, and expected length}
\label{subsec:integrability_expected_degree}
In this section we prove that the expected degree and the expected total edge-length are finite for $\alpha > 1$ and $\alpha > 1 + 1/d$ respectively, without the use of sphere volume bounds. 
The proofs rely on the elementary observation that ``spheres partition balls''.
We group spheres into annuli, whose volume we then upper bound by the volume of the corresponding ball, for which sharp order estimates are known.
In this way, the decay of the kernel and the error we introduce for estimating the volume of spheres ``average out'' as we sum over $r$.

\begin{lemma}[Expected degree]
    \label[lemma]{lem:finite_expdeg}
    Let $G$ be a transitive graph of polynomial growth with $d \geq 1$, and let $J : V \times V \to \R_+$ be a transitive kernel with $J(x,y) = O(d_G(x,y)^{-d \alpha})$ with $\alpha > 1$. 
    Then $\sum_{y \in G \setminus \{x\}} J(x,y) < \infty$ and $\expdeg{x} < \infty$ for all $x \in G$.
\end{lemma}

\begin{proof}
    Recall from \eqref{eq:kernel_integrability_spheres} that
    \begin{equation}
        \sum_{y \in G \setminus \{x\}} 
        J(x,y)
        =
        \sum_{r = 1}^{R_J - 1}
        \sum_{y \in S(x,r)}
        J(x,y)
        +
        c_J
        \sum_{r = R_J}^{\infty}
        \# 
        S(x,r)
        r^{-d \alpha}.
    \end{equation}
    The first sum is finite.
    Decomposing the infinite sum over powers of two, and by the volume bounds for balls and the upper bound for the kernel in \eqref{eq:kernel_bounds}, 
    \begin{align}
    \label{eq:spheres_infinite_sum_1}    
        \sum_{r = R_J}^{\infty}
        \# 
        S(x,r)
        r^{- d \alpha}
        & \leq
        \sum_{i=0}^{\infty}
        \sum_{r = 2^i R_J}^{2^{i+1} R_J}
        \# S(x,r)
        r^{- d \alpha}
        \leq 
        \sum_{i=0}^{\infty}
        \# B(x,2^{i+1} R_J)
        (2^i R_J)^{-d \alpha} 
        \\
        & \leq  
        2^d
        \Cvol{}
        R_J^{d(1-\alpha)} 
        \sum_{i=0}^{\infty}
        (2^{d(1-\alpha)})^i.
    \end{align}
    The geometric series converges for $\alpha > 1$ and hence $\sum_{y \in G \setminus \{x\}} J(x,y) < \infty$ for all $x \in G$.
    Together with \eqref{eq:expdeg_first_bound}, it follows that $\expdeg{x} < \infty$ for all $x \in G$.
\end{proof}

In the proof of the decay of the truncated one-arm, it will be useful to have bounds for the expected degree of $x \in G$ with edges of length at least $r$.
We write
\begin{equation}
    \label{eq:def_expdegrad}
    \expdegrad{x}{L} 
    = 
    \E_{\beta} 
    [\# 
    \{
        y \in 
        G 
        : 
        d_G(x,y) = L,
        x \sim y
    \}    ]
\end{equation}
for the expected degree of $x \in G$ with edges of length exactly $L \in \N$.

\begin{lemma}[Expected degree $\geq r$]
    \label[lemma]{lem:exp_deg_geq_r}
    Let $G$ be a transitive graph of polynomial growth with $d \geq 1$, and let $J : V \times V \to \R_+$ be a transitive kernel with $J(x,y) = O(d_G(x,y)^{-d \alpha})$ with $\alpha > 1$. 
    Then there exist $c = c(G,J)$ and $C = C(G,J) > 0$ such that for all $x \in G$ and $r \in \N$
    \begin{equation}\label{eq:connection-to-far}
        c
        \beta
        r^{d(1-\alpha)}
        \leq 
        \sum_{L = r}^{\infty}
        \expdegrad{x}{L}
        \leq 
        C
        \beta
        r^{d(1-\alpha)}.
    \end{equation}
\end{lemma}

\begin{proof}
    Let $\rmax = \max(r,R_J)$, so that
    \begin{equation}
        \sum_{L=r}^{\infty} 
        \expdegrad{x}{L}
        \leq
        \charf
        \left(
            \rmax
            > 
            r
        \right)
        \sum_{L = r}^{\rmax}
        \expdegrad{x}{L}
        +
        c_J
        \beta
        \sum_{L = \rmax}^{\infty}
        \# 
        S(x,L)
        L^{- d \alpha}.
    \end{equation}
    The first sum is finite. The same ``spheres partition balls'' argument as in Lemma \ref{lem:finite_expdeg} yields
    \begin{align}
        \sum_{L = \rmax}^{\infty}
        \# 
        S(x,L)
        L^{- d \alpha}
        \leq 
        2^d
        \Cvol{}
        \rmax^{d(1-\alpha)} 
        \sum_{i=0}^{\infty}
        (2^{d(1-\alpha)})^i
        \leq
        c_1
        \rmax^{d(1-\alpha)} 
        \leq 
        c_1
        r^{d(1-\alpha)} 
    \end{align}
    for some $c_1 = c_1(G) > 0$, where the geometric series converges since $\alpha > 1$ and where we use that $m \geq r$.
    It follows that
    $
        \sum_{L=r}^{\infty} 
        \expdegrad{x}{L} 
        \leq 
        C
        \beta
        r^{d(1-\alpha)}
    $
    for some $C = C(G,J) > 0$.
    Similarly, using the lower bound for the kernel in \eqref{eq:kernel_lower_bound} and the known lower bounds for the volume of spheres,
    \begin{align}
        \sum_{L=r}^{\infty} 
        \expdegrad{x}{L}
        & \geq 
        \charf
        \left(
            \rmax
            > 
            r
        \right)
        \sum_{L = r}^{\rmax}
        \expdegrad{x}{L}
        +
        \frac\beta2
        \sum_{L=m}^{\infty}
        \# S(x,L) 
        L^{-d \alpha}
        \\
        &\ge  \frac\beta2
        \sum_{L=m}^{\infty} 
        L^{d-1-d \alpha}
         \geq 
        c \beta 
        r^{d(1-\alpha)}
    \end{align}
    for some $c = c(G,J) > 0$.
\end{proof}

\begin{lemma}[Expected length]
    \label[lemma]{lem:finite_exp_leng}
    Let $G$ be a transitive graph of polynomial growth with $d \geq 1$, and let $J : V \times V \to \R_+$ be a transitive kernel with $J(x,y) = O(d_G(x,y)^{-d \alpha})$ with $\alpha > 1$.
    Then for any $x \in G$ we have $\expleng{x} < \infty$ if and only if $\alpha > 1 + 1/d$, in which case there exists $c = c(G,J) > 0$ such that $\expleng{x} \leq c \beta$.
\end{lemma}

\begin{proof}
    By the upper bound on the kernel in \eqref{eq:kernel_bounds},
    \begin{align}
        \expleng{x}
        \leq
        \beta
        \left(
            \sum_{L = 1}^{R_J - 1}
            \sum_{y \in S(x,L)}
            L
            J(x,y)
            +
            c_J
            \sum_{L = R_J}^{\infty}
            \#
            S(x,L)
            L^{1 - d \alpha}
        \right).
    \end{align}
    The first sum is finite. The same ``spheres partition balls'' argument as in Lemma \ref{lem:finite_expdeg} yields
    \begin{align}
        \sum_{L = R_J}^{\infty}
        \#
        S(x,L)
        L^{1 - d \alpha}
        \leq
        \sum_{i = 0}^{\infty}
        \# 
        B(x,2^{i+1} R_J)
        (2^i R_J)^{1 - d \alpha}
        \leq
        c
        \sum_{i = 0}^{\infty}
        \left(
            2^{1 - d\alpha + d}
        \right)^i
    \end{align}
    for some $c = c(G,J)$. 
    The geometric series converges if and only if $\alpha > 1 + 1/d$.
    It follows that for all $x \in G$, $\expleng{x} < \infty$ and in particular $\expleng{x} \leq c \beta$ for some $c = c(G,J) > 0$ if and only if $\alpha > 1 + 1/d$. 
\end{proof}

\subsection{Rotationally symmetric kernel}
\label{subsec:rotation_inv_ker}
In the previous section we noted that spheres appear naturally in \LRP{}.
Similarly, it will sometimes be useful to assume that a transitive kernel $J : V \times V \to \R_+$ is \textbf{rotationally symmetric}, meaning that for any $x,y,z \in V$ with $d_G(x,y) = d_G(x,z)$ we have $J(x,y) = J(x,z)$. 
An arbitrary transitive kernel $J : V \times V \to \R_+$ is not necessarily rotationally invariant, so we define the kernel $\jmax$ given by
\begin{equation}
    \label{eq:jmax_def}
    \jmax : V \times V \to \R_+,
    \
    \jmax(x,y)  
    = 
    \max\left\{ 
        J(o,z) 
        : 
        z \in S(d_G(x,y)) 
    \right\}.
\end{equation}
By construction the kernel $\jmax$ is transitive, it is rotationally symmetric, and it satisfies $\jmax(x,y) \geq J(x,y)$ for all $x,y \in V$.
If additionally we assume that $J$ satisfies $J(x,y) = \Theta(d_G(x,y)^{-d \alpha})$ with $\alpha > 1$, it is immediate that $\jmax$ satisfies the same asymptotics.
In the arguments to come we will use $\jmax$ as a rotationally symmetric upper bound for an arbitrary kernel $J$.
When we need to specify the kernel, we write $\E_{\beta,\jmax}$ for the expectation operator of \LRP{} on $G$ with parameter $\beta$ and kernel $\jmax$.

\subsection{Isoperimetric inequality for graphs of polynomial growth}
For a finite set $A \subset V(G)$, the \textbf{boundary} of $A$ is defined to be
\begin{equation}\label{eq:boundary-def}
    \partial A 
    = 
    \{
        x \in A 
        : 
        \exists 
        y 
        \not\in 
        A, 
        d_G(x,y) = 1 
    \}.
\end{equation}
It is a consequence of a theorem of Coulhon and Saloff-Coste \cite{coulhon_isoperimetrie_1993} that a graph of polynomial growth of dimension $d$ satisfies a \textbf{$d$-dimensional isoperimetric inequality}: there exists a constant $\ciso > 0$ such that 
\begin{equation}
    \label{eq:isoperimetric_inequality}
    \# \partial A 
    \geq 
    \ciso
    (\# A)^{(d-1)/d}
\end{equation}
for every finite set $A \subset G$. 
The result of Coulhon and Saloff-Coste is generalised to locally finite transitive graphs in \cite[Proposition 5.1]{tessera_sharp_2020}.

\subsection{Bounding the placement of blocks}
In our arguments to follow, we will be faced with the following situation. 
For a finite connected set $A \subset V$ and $r \in \N$, we want to bound the cardinality of the set
$
    P(A,r) 
    =
    \{
        (x,y) 
        \in 
        A 
        \times 
        A^c
        :
        d_G(x,y) = r 
    \}
$, namely the number of pairs $x \in A$ and  $y \not\in A$   at distance exactly $r$,
in terms of $\# \partial A$.
We can write this quantity as
\begin{equation}
  \# P(A,r)=  \sum_{x \in A}
    \sum_{y \not\in A}
    \charf
    \left(
        d_G(x,y)
        =
        r
    \right).
\end{equation} 
The following bound is an immediate consequence of the $d$-dimensional isoperimetric inequality.

\begin{lemma}
    \label[lemma]{lem:trivial_isoperimetry}
    Let $G$ be a transitive graph of polynomial growth with $d \geq 2$, let $A \subset G$ be a finite subset, and let $r \in \N$.
    Then
    \begin{equation}\label{eq:pair-easy}
        \# 
        P(A,r)
        =
        \sum_{x \in A}
        \sum_{y \not\in A}
        \charf
        \left(
            d_G(x,y)
            =
            r
        \right)
        \leq 
        \ciso
        \left(
            \# 
            \partial A
        \right)^{d/(d-1)}
        \# S(r).
    \end{equation}
\end{lemma}

\begin{proof}
    It is immediate that $\# P(A,r) \leq \# A \# S(r)$ and the result follows from the $d$-dimensional isoperimetric inequality in \eqref{eq:isoperimetric_inequality}. 
\end{proof}

A mass-transport argument gives the following bound in the general setting of transitive unimodular graphs. 
We briefly review the notion of unimodularity and the mass-transport principle.
Further background may be found in \cite[Chapter 8]{lyons_probability_2016}.
Let $G = (V,E)$ be a connected, locally finite transitive graph with automorphism group $\aut(G)$.
The \textbf{modular function} is $\Delta : V^2 \to \R$ is defined to be $\Delta(u,v) = \# \stab_v u / \# \stab_u v$, where $\stab_u = \{\gamma \in \aut(G) : \gamma u = u\}$ is the stabilizer of $u$ under $\aut(G)$ and $\stab_u v = \{\gamma v : \gamma \in \stab_u\}$ is the orbit of $v$ under $\stab_u$. 
We say that $G$ is \textbf{unimodular} if $\Delta(u,v) \equiv 1$. 
It follows from \cite{soardi_amenability_1990} that transitive amenable graphs, and in particular transitive graphs of polynomial growth, are unimodular.
Suppose that $G$ is a connected, locally finite, transitive unimodular graph. 
The graph $G$ satisfies the \textbf{mass-transport principle}, which states that for every $F : V^2 \to [0,\infty]$ that is diagonally invariant in the sense that $F(\gamma u, \gamma v) = F(u,v)$ for every $u,v \in V$ and $\gamma \in \aut(G)$, we have that 
$
    \sum_{v \in V} 
    F(u,v) 
    = 
    \sum_{v \in V}
    F(v,u)
$
for any $u \in V$.

\begin{proposition}
\label[proposition]{lem:want_to_show}
    Let $G$ be a locally finite transitive unimodular graph, let $A \subset G$ be a finite subset, and let $r \in \N$. Then
    \begin{equation}
    \label{eq:want_to_show}
        \# 
        P(A,r)
        =
        \sum_{x \in A}
        \sum_{y \not\in A}
        \charf
        \left(
            d_G(x,y)
            =
            r
        \right)
        \leq  
        r
        \# \partial A
        \#S(r).
    \end{equation}
\end{proposition}

The bound in \eqref{eq:want_to_show} is sharper than that in \eqref{eq:pair-easy} for small $r$, specifically, whenever $r < (\# \partial A)^{1/(d-1)}$.
This dependence on the length of the edge and the additional factor of $r$ in \eqref{eq:want_to_show}  as compared to \eqref{eq:pair-easy} explains the transition in the driving exponent for the cluster-size decay at $\alpha=1+1/d$, see Section \ref{subsubsec:driving_exponent}.
We discuss this in detail in Remark \ref{rem:crucial_alpha_assumption}.
When $G$ is the Cayley graph of a group, the mass transport arguments involved in proving Proposition \ref{lem:want_to_show} are typically attributed to Varopoulos \cite[p.~348 Remarks]{gromov_metric_2007}, and the inequality in \eqref{eq:want_to_show} is, in at least one approach, the first step to proving the $d$-dimensional isoperimetric inequality. 
In the general setting of transitive unimodular graphs, the proof of Proposition \ref{lem:want_to_show} is implicit in \cite[Lemma 10.46]{lyons_probability_2016}.
We include a proof for completeness.

\begin{proof}[Proof of Proposition \ref{lem:want_to_show}]
    Let $A \subset G$ be finite and let $r \in \N$.
    For $x,y,z \in V(G)$, let $f_k(x,y,z)$ be the proportion of geodesic paths from $x$ to $z$ where the $k$th vertex is $y$ and let
    $
        F_{r,k}(x,y)
        = 
        \sum_{z \in S(x,r)} 
        f_k(x,y,z).
    $
    Since  $\sum_{y \in V(G)} f_k(x,y,z) = 1$, it follows that
    $
        \sum_{y \in V(G)} F_{r,k}(x,y) = \# S(x,r)
    $
    for every $x \in V(G)$ and $r \in \N$.
    Since $F_{r,k}$ is invariant under the diagonal action of $\aut(G)$, the mass-transport principle gives that 
    \begin{equation}
        \label{eq:mass_transport_applied}
        \sum_{x \in V(G)} F_{r,k}(x,y) = \# S(x,r)    
    \end{equation}
    for every $y \in V(G)$ and $r \in \N$.
    Let 
    \begin{equation}
        Z_r 
        = 
        \sum_{x \in A} 
        \sum_{z \in S(x,r) \setminus A}
        \sum_{y \in \partial A}
        \sum_{k=0}^{r-1}
        f_k(x,y,z).
    \end{equation}
    For fixed $x \in A$ and $z \in S(x,r) \setminus A$, any geodesic path from $x$ to $z\in A^c$ must pass through $\partial A$ and so the inner double sum is at least $1$. This gives that 
    \begin{equation}
        Z_r 
        \geq 
        \sum_{x \in A} 
        \sum_{z \in S(x,r) \setminus A}
        1
        =
        \# 
        P(A,r).
    \end{equation}
    Reversing the order of the summation and using the equality in \eqref{eq:mass_transport_applied} we have
    \begin{align}
        Z_r 
        &=
        \sum_{y \in \partial A}
        \sum_{k=0}^{r-1}
        \sum_{x \in A} 
        \sum_{z \in S(x,r) \setminus A}
        f_k(x,y,z)
         \leq
        \sum_{y \in \partial A}
        \sum_{k=0}^{r-1}
        \sum_{x \in V(G)} 
        \sum_{z \in S(x,r)}
        f_k(x,y,z)
        \\
        & =
        \sum_{y \in \partial A}
        \sum_{k=0}^{r-1}
        \sum_{x \in V(G)} 
        F_{r,k}(x,y)
     \leq 
        r 
        \# 
        \partial 
        A
        \# 
        S(r)
    \end{align}
    and hence $\# P(A,r) \leq r \# \partial A \# S(r)$.
\end{proof}

\subsection{Coarse connectivity in transitive graphs of polynomial growth}
\label{subsec:boundary_connectivity}
Understanding the connectedness of the boundaries of finite sets is essential to setting up Peierls-type arguments.
 Babson and Benjamini \cite{babson_cut_1999} take a homological approach to this question, and Tim\'ar \cite{timar_cutsets_2007,timar_boundary-connectivity_2013} introduces a general framework for boundary connectivity in terms of graph theory.
These results have been used in \cite{georgakopoulos_analyticity_2023} to prove the analyticity of percolation characters, in \cite{contreras_supercritical_2024} to prove the supercritical sharpness of \BNNP{}, in \cite{easo_counting_2025} to obtain (amongst other things) a new proof of Benjamini and Schramm's $p_c < 1$ conjecture, or in \cite{cannon_pirogov--sinai_2024} to study Pirogov–Sinai theory beyond the Euclidean lattice.
 A \textit{local} version of these results is obtained in \cite{martineau_percolation_2025}.
Let $G$ be a graph. 
For $x \in V(G)$, we say that $A \subset V(G)$ is a \textbf{cutset} from $x$ to $\infty$ if $x$ belongs to a finite connected component of $(V \setminus A,E)$.
If no strict subset of $A$ is a cutset between $x$ and $\infty$, we say that $A$ is a \textbf{minimal cutset} between $x$ and $\infty$.
Given $R \in \N$, we say that two vertices $x,y \in G$ are \textbf{$R$-adjacent} if $d_G(x,y) \leq R$, and we say that $A \subset V$ is \textbf{$R$-connected} if $A$ is connected as a subgraph of $G$ with the $R$-adjacency graph structure.
We will use the following consequence of Tim\'ar's result \cite[Theorem 5.1]{timar_cutsets_2007} in the setting of transitive graphs of polynomial growth, given in \cite{contreras_supercritical_2024}.
 
\begin{proposition}[{\cite[Lemma 2.1]{contreras_supercritical_2024}}]
\label[proposition]{prop:boundary_connectivity}
    Let $G$ be a transitive graph of polynomial growth with $d \geq 2$.
    There exists $R = R(G) \in \N$ such that any minimal cutset $A \subset V(G)$ from some $x \in V(G)$ and $\infty$ is $R$-connected.
\end{proposition}
\section{Lower bound and a short proof of a weaker result}
\label{sec:warm_up}
In this section we prove Lemma  \ref{lem:one_arm_lower_bound} on the lower bound the truncated one-arm and Theorem \ref{thm:lazy_one_arm}.
\Cref{thm:lazy_one_arm} follows from an elementary pigeonhole argument and gives a preliminary bound for the decay of the truncated one-arm.
In the next section we use a more sophisticated pigeonhole argument to prove the sharp bounds in \Cref{thm:arm_exponent_1}.

\begin{proof}[Proof of Lemma \ref{lem:one_arm_lower_bound}]
    Let $\mathcal{A}$ denote the event where the cluster of the origin consists of a single long edge out to $B(r)^c$, so that $\P_{\beta} \left(o \leftrightarrow B(r)^c, \# \cluster < \infty \right) \geq \P_{\beta}(\mathcal{A})$.
    By the independence of edges in \LRP{} and by the transitivity of the kernel $J$,
    \begin{align}
        \P_{\beta}
        \left(
            \mathcal{A}
        \right)
        & =
        \sum_{x \in B(r)^c} 
        \Bigg(  
            \P_{\beta}
            \left(
                o
                \sim
                x
            \right)
            \prod_{y \in G \setminus \{x\}}
            \prod_{z \in G \setminus \{o\}}
            \P_{\beta}
            \left(
                o 
                \not\sim
                y
            \right)
            \P_{\beta}
            \left(
                x 
                \not\sim
                z
            \right)
        \Bigg)
        \\
        \label{eq:one_arm_lower_bound_1}
        & =
        \exp
        \left(
            -
            2 
            \beta 
            \sum_{y \in G \setminus \{o\}}
            J(o,y)
        \right)
        \sum_{L = r}^{\infty}
        \expdegrad{x}{L}
        \geq
        c r^{d(1-\alpha)}
    \end{align}
    for some $c = c(\beta,G,J) > 0$, where in \eqref{eq:one_arm_lower_bound_1} we use that the kernel is integrable by Lemma \ref{lem:finite_expdeg} and we use Lemma \ref{lem:exp_deg_geq_r} to bound the expected degree of edges of length at least $r$.
\end{proof}

\begin{proof}[Proof of \Cref{thm:lazy_one_arm}]
    By Lemma \ref{lem:one_arm_decomp}, we have that
    \begin{equation}
        \P_{\beta}
        \left(
            o 
            \leftrightarrow 
            B(r)^c,
            \# K < \infty
        \right)
        = 
        o
        \left(
            r^{d(1-\alpha)} 
        \right) 
        + 
        \P_{\beta}
        \left(
            o 
            \leftrightarrow B(r)^c, 
            \# \cluster
            \leq 
            k
        \right)
    \end{equation}
    for $k = k(r) = \mathrm{polylog}(r)$.
    By the pigeonhole principle, the event $\{o \leftrightarrow B(r)^c, \# \cluster \leq k\}$ implies that there is at least one edge of length at least $\lceil r / k \rceil$.
    By the transitivity of the kernel $J$, by Lemma \ref{lem:exp_deg_geq_r}, and by the choice of $k$, we have
    \begin{align}
        \P_{\beta}
        \left(
            0 
            \leftrightarrow B(r)^c, 
            \# \cluster
            \leq 
            k
        \right)
        & \leq 
        k 
        \sum_{y \in B(\lceil r / k \rceil)^c}
        \P_{\beta}
        \left(
            o 
            \sim y
        \right)
        =
        k
        \sum_{L = \lceil r / k \rceil}^{\infty}
        \expdegrad{o}{L}
        \\
        & \leq 
        c
        k
        (r/k)^{d(1-\alpha)}
        =
        c 
        \mathrm{polylog}(r)
        r^{d(1-\alpha)}
    \end{align}
    for some $c = c(\beta,G,J) > 0$.
\end{proof}
\section{Cluster-size decay lower bound}
\label{sec:cluster_size_lower_bound}
In this section we prove \Cref{thm:cluster_lower_bound} on the lower bound for the distribution of finite clusters. 
The proof is by a stochastic domination argument where we renormalise to \BNNP{}.
Before we start with the proof, we need a preliminary estimate. 
Recall from \eqref{eq:jab_abbreviation} that for two disjoint sets $A,B \subset V$ we write $J(A,B) = \sum_{x \in A} \sum_{y \in B} J(x,y)$.

\begin{lemma}
\label[lemma]{lem:kernel_connection_estimate}
    Let $G$ be a transitive graph of polynomial growth with $d \geq 1$, and suppose that $J : V \times V \to \R_+$ is a transitive kernel with $J(x,y) = O(d_G(x,y)^{-d \alpha})$ with $\alpha > 1$.
    Let $\delta>0$ and $c > 0$ be such that $\#S(r)\leq c r^{d-\delta}$ for all $r \in \N$. 
    Then there exists $C = C(G,J) > 0$ such that for all $r \in \N$,
    \begin{equation}
    \label{eq:Jrr}
        J(B(r), B(r)^{c})
        \leq 
        C
        r^{
        \max
        \left(
            d(2-\alpha),
            d - 1
        \right)  
        + 
        1
        -
        \delta
        }
        \log(r)^{
            \charf
            \left(
                \alpha=1+1/d
            \right)
        }.
    \end{equation}
\end{lemma}

\begin{proof}
    Decomposing the sum over $x \in B(r)$ according to the distance of $x$ from $B(r)^c$,
    \begin{equation}
    \label{eq:J-Br-connect}
        \sum_{x \in B(r)}
        \sum_{y \in B(r)^c}
        J(x,y)
        = 
        \sum_{z = 0}^{r} 
        \sum_{x\in S(r-z)} 
        \sum_{y \in B(r)^c} 
        J(x,y)
        \leq 
        \sum_{z = 0}^{r} 
        \# S(r-z)
        \sum_{y \in B(z)^c} 
        J(o,y),
    \end{equation}
    where we use the transitivity of the kernel $J$ and the fact that for any $x \in S(r-z)$ we have $d_G(x,B(r)^c) \geq z$.
    For $z \in \N$, it follows from Lemma \ref{lem:exp_deg_geq_r} that there exists $c_1 = c_1(G,J) > 0$ such that
    \begin{equation}
        \sum_{y \in B(z)^c} 
        J(o,y)
        =
        \sum_{L = z + 1}^{\infty}
        \sum_{y \in S(L)}
        J(o,y)
        \leq 
        c_1 
        (z + 1)^{d(1 - \alpha)}.
    \end{equation}
    Together with the assumed upper bound on the volumes of spheres we have shown that
    \begin{equation}
        \label{eq:Br-connection-2}
        J(B(r), B(r)^{c}) 
        \leq
        c_1 
        \sum_{z = 0}^{r} 
        \#S(r-z) 
        (z+1)^{d(1 - \alpha)} 
        \leq  
        c_1 
        r^{d - \delta}
        \sum_{z = 0}^{r} z^{d(1 - \alpha)},
    \end{equation}
    where we use that $(r-z)^{d - \delta} \leq r^{d - \delta}$ since $d - \delta > 0$, and that $(z+1)^{d(1 - \alpha)} \leq z^{d(1 - \alpha)}$ since $d(1 - \alpha) < 0$.
    There exists $c_3 > 0$ such that
    \begin{equation}
    \label{eq:binomial-sum-2}
        \sum_{z = 0}^{r} 
        z^{d(1 - \alpha)} 
        \leq
        \begin{cases} 
            c_3 
            /
            (d(1 - \alpha)) 
            & 
            \text{if } 
            d(1 - \alpha) 
            <
            -1,
            \\
            c_3
            \log r 
            & 
            \text {if } 
            d(1 - \alpha)
            =
            -1,
            \\
            c_3 
            r^{d(1 - \alpha) + 1}
            / 
            (d(1 - \alpha) + 1) 
            & \text {if } 
            d(1 - \alpha)
            >
            - 1,
        \end{cases} 
    \end{equation}
    and the result follows.
\end{proof}

Although we do not need it for our arguments, we can also show a matching lower bound on $J(B(r),B(r)^c)$, conditionally on the conjectured volume bounds for spheres (see Section \ref{subsubsec:spheres_in_graphs}).

\begin{lemma}
     Let $G$ be a transitive graph of polynomial growth with $d \geq 1$, and suppose that $J : V \times V \to \R_+$ is a transitive kernel with $J(x,y) = \Omega(d_G(x,y)^{-d \alpha})$ with $\alpha > 1$.
     Then there exists $c = c(G,J) > 0$ such that for all $r \in \N$
     \begin{equation}
    \label{eq:Jrr-lower}
        J(B(r), B(r)^c) 
        \geq
        c
        r^{
            \max
            \left(
                d(2-\alpha),
                d-1
            \right)  
        }
        \log(r)^{
            \charf
            \left(
                \alpha=1+1/d
            \right)
        }.
\end{equation}
\end{lemma}

\begin{proof}
    Suppose that $r \geq 2 R_J$. 
    For $z \leq r$ and $x \in S(r-z)$, let $x_z$ be the vertex on $S(r + 2z)$ closest to $x$, so that $d_G(x,x_z) = 3 z$ and $B(x_z,z) \subset B(r)^c$.
    If additionally $2 z \geq r/2$ so that $R_J \leq d_G(x,y) \leq 4z$ for $y \in B(x_z,z)$, then by the lower bound on the kernel
    \begin{equation}
        \sum_{y \in B(r)^c} 
        J(x,y)
        \geq 
        \sum_{y \in B(x_z,z)} 
        J(x,y)
        \geq 
        c_1 
        z^{d(1-\alpha)}
    \end{equation}
    for some $c_1 > 0$.
    Decomposing $J(B(r),B(r)^c)$ as in \eqref{eq:J-Br-connect} we have
    \begin{align}
        J(B(r),B(r)^c)
        & \geq
        \sum_{z = r/4}^{r} 
        \sum_{x\in S(r-z)} 
        \sum_{y \in B(r)^c} 
        J(x,y)
        \geq 
        c_2
        \min_{r/4 \leq z \leq r} 
        \left(
            \# 
            S(r-z)
        \right)
        \sum_{z = r/4}^{r}
        z^{d(1-\alpha)}
        \\
        & \geq
        c_2
        r^{d-1}
        \sum_{z = r/4}^{r}
        z^{d(1-\alpha)}
    \end{align}
    for some $c_2 > 0$. 
    The statement follows by bounding the sum over $z$ as in \eqref{eq:binomial-sum-2}.
\end{proof}

The next proposition proves that in the ball of radius $r$ we find that the origin is in a linear sized component with large enough probability.

\begin{proposition}
    \label[proposition]{prop:giant_concentration}
    Let $G$ be a transitive graph of polynomial growth with $d \geq 1$, and suppose that $J : V \times V \to \R_+$ is a transitive kernel with $J(x,y) = \Omega(d_G(x,y)^{-d \alpha})$ with $\alpha > 1$. 
    Let $\beta > \beta_c$ be sufficiently large.
    Then there exist $\rho = \rho(\beta,G,J) > 0$ and $c = c(\beta, G,J) > 0$ such that for all $r$ sufficiently large,
    \begin{equation}
        \P_{\beta}
        \left(
            \# 
            \cluster
            (r)
            \geq 
            \rho
            \# B(r) 
        \right)
        \geq 
        \exp
        \left(
            - c  
            \# S(r)
        \right).
    \end{equation}      
\end{proposition}

\begin{proof}
    Note that $S(r) = \partial B(r-1)$ is a minimal cutset, and by Proposition \ref{prop:boundary_connectivity} there exists $R = R(G) > 0$ such that $S(r)$ is $R$-connected.
    Let $T$ be a spanning tree of the graph of $S(r)$ with respect to $R$-adjacency, and let $H$ be the graph obtained by taking the union of $G$ and $T$. 
    Clearly $G \subset H$, and we write $d_H$ for the graph metric associated to $H$.
    Note that the vertex set of the ball of radius $r$ centred at the origin in $G$ is the same as that in $H$, and we write $B(r)$ irrespectively of the metric.
    Let $c_J,R_J > 0$ be such that $J(x,y) \geq c_J d_G(x,y)^{-d\alpha}$ for all $x,y \in G$ with $d_G(x,y) \geq R_J$.
    To each pair of vertices in $B(r)$ with $d_H(x,y) = 1$, we assign a vertex $z_{xy} \in B(r)$ satisfying that $d_G(x,z_{xy}) = R_J + 1$ and $ R_J \leq d_G(y,z_{xy}) \leq 3 \max(R,R_J)$, as follows.
    If $d_G(x,y) = 1$, then either $d_G(o,x)\le R_J+1$ or $d_G(o,x)>R_J+1$. 
    For the first case let $z_{xy}$ be any vertex with $d_G(x,z_{xy}) = R_J + 1$. 
    For the second case let $z_{xy}$ be a vertex closest to $x$ on the sphere $S(o, d_G(o,x)-(R_J+1))$.
    Since $r \geq 2(R_J + 1)$ we have that $z_{xy} \in B(r)$ for both cases, and by the triangle inequality we have that $R_J \leq d_G(y,z_{xy}) \leq 3 R_J$.
    If $d_G(x,y) \neq 1$, so that by construction $x,y \in S(r)$ with $d_G(x,y) \leq R$, let $z_{xy}$ be the vertex in $S(r - (R_J + 1))$ nearest to $x$.
    It follows that $z_{xy} \in B(r)$, $d_G(x,z_{xy}) = R_J + 1$, and by the the fact that $y \in S(r)$ together with the triangle inequality we have that $R_J + 1 \leq d_G(y,z_{xy}) \leq 3\max(R,R_J)$.

    Given a \LRP{} configuration $\omega$ on $G$ with parameter $\beta$, we define a $1$-dependent bond percolation $\eta$ on the subgraph of $H$ induced by $B(r)$ as follows. 
    For $x,y \in B(r)$ with $d_H(x,y) = 1$, we say that the bond $(x,y)$ is open and write $x \sim_{\eta} y$ if and only if both $(x,z_{xy})$ and $(z_{xy},y)$ are open in $\omega$.
    It follows from the choice of $z_{xy}$ that
    \begin{equation}
        \P_{\eta}
        \left(
            x 
            \sim_{\eta}
            y
        \right)
        = 
        \P_{\beta}
        \left(
            x 
            \sim 
            z_{xy}
            \sim
            y
        \right)
        \geq 
        \Big(1 - e^{-\beta c_J (3\max(R,R_J))^{-d \alpha}}\Big)^2
        =: 1
        -
        q(\beta).
    \end{equation}
    To see that $\eta$ is $1$-dependent, let $x_1,x_2,y_1,y_2 \in B(r)$ be distinct vertices with $d_H(x_i,y_j) = 1$.
    The events $x_1 \sim_{\eta} y_1$ and $x_1 \sim_{\eta} y_2$ are not independent, since it is possible that $z_{x_1 y_1} = z_{x_1 y_2}$ and then both $x_1 \sim_{\eta} y_1$ and $x_1 \sim_{\eta} y_2$ depend on the same edge in $\omega$.
    The events $x_1 \sim_{\eta} y_1$ and $x_2 \sim_{\eta} y_2$ are independent however, since even if $z_{x_1 y_1} = z_{x_2 y_2}$ the events $x_1 \sim_{\eta} y_1$ and $x_2 \sim_{\eta} y_2$ depend on four different edges in $\omega$.
        Note that the graph $H$ has bounded degree.
    By \cite{liggett_domination_1997}, $\eta$ can be stochastically dominated from below by an independent bond percolation $\eta^\star$ on the subgraph of $H$ induced by $B(r)$ with parameter $p = p(\beta) = 1 - c q(\beta)$ for some $c > 0$. 
        In fact, $c$ can be chosen to be the twice the maximal degree in $H$ in $B(r)$, which is at most $2(\# B(1)+\#B(R))$.
    The process $\eta^\star$ naturally induces an independent bond percolation $\eta^{\star \star}$ on $G$ with parameter $p(\beta)$.
    Further, the critical parameter $p_c(G)$ of the process $\eta^{\star \star}$ satisfies $p_c(G) < 1$ since $G$ is a transitive graph of polynomial growth with $d \geq 2$. 
    Let $\beta$ be sufficiently large so that $p = p(\beta) > p_c(G)$, and let $\theta^{\star \star} > 0$ be the percolation probability in $\eta^{\star \star}$.
    Let $\mathcal{A}(\eta^*)$ denote the event that at least $\# B(r) \theta^{\star\star} / 2$ vertices in $B(r)$, including the origin $o$, have a path of open bonds to $S(r)$ in $\eta^*$, and let $\mathcal{A}(\eta^{\star \star})$ denote the event 
    $
    \{
        \# 
        \left( 
            \cluster_\infty(\eta^{\star \star}) 
            \cap 
            B(r) 
        \right) 
        \geq 
        \# 
        B(r) 
        \theta^{\star \star} 
        / 
        2 
    \}
    $ in $\eta^{\star \star}$, where we write $\cluster_\infty(\eta^{\star \star})$ for the almost surely unique infinite cluster in $\eta^{\star \star}$.
    By construction we have that
    \begin{equation}
        \P_{\eta^{\star}}
        \left(
            \mathcal{A}
            (\eta^{\star})
        \right) 
        \geq 
        \P_{\eta^{\star\star}}
        \left(
            \mathcal{A}(\eta^{\star\star}) 
            \cap 
            \{o \in \cluster_\infty(\eta^{\star \star})\} 
        \right),
    \end{equation}
    and it follows from \cite[Lemma 8.68]{grimmett_percolation_1999} that
    $
        \P_{\eta^{\star \star}}
        \left(
            \mathcal{A}(\eta^{\star \star})
        \right)
        \geq 
        \theta^{\star \star}
        /
        2
    $.
    Both $\mathcal{A}(\eta^{\star\star})$ and $\{o \in \cluster_\infty(\eta^{\star \star})\}$ are increasing events in $\eta^{\star \star}$, and by the FKG inequality \cite[Theorem 2.4]{grimmett_percolation_1999},
    \begin{align}
        \P_{\eta^{\star\star}}
        \left(
            \mathcal{A}(\eta^{\star\star})
            \cap 
            \{
                o \in \cluster_\infty(\eta^{\star \star})
            \}
        \right)
        \geq 
        \P_{\eta^{\star\star}}
        \left(
            \mathcal{A}(\eta^{\star\star})
        \right)
        \P_{\eta^{\star\star}}
        \left(
                o \in \cluster_\infty(\eta^{\star \star})
        \right)
        \geq 
        (\theta^{\star \star})^2/2.
    \end{align}
    Let $\mathcal{E}(\eta^\star)$ denote the event that all edges in $T$ are open in $\eta^\star$, so that in particular $S(r)$ is connected in $\eta^*$.
    The event $\mathcal{E}(\eta^\star) \cap \mathcal{A}(\eta^{\star})$ implies that $\# \cluster(r) \geq \# B(r) \theta^{\star \star} / 2$ in $\eta^\star$.
    Both the events $\mathcal{E}(\eta^\star)$ and $\mathcal{A}(\eta^{\star})$ are increasing events in $\eta^\star$, and by the FKG inequality we have
    $
        \P_{\eta^\star}
        \left(
            \mathcal{E}(\eta^\star) 
            \cap 
            \mathcal{A}(\eta^{\star}) 
        \right)
        \geq 
        p^{\#S(r) - 1}
        (\theta^{\star \star})^2
        /
        2
    $.
    It follows from the construction of $\eta$ from $\omega$ and from the coupling given by the stochastic domination between $\eta$ and $\eta^\star$ that $\# \cluster(r,\omega) \geq \# \cluster (r,\eta) \geq \# \cluster(r,\eta^\star)$.
    We have shown that
    \begin{align}
        \P_\beta
        \left(
            \# 
            \cluster(r,\omega)
            \geq
            \# 
            B(r)
            \theta^{\star \star}
            /
            2 
        \right)
        \geq 
        \P_{\eta}
        \left(
            \# 
            \cluster(r,\eta)
            \geq
            \# 
            B(r)
            \theta^{\star \star}
            /
            2 
        \right)
        \geq 
        p^{\#S(r) - 1}
        (\theta^{\star \star})^2
        /
        2,
    \end{align}
    and the result follows by setting $\rho = \rho(\beta, G, J) > 0$ and $c = c(\beta, G,J) > 0$.
\end{proof} 

In the proof of Proposition \ref{prop:giant_concentration} above, we renormalise to a bond percolation process. 
We do this since the argument where we apply \cite[Lemma 8.68]{grimmett_percolation_1999} to construct a linear-sized cluster in $B(r)$ cannot be applied directly in \LRP{} as edge lengths may be unbounded and a path from some vertex in $B(r)$ to $\infty$ must not necessarily cross  $S(r)$, or any larger sphere of fixed radius for that matter.

\begin{proof}[Proof of \Cref{thm:cluster_lower_bound}]
    Let $r = (Mk/\cvol{})^{1/d}$ for some $M > 0$ to be chosen later, so that $\# B(r) \geq M k$ by the volume bounds for balls in \eqref{eq:ball_volume_bounds}. 
    Writing $\cluster(r)$ for the cluster of the origin in the ball of radius $r$,
    \begin{align}
    \label{eq:lower-bound-start}
        \P_{\beta}
        \left(      
            k 
            \leq 
            \# 
            \cluster
            < 
            \infty
        \right) 
        & \geq 
        \P_{\beta}
        \left(
            \# 
            \cluster(r)
            \geq 
            k,
            B(r)
            \not\sim
            B(r)^c
        \right)
        \\
        & \geq 
        \P_{\beta}
        \left(
            \# \cluster(r)
            \geq 
            k
        \right)
        \P_{\beta}
        \left(
            B(r)
            \not\sim
            B(r)^c
        \right)\label{eq:lower-bound-product}
    \end{align}    
    where we use the independence of edges in \LRP{}.
    Using the upper bound on $J(B(r), B(r)^c)$ in  \Cref{lem:kernel_connection_estimate} gives that
    \begin{equation}
    \begin{aligned}
        \P_{\beta}
        \left(
            B(r)
            \not\sim
            B(r)^c
        \right)
        &=
        \exp
        \left( 
            - 
            \beta 
            J
            \left(
                B(r),
                B(r)^c
            \right)
        \right)
    \\
        &\geq
        \exp
        \left(
            - 
            c_1
            r^{
                \max(d(2-\alpha), d-1)+1-\delta
            } (\log r)^{\charf(\alpha=1+1/d)}
        \right)
        \end{aligned}
    \end{equation}     
    for some $c_1 = c_1(\beta, G,J) > 0$.
    It follows from \Cref{prop:giant_concentration} that there exist $\rho = \rho(\beta,G,J) > 0$ and $c = c(\beta,G,J) > 0$ such that for $M=1/\rho$ and $k$ sufficiently large,
    \begin{equation}
        \P_{\beta}
        \left(
            \# \cluster(r)
            \geq 
            k
        \right)
        \geq 
        \P_{\beta}
        \left(
            \# 
            \cluster(r)
            \geq 
            \rho 
            \# 
            B(r)
        \right)
        \geq
        \exp
        \left(
            - 
            c_2
            \# 
            S(r)
        \right),
    \end{equation}
    By the choice of $r$ together with the assumption on the volume of spheres we arrive at
    \begin{align}
        \P_{\beta}
        \left(      
            k 
            \leq 
            \# 
            \cluster
            < 
            \infty
        \right)&\geq
        \exp
        \left(
            -c_3
            \left(
                r^{
                    \max(d(2-\alpha), d-1)+1-\delta
                } 
                (\log r)^{
                    \charf(\alpha=1+1/d)
                }
                +
                r^{d-\delta}
            \right)
        \right)
        \\
        & \geq
        \exp
        \left(
            -
            A
            k^{\max(2-\alpha, (d-1)/d) +(1-\delta)/d} 
            (\log k)^{
                \charf(\alpha=1+1/d)
            }
        \right)
    \end{align}
    for some $c_3 = c_3(\beta,G,J) >0$ and $A = A(\beta,G,J) > 0$ for all $k \in \N$.    
\end{proof}

\begin{remark}
    \label[remark]{rem:sphere_inevitability}
    The event that realises the lower bound demands a localised cluster in a ball of size $\Theta(k)$ and radius $r=\Theta(k^{1/d})$. Near the boundary of the ball we need $\#S(r)$ many short edges to be closed, which means that in this argument it is inevitable to use upper bounds on the sphere-sizes.
\end{remark}

\section{Anchored isoperimetric inequality}
\label{sec:anchored_isoperimetric_ineq}
In this section we prove \Cref{thm:cluster_tail_implies_iso}.
In \eqref{eq:anchored_isop_ineq} we defined the anchored isoperimetric inequality in terms of $\# W$, whereas \cite{pete_note_2008,hutchcroft_transience_2023} define it in terms of $\sum_{w \in W} \deg(w)$. 
Since percolation clusters in \LRP{} can have unbounded degree, this choice seems to be the most natural one in our setting. 
For \BNNP{} on bounded degree graphs $G$, this difference amounts to modifying the constant $c$.

\begin{proof}[Proof of \Cref{thm:cluster_tail_implies_iso}]
In \eqref{eq:anchored_isop_ineq}, for sets $W$ of size $n$, one can replace the function $n^{(d-1)/d}$ with an arbitrary function $\phi(n)$, then we say that that $\mathcal C_\infty$ satisfies an anchored $\phi(n)$-isoperimetric inequality. 
We show that $\P_{\beta} \left(k < \# \cluster <\infty \right) \leq \exp \left(-c k^{\zeta}\right)$ implies that the infinite component of long-range percolation satisfies a $n^{\zeta}/\log n$ isoperimetric inequality a.s. and thus it also satisfies a $n^{\zeta'}$-inequality for all $\zeta'<\zeta$ a.s. This implies, by its definition, that the anchored isoperimetric dimension is $\zeta$ a.s. 

Let us define $\mathcal V_n^{(o)}:=\{H\ni o, H=(V(H), E(H)), \#H = n \}$ the set of possible realisations of a finite component of $o$ containing $n$ vertices in long-range percolation. Note that $H$ is a connected graph in $\mathrm{LRP}_\beta$ but it may contain long edges so its vertex set is not necessarily connected in $G$. 
Let $K$ be a component of $o$ in $\mathrm{LRP}_\beta$. 
Our goal is to show that for some $\delta=\delta(\beta, G,J)>0$
\begin{equation}\label{eq:anchored-complement}
\P_\beta\Big(\exists H \in \mathcal V_n^{(o)}: H\subseteq K, \#\partial_{E(K)} H \le \delta\frac{n^{\zeta}}{\log (n)} \Big)\le C \exp\Big(-\frac{c(\beta)}{2} n^{\zeta}\Big).
\end{equation}
In the event on the left hand side, $K$ may also be finite. In that case, the bound still holds, as then the probability that $K$ has size at least $n$ is bounded by the right-hand side. 
Let 
\begin{equation}\label{eq:new-iso}
\psi(n):=\frac{n^{\zeta}}{\log (n)}
\end{equation}
As the right hand-side of \eqref{eq:anchored-complement} is summable in $n$, by the Borel Cantelli lemma there is a.s.\ an $n_o=n_o(\beta, G,J)$ such that the complement of the event in \eqref{eq:anchored-complement} holds for all sets of size at least $n_o$. On  the event that $K$ is infinity, which happens with strictly positive probability $\theta(\beta)>0$,  any finite graph $H$ subgraph of $K$ has at least one edge on its edge-boundary $\partial_{E(K_\infty)}H$, by choosing $\delta$ smaller if necessary to account for sets of size at most $n_o$, the $\psi(n)$-anchored isoperimetric inequality thus holds also  for all smaller subsets of $K_\infty$ containing $o$.

Given $H\in \mathcal V_n^{(o)}$, let $S$ be a finite set of $m$ pairs of vertices (edges) so that each edge in $S$ is adjacent to $V(H)$.  For an edge $e=(x,y)$ let us write $J(e):=J(x,y)$ for brevity. Let us also write $\|J\|_1:=\sum_{y\in G\setminus\{o\}} J(o,y)$.
Then we can compute using the transitivity of the kernel $J$ that 
\begin{equation}
\P_\beta\big(H\subseteq K, \partial_{E(K)}H =S\big) =  e^{-\beta n\|J\|_1 + \beta \sum_{e\in E(H), e\in S} J(e) } \prod_{e\in E(H), e\in S} \Big(1-e^{-\beta J(e)}\Big),
\end{equation}
indeed, all edges adjacent to $V(H)$ except those in $E(H)$ or in $S$ need to be closed, while edges in $E(H)$ and in $S$ need to be all open. 
We can rearrange this formula to
\begin{equation}
\begin{aligned}
\P_\beta\big(H\subseteq K, \partial_{E(K)}H =S\big) &=  e^{-\beta n\|J\|_1 + \beta \sum_{e\in E(H)} J(e)} \prod_{e\in E(H)} \Big(1-e^{-\beta J(e)}\Big) \\
&\quad \cdot \prod_{e\in S}  e^{ \beta J(e)}   \Big(1-e^{-\beta J(e)}\Big)\\
&=\P_\beta(K=H)\cdot \prod_{e\in S}  e^{ \beta J(e)}   \Big(1-e^{-\beta J(e)}\Big),
\end{aligned}
\end{equation}
The last formula closely resembles \cite[Eq.\,(3.2)]{hutchcroft_transience_2023} except here the edge-probabilities vary. We now classify $S$ according to these varying edge-probabilities.
For each vertex $x$ in $G$, let us fix an automorphism $\gamma_{x}$ so that $\gamma_{x} x=o$ and then $\gamma_{x} y \in S(o,d_G(x,y))$. Let us then call the edge $(x,y)$ of type $z=\gamma_{x}y$. Note that for each vertex on $y\in S(x,r)$, there is exactly one vertex $z$ in $S(o,r)$ so that the edge $(x,y)$ is of type $z$. Let us define the profile of $S$ as a vector $\mathbf m=(m_{z})_{z\in G\setminus\{o\}}$, so that $S$ contains $m_z$ edges of type $z$. Clearly  $\|\mathbf m\|_1=\sum_{z\in G\setminus \{o\}} m_z=m$ as there are $m$ edges in $S$ per our assumption. 
The quantity $\|J\|_\infty=\max_{y\in G\setminus\{o\}}J(o,y)$ is finite by the integrability of the kernel $J$ on $G$.
Then by the transitivity of the kernel $J$, it holds that
\begin{equation}
    \prod_{e\in S}  e^{ \beta J(e)}   \Big(1-e^{-\beta J(e)}\Big) \le   \prod_{z\in G\setminus\{o\}}   \Big( e^{\beta \|J\|_\infty}\Big(1-e^{-\beta  J(o,z)}\Big)\Big)^{m_z}.
\end{equation}
We thus obtain a uniform upper bound for all $S$ with the same profile $\mathbf m$.
Still fixing the graph  $H$ with $\#H=n$, now we sum over the possible choices of $S$ of the same size $m$. Note that if $S$ contains $m_z$ edges of type $z$, then these edges must be chosen from at most $n$ many possible choices, as for each vertex $x\in H$ there is exactly one vertex $y$ of type $z$. So we obtain:
\begin{equation}
\begin{aligned}
    \P_\beta
    \Big(
        H\subseteq K, 
        &
        \#\partial_{E(K)}H = m
    \Big) 
    = 
    \P_\beta(K=H)
    \sum_{S: \#S=m} 
    \prod_{e\in S}  
    e^{ \beta J(e)}   
    \Big(
        1-e^{-\beta J(e)}
    \Big)
    \\
    &\le
    \P_\beta(K=H) 
    \sum_{
        \mm
        : 
        \normone{\mm} = m
    }
    \prod_{z \in G \setminus \{o\}} 
    \binom{n}{m_z}
    \Big(
        e^{ \beta\| J\|_\infty}   
        \Big(
            1-e^{-\beta  J(o,z)}
        \Big)
    \Big)^{m_z}
    \end{aligned}
\end{equation}
We now use the upper bound that $\binom{n}{m_z} \le n^{m_z}$, which after elementary rearrangement and $\|\mathbf m\|_1=m$ yields that 
\begin{equation}
\begin{aligned}
    \P_\beta\Big(H\subseteq K, &\#\partial_{E(K)}H = m\Big) \le \P_\beta(K=H) \big(n e^{\beta \|J\|_\infty}\big)^m \\
    &\qquad \cdot  \sum_{\mathbf m: \|\mathbf m\|_1=m}  \prod_{z\in G\setminus\{o\}}  \Big(1-e^{-\beta J(z)}\Big)^{m_z}\\
&\le \P_\beta(K=H) \big(n e^{\beta \|J\|_\infty}\big)^m \bigg(\sum_{z\in G\setminus\{o\}}\Big(1-e^{-\beta J(o,z)}\Big)\bigg)^m,
\end{aligned}
\end{equation}
where we noticed that the product in the middle row is at most the distribution of the power in the last row. The expression in the sum is exactly $\E_{\beta}[\deg(o)]$, the expected degree of $o$ in LRP, which is assumed to be finite. 
We arrive at the simple expression
\begin{equation}\label{eq:simple-bound}
     \P_\beta\Big(H\subseteq K, \#\partial_{E(K)}H = m\Big) \le \P_\beta(K=H) \Big( n e^{\beta \|J\|_\infty} \E_{\beta}[\deg(o)]\Big)^m. 
\end{equation}
Recall $\psi(n)$ from \eqref{eq:new-iso} and sum over $m\le \delta\psi(n)$. Choosing $n$ large enough so that the base of $m$ in \eqref{eq:simple-bound}  
is at least $2$, and writing $e^{\beta \|J\|_\infty} \E_{\beta}[\deg(o)]:=M_{\beta, J}$, we arrive at the simple expression using the bound on the geometric series that
\begin{equation}\label{eq:h-k-intermed}
\begin{aligned}
\P_\beta\Big(H\subseteq K, \#\partial_{E(K)}H \le \delta \psi(n)\Big)&\le \P_\beta(K=H)\sum_{m=0}^{\delta\psi(n)} \Big( n e^{\beta \|J\|_\infty} \E_{\beta}[\deg(o)]\Big)^m\\
&\le \P_\beta(K=H) 2 \Big(n M_{\beta,J}\Big)^{\delta\psi(n)}.
\end{aligned}
\end{equation}
Using now the from of $\psi(n)$ from \eqref{eq:new-iso} 
\[ 
\begin{aligned}
\Big(n M_{\beta,J}\Big)^{\delta\psi(n)} &= \exp\Big(\delta \frac{n^{\zeta}}{\log n} \log n + \delta \frac{n^{\zeta}}{\log n}\log M_{\beta,J}\Big)
\le \exp( 2\delta n^{\zeta})
\end{aligned}
\]
for all $n$ with $n\ge M_{\beta,J}=e^{\beta \|J\|_\infty} \E_{\beta}[\deg(o)]$.
Returning to \eqref{eq:h-k-intermed}, we obtain that for all graphs $H$ with $\#H=n$,
\[ 
\P_\beta\Big(H\subseteq K, \#\partial_{E(K)}H \le \delta \psi(n)\Big) \le \P_\beta(K=H)\cdot e^{2\delta n^{\zeta}}.
\]
Taking now a union bound over all $H$ in $V_n^{(o)}$, we arrive at:
\begin{equation}
    \begin{aligned}
\P_\beta\Big(\exists H \in \mathcal V_n^{(o)}: H\subseteq K, \#\partial_{E(K)}H &\le \delta \psi(n) \Big) \le e^{2\delta n^{\zeta}}\sum_{H\in  \mathcal V_n^{(o)}} \P_\beta(K=H)\\
    &=e^{2\delta n^{\zeta}} \P_\beta(\#K(o)=n).
    \end{aligned}
\end{equation}
Choosing now $\delta:=c(\beta)/4$ where $c(\beta)$ is as in the tail of the distribution of finite clusters finishes the proof of \eqref{eq:anchored-complement} and thus also that of the theorem.
\end{proof}
\section{Block decomposition and the block graph construction}
\label{sec:block_decomp}
In this section we define the block decomposition and the block graph construction, due to \cite{jorritsma_cluster-size_2024}. 

\subsection{Block decomposition}
\label{subsec:block_decomp}
Let $G = (V,E)$ be a connected graph. 
We say that a non-empty finite set $A \subset V$ of vertices is \textbf{$1$-connected} if the induced graph $G(A)$ consists of a single connected component.
Note that being $1$-connected is equivalent to being connected, and the definition of $1$-connected is coherent with the definition of $R$-connected in Section \ref{subsec:boundary_connectivity}.
For $b \geq 2$, a sequence of $1$-connected sets $(A_0,\ldots,A_{b-1})$ is \textbf{$1$-disconnected} if $d_G(A_i, A_j) > 1$ for all $i \neq j$.
The following lemma establishes that a finite set has a unique decomposition into finite connected components. 

\begin{lemma}[Uniqueness of blocks]
\label[lemma]{lem:uniqueness_and_connectedness}
    Let $C \subset V$ be a non-empty finite set of vertices.
    Then there exists $b \in \N$ and a unique $1$-disconnected sequence of $1$-connected sets $(A_0,\ldots,A_{b-1})$ such that $C = A_0 \sqcup \ldots \sqcup A_{b-1}$.
\end{lemma}

\begin{proof}
    Define the equivalence relation $\leftrightarrow_C$ on $C$ where $x \leftrightarrow_C y$ if and only if there is a connected path from $x$ to $y$ in $C$.
    The set $C$ consists of the blocks given by the equivalence classes of $\leftrightarrow_C$. 
    The uniqueness of the sequence follows since $C$ is finite and $\leftrightarrow_C$ is an equivalence relation.
\end{proof}

Recall from Section \ref{subsec:boundary_connectivity} that $A \subset V$ is a cutset from $x$ to $\infty$ if $x$ belongs to a finite connected component of $(V \setminus A,E)$.
We define the \textbf{closure} of $A$ as  
\begin{equation}
    \label{eq:closed_blocks}
    \overline{A} 
    =
    A \cup
    \{
        x \in G
        :
        A 
        \text{ is a cutset from } 
        x 
        \text{ to } 
        \infty
    \}.
\end{equation}
We say that a non-empty finite set $A \subset V$ is a \textbf{block} if it is $1$-connected and $\overline{A} = A$, we let 
\begin{equation}
    \label{eq:blocks}
    \AA
    =
    \{
        A \subset V
        : 
        A \text{ is } 1 \text{-connected},
        \overline{A} = A
    \}
\end{equation}  
be the set of blocks, and we let 
\begin{equation}
    \label{eq:disconnected_blocks}
    \AA(b)
    =
    \{
        \ablocks
        =
        (A_0,\ldots,A_{b-1}) \in \AA^b
        :
        \ablocks
        \text{ is }
        1 
        \text{-disconnected}
    \}
\end{equation}  
be the set of sequences of $1$-disconnected blocks.
In the following lemma we show that a $1$-disconnected sequence of $1$-connected sets is always contained in a unique sequence of $1$-disconnected blocks with nesting boundaries.

\begin{lemma}[Closed blocks]
\label[lemma]{lem:closed_blocks}
    Let $b \in \N$ and let $A_0,\ldots,A_{b-1}$ be a $1$-disconnected sequence of $1$-connected sets.
    Then there exists $\ell \leq b$ and $(B_0,\ldots,B_{\ell-1}) \in \AA(\ell)$ such that $\cup_{i \leq b} A_i \subset \cup_{j \leq \ell} B_j$ and $\cup_{j \leq \ell} \partial B_j \subset \cup_{i \leq b} \partial A_i$.
\end{lemma}

\begin{proof}
    We begin by showing that for any $1$-connected set $B \subset V$ we have $\partial \overline{B} \subseteq \partial B$.
    Let $y \in \partial \overline{B}$.
    It follows from the definitions of the boundary and the closure of a set that $y \in B \cup \{x \in G : B \text{ is a cutset from } x \text{ to } \infty\}$ with $d_G(y,z) = 1$ for some $z \not\in B \cup \{x \in G : B \text{ is a cutset from } x \text{ to } \infty\}$.
    Suppose that $y \in \{x \in G : B \text{ is a cutset from } x \text{ to } \infty\}$.
    By the definition of a cutset, there is no $z \not\in B \cup \{x \in G : B \text{ is a cutset from } x \text{ to } \infty\}$ such that $d_G(y,z) = 1$.
    It follows that $y \in B$ with $d_G(y,z) = 1$ for some $z \not\in B \cup \{x \in G : B \text{ is a cutset from } x \text{ to } \infty\}$.
    In particular, $y \in B$ with $d_G(y,z) = 1$ for some $z \not\in B$, so that $y \in \partial B$ and hence $\partial \overline{B} \subseteq \partial B$.

    Suppose now that $B_i$ and $B_j$ are two distinct $1$-disconnected, $1$-connected sets.
    We show that the closures $\overline{B}_i$ and $\overline{B}_j$ are either $1$-disconnected blocks, or one contains the other.
    Suppose in the first instance that $\overline{B}_i \cap \overline{B}_j \neq \varnothing$, and suppose further that $B_i \subseteq \overline{B}_j$. 
    Since the $B_i$ and $B_j$ are $1$-disconnected, we have $B_i \cap B_j = \varnothing$. 
    If $x \in (\overline{B}_j \setminus B_j) \cap B_i$, then $x \in \overline{B}_j$.
    If $x \in (\overline{B}_j \setminus B_j) \cap (\overline{B}_i \setminus B_i)$, then $B_i$ is a cutset from $x$ to $\infty$ and since $B_i \subseteq \overline{B_j}$ then $B_j$ is also a cutset from $x$ to $\infty$ so that $x \in \overline{B}_j$.
    It follows that $\overline{B}_i \subseteq \overline{B}_j$.
    Similarly, if $B_j \subseteq \overline{B}_i$ then $\overline{B}_j \subseteq \overline{B}_i$.
    Suppose now that $\overline{B}_i \cap \overline{B}_j \neq \varnothing$. 
    Since the blocks are $1$-disconnected we have $d_G(B_i,B_j) \geq 2$.
    If $x \in \overline{B}_i \setminus B_i$ then since $B_i$ is a cutset from $x$ to $\infty$ and since $\partial \overline{B_i} \subseteq \partial B_i$ we have 
    \begin{equation}
        d_G(x,B_j)
        \geq 
        d_G(x,\partial \overline{B}_i)
        +
        d_G(\partial \overline{B}_i, B_j)
        \geq 
        d_G(x,\partial \overline{B}_i)
        +
        d_G(\partial B_i, B_j)
        \geq 2.
    \end{equation}
    If $y \in \overline{B}_j \setminus B_j$ we have similarly that $d_G(B_i,y) \geq 2$ and $d_G(x,y) \geq 2$. 
    It follows that
    \begin{equation}
        d(\overline{B}_i,\overline{B}_j)
        =
        \min_{
            x \in \overline{B}_i \setminus B_i,
            y \in \overline{B}_j \setminus B_j
        }
        \left(
             d_G(B_i,B_j),
             d_G(x,B_j),
             d_G(B_i,y),
             d_G(x,y)
        \right)
        \geq 2
    \end{equation}
    and hence $\overline{B}_i$ and $\overline{B}_j$ are $1$-disconnected. 
    It follows that there is a minimal $\ell \leq b$ such that (up to re-indexing and re-labeling) we can choose a sequence of $1$-disconnected blocks $(B_0,\ldots,B_{\ell - 1})$ from $(\overline{A}_0,\ldots, \overline{A}_{b-1})$ with $\cup_{i \leq b - 1} A_i \subseteq \cup_{j \leq \ell - 1} B_j$ and $\partial_{j \leq l - 1} B_j \subseteq \partial \cup_{i \leq b - 1} A_i$.
\end{proof}

In the following lemmas we show that the boundary of a block is coarsely connected, and similarly that the graph $G$ minus a sequence of blocks is also coarsely connected. 
These statements follow from \Cref{prop:boundary_connectivity}.

\begin{lemma}
    \label[lemma]{lem:block_connectivity}
    Let $G$ be a transitive graph of polynomial growth with $d \geq 2$, and let $R \in \N$ be as in \Cref{prop:boundary_connectivity}.
    For any block $A \in \AA$, the boundary $\partial A$ is $R$-connected.
\end{lemma}

\begin{proof}
    We argue by induction on $n = \# A$ the size of $A$. 
    When $n = 1$, $A = \partial A$ and the boundary is trivially $1$-connected.
    We may assume inductively that the statement holds for all $k \leq  n - 1$.
    Suppose now that $A \in \AA$ with $\#A = n$. 
    Let $S_1\subset \ldots\subset S_n$ be a sequence of nested subsets of $A$ with $S_i = S_{i-1} \cup \{x\}$ for $x \in A \setminus S_{i-1}$ chosen such that $S_i \in \AA$.
    In words, at each step of the sequence we append a vertex such that the new set is a block and such that the sequence exhausts $A$.
    It follows that $\# S_{i} = i$ and $A = S_{n-1} \cup \{x\}$ for some $x \in A$.
    If $x \not\in \partial A$ then $\partial A = \partial S_{n-1}$ and by the induction hypothesis $\partial A$ is $R$-connected. 
    Suppose that $x \in \partial A$.
    It follows that there exists $y \in \partial S_{n-1}$ with $d_G(x,y) = 1$, and either $y \in \partial A$ or $y \not \in \partial A$. 
    In the first case, we have $\partial A = \partial S_{n-1} \cup \{x\}$ with $d_G(v,S_{n-1}) = 1$ and by the induction hypothesis $\partial A$ is $R$-connected.
    In the second case, the set $\partial S_{n-1} \cup \{x\}$ is a cutset from $y$ to $\infty$ and we can choose $W \subset \partial S_{n-1} \cup \{x\}$ the minimal cutset from $y$ to $\infty$.
    By \Cref{prop:boundary_connectivity}, the set $W$ is $R$-connected. 
    Further, since $y \in \partial S_{n-1}$ it must be the case that $x \in W$, and clearly $W \cap \partial S_{n-1} \neq \varnothing$, so that $d_G(x, \partial S_{n-1}) \leq R$.  
    By the induction hypothesis $\partial S_{n-1}$ is $R$-connected, and we have shown that $\partial A = \partial S_{n-1} \cup \{x\}$ with $d_G(x,\partial S_{n-1}) \leq R$ so that $\partial A$ is $R$-connected.
    This concludes the induction.
\end{proof}

\begin{lemma}
    \label[lemma]{lem:complement_r_connected}
    Let $G$ be a transitive graph of polynomial growth with $d \geq 2$, and let $R \in \N$ be as in \Cref{prop:boundary_connectivity}.
    Let $b \in \N$ and let $\ablocks = (A_0,\ldots,A_{b-1}) \in \AA(b)$ be a sequence of $1$-disconnected blocks. Then the set $G \setminus \ablocks$ is $R$-connected.
\end{lemma}

\begin{proof}
    Let $x \in G \setminus (\cup_{i} A_i)$.
    If $\cup_{i} A_i$ is not a cutset from $x$ to $\infty$, then $x$ is in the infinite connected component of $G \setminus (\cup_{i} A_i)$.
    Suppose now that $\cup_{i} A_i$ is a cutset from $x$ to $\infty$. 
 We define $\{x\in A^c, \exists y\in A, d_G(x,y)=1\}$ the external boundary of a block $A$. Note that the external boundaries of $A_0, \dots, A_{b-1}$ all belong to $G\setminus (\cup_{i} A_i)$ as the distance between any $A_i, A_j$ is at least $2$.
    Since by definition blocks are closed, the external boundary of each block is $R$-connected in $G \setminus (\cup_{i} A_i)$ by \Cref{prop:boundary_connectivity}.
    As there are finitely many blocks, it follows that we can find an $R$-connected path from $x$ to a vertex $y \in G \setminus (\cup_{i} A_i)$ for which $\cup_i A_i$ is not a cutset.
    This concludes the proof.
\end{proof}

We now consider the percolation configuration of $\cluster$ the finite cluster of the origin.
Due to the presence of long-range edges and in light of \Cref{lem:uniqueness_and_connectedness,lem:closed_blocks}, the vertex set of $\cluster$ is contained in a unique (up to permutation) sequence of $1$-disconnected blocks.
We call this sequence of $1$-disconnected closed blocks the \textbf{block decomposition} of $K$ and we write $\BB(\cluster)$ for the block decomposition of $K$.
For $b \in \N$, $\mm \in \N^b$, and $\AA(b)$ as in \eqref{eq:disconnected_blocks}, we write
\begin{equation}
    \label{eq:possible_block_decompositions}
    \AA(b,\mm)
    =
    \{
        \ablocks
        =
        (A_0,\ldots,A_{b-1})
        \in
        \AA(b)
        :
        o \in A_0,
        \# \partial A_i = m_i
    \}
\end{equation}
for the set of possible block decompositions of $\cluster$ with $b$ many blocks and fixed boundary sizes.
After revealing the block decomposition of $\cluster$, the vertices on the boundaries of the blocks may not connect in the percolation configuration to vertices in $G \setminus \BB(\cluster)$. 
To make this precise we say that a subset of blocks $A_{i_0},\ldots,A_{i_{\ell - 1}} \subseteq \BB(\cluster)$ is \textbf{$\BB(\cluster)$-isolated} if $\cup_{j \leq \ell - 1} \partial A_{i_j} \not\sim G \setminus \BB(\cluster)$.
The \textbf{block graph} $\HH(\BB(\cluster))$ is obtained from the block decomposition by contracting each block: if $\BB(\cluster) = \ablocks \in \AA(b)$, the block graph $\HH(\ablocks)$ has vertex-set $V = \{0,\ldots,b-1\}$ and edge-set $E = \{(i,j) : A_i \sim A_j\}$. 
For two blocks $A_i,A_j \in \AA$ and $r \in \N$, we write $A_i \jump{} A_j$ for the event 
\begin{equation}
    \label{eq:connection_event_definition}
    \{ 
        \exists x \in A_i,
        \exists y \in A_{j},
        d_G(x,y)=r,
        x \sim y
    \}.
\end{equation}
The \textbf{weighted block-graph} $\HH^* (\BB(\cluster))$ is obtained from the block decomposition by contracting the blocks and keeping track of the length of edges: if $\BB(\cluster) = \ablocks \in \AA(b)$, the weighted block-graph $\HH^* (\ablocks)$ is the weighted multi-graph (without loops) with vertex-set $V = \{0,\ldots,b-1\}$ and edge-set $E = \{(i,j,r) : A_i \jump{} A_j\}$.
Note that two blocks may connect several times.
The block graph and the weighted block graph associated to the percolation configuration in Figure \ref{fig:lrp_configuration} are illustrated in Figures \ref{fig:block_and_spanning_graph} and \ref{fig:weighted_block_graph}.
We summarise this construction in the following lemma.

\begin{lemma}
    \label[lemma]{lem:lattice_animal}
    Let $G$ be a transitive graph of polynomial growth with $d \geq 1$, and let $J : G \times G \to \R_+$ be transitive, and suppose that $J(x,y) = \Theta(d_G(x,y))^{-d \alpha}$ with $\alpha > 1$.
    The finite cluster of the origin $\cluster$ has a unique block decomposition $\BB(\cluster)$, the block graph $\HH(\BB(\cluster))$ is connected, and the block decomposition is $\BB(\cluster)$-isolated. 
    In particular,
    \begin{equation}
        \P_{\beta}
        \left(
            \#K <\infty 
        \right)
        \leq
        \sum_{\substack{
            b \in \N,
            \mm \in \N^{b}
            \\
            \bblocks \in \AA(b,\mm)
        }}
        \P_{\beta}
        \left(
            \BB(\cluster)
            =
            \bblocks,
            \HH(\bblocks)
            \text{ is connected},
            \bblocks
            \text{ is }
            \BB(\cluster)
            \text{-isolated}
        \right).
    \end{equation}
\end{lemma}

\begin{proof}
    Let $A$ be a finite connected subgraph of the complete graph on $G$ containing the origin and suppose that $K = A$. 
    By \Cref{lem:uniqueness_and_connectedness} there exists $\ell \in \N$ and a unique $1$-disconnected sequence of $1$-connected set $(A_0,\ldots,A_{b-1})$ such that $V(A) = A_0 \sqcup \ldots \sqcup A_{\ell-1}$. 
    By \Cref{lem:closed_blocks} there exists $b \in \N$ and a unique $1$-disconnected sequence of blocks $(B_0,\ldots,B_{b - 1})$ such that $\cup_{i \leq \ell - 1} A_i \subset \cup_{j \leq b - 1} B_j$ and $\cup_{j \leq \ell - 1} \partial B_j \subset \cup_{i \leq b - 1} \partial A_i$.
    We have $\bblocks = (B_0,\ldots,B_{b-1}) \in A(b)$ and $\BB(\cluster) = \bblocks$.
    The cluster $K$ is connected and this implies that the block graph $\HH(\bblocks)$ is connected. 
    Further, if $K = A$, all edges from vertices in $A$ to vertices not in $A$ must be closed. 
    By construction vertices in $\cup_{j \leq b - 1} \partial B_j$ are in $K$ and vertices not in $\BB(K)$ are not in $K$, and all such edges must be closed. 
    This is precisely the definition of $\bblocks$ being $\BB(\cluster)$-isolated.
\end{proof}

Note that the block decomposition of $\cluster$ is far from entirely specifying the realisation of $\cluster$. For instance, it is not concerned with the status of edges inside the blocks.

\subsection{Isolation}
In the following lemma we identify a probability cost for the isolation of a sequence of blocks in the block graph in terms of the total boundary size of the blocks.
The argument combines the key properties of the block graph construction (that blocks are closed, the distance between any two blocks is at least two, and boundary vertices are in $K$) with Tim\'ar's \cite{timar_cutsets_2007,timar_boundary-connectivity_2013} notion of coarse connectivity.
The key idea of this argument is sketched in Figure \ref{fig:isolation}.

\begin{lemma}[Isolation]
    \label[lemma]{lem:one_arm_isolation}
    Let $G$ be a transitive graph of polynomial growth with $d \geq 2$, and suppose that $J : V \times V \to \R_+$ is a transitive kernel with $J(x,y) = \Omega(d_G(x,y)^{-d \alpha})$ with $\alpha > 1$.
    Let $b \in \N$, $\mm \in \N^b$, and let $\bblocks = (B_0,\ldots,B_{b-1}) \in \AA(b, \mm)$ be a set of blocks.
    For $\ell \leq b$ and $I = \{i_0,\ldots,i_{\ell - 1}\} \subseteq [0,b-1]$, let $\bblocks' = (B_{i_0},\ldots,B_{i_{\ell - 1}}) \subseteq \bblocks$ be a subset of blocks.
    Then there exists $c = c(G,J) > 0$ such that
    \begin{equation}
        \label{eq:one_arm_isolation}
        \P_{\beta}
        \left(
            \bblocks'
            \text{ is }
            \BB(\cluster)
            \text{-isolated}
            \mid
            \BB(\cluster)
            =
            \bblocks
        \right) 
        \leq
        \prod_{i\in I} \exp(-c\beta m_i).
    \end{equation}
\end{lemma}

\begin{proof}
    Since the blocks $\bblocks = (B_0,\ldots,B_{b-1})$ are $1$-disconnected, and by the definition of the boundary, for each boundary vertex $x \in \cup_{i \in I} \partial B_i$ we can choose an adjacent (in the underlying graph $G$) vertex $y \in G \setminus (\cup_{i} B_i)$ to make the pair $(x,y)$.
    Let $c_J$ and $R_J >0$ be such that $J(x,y) \geq c_J d_G(x,y)^{-d\alpha}$ for all $x,y \in G$ with $d_G(x,y) \geq R_J$.
    By \Cref{lem:complement_r_connected} the set $G \setminus (\cup_{i} B_i)$ is $R$-connected for some $R \in \N$.
    Together with the fact that $G$ is locally finite, it follows that for each pair $(x,y)$ we can find an infinite self-avoiding $R$-connected path $\left(x,y,z_{x,2},z_{x,3}, \ldots \right)$ with $z_{x,n} \in G \setminus (\cup_{i} B_i)$ for every $n \geq 2$ and such that for some $k_x \geq 2$ we have $z_{x,i} \in B(x,R_J + 1)^c$ for every $n \geq k_x$.
    In words, the path travels avoiding the cluster $K$ altogether, and after the first $k_x$ steps the distance between the path and $x$ is at least $R_J + 1$.
    Since the path is $R$-connected, for every $j \geq 0$ we have that $R_J + 1 \leq d_G(x,z_{x,k_x + j}) \leq R_J + R j$ and hence
    \begin{equation}    
        \label{eq:path_kernel_bounds}
        J(x,z_{x, k_x + j})
        \geq 
        c_J
        (R_J + R j)^{-d\alpha}.
    \end{equation}
    Further, for every $j \geq 0$ we have $z_{x,j} \in G \setminus (\cup_{i} B_i)$ and it follows from the definition of $\BB(\cluster)$-isolation that the edge between $x$ and $z_{x,j}$ must be closed.
    After repeating this procedure for all boundary vertices $x \in \cup_{i \in I} \partial B_i$, a given prescribed closed edge is counted exactly once. 
    Indeed, the paths are chosen to be self-avoiding, and prescribed closed edges associated to different paths have distinct endpoints $x\in \cup_{i \in I} \partial B_i$.
    By the bound in \eqref{eq:path_kernel_bounds}, we have shown that
    \begin{multline}
        \P_{\beta}
        \left(
            \bblocks'
            \text{ is }
            \BB(\cluster)
            \text{-isolated}
            \mid
            \BB(\cluster)
            =
            \bblocks
        \right) 
        \leq 
        \prod_{x \in \cup_{i \in I} \partial B_i}
        \prod_{j \geq 0}
        \P_{\beta}
        \left(
            x
            \not\sim
            z_{x,j}
        \right)
        \\
        \leq
        \prod_{x \in \cup_{i \in I} \partial B_i}
        \exp
        \Big(
            -
            c_J
            \beta
            \sum_{j = 0}^{\infty}
            \big(
                R_J
                +
                R
                j
            \big)^{-d \alpha}
        \Big)
        \leq
        \prod_{i \in I}
        \exp
        \left(
            -
            c
            \beta
            m_i
        \right)
   \end{multline}
    for some $c = c(G,J) > 0$, as required.
\end{proof}

\subsection{Peierl's bounds}
\label{sec:integrating}
We use Tim\'ar's \cite{timar_cutsets_2007,timar_boundary-connectivity_2013} notion of coarse connectivity to establish a standard Peierls argument in the setting of transitive graphs of polynomial growth. 
For the sake of completeness, we provide a proof.
Recall that $\AA$ was defined in \eqref{eq:blocks} to be the set of blocks in $G$.

\begin{lemma}[Peierls argument]
\label[lemma]{lem:peierls_argument}
    Let $G$ be a transitive graph of polynomial growth with $d \geq 2$.
    There exists $\cpei = \cpei(G) > 0$ such that
    \begin{equation}
        \label{eq:peierls_argument}
        \#
        \{
            A \in \AA 
            : 
            A \ni x, 
            \# \partial A = m
        \}
        \leq 
        \exp(\cpei m)
    \end{equation}  
    for all $x \in G$ and $m \in \N$.
\end{lemma}

\begin{proof}
    Let $m \in \N$, let $x \in G$, and consider a block $A \in \AA$ containing $x$ with $\partial A = m$. 
    By \Cref{lem:block_connectivity}, there exists $R = R(G) \in \N$ such that the boundary $\partial A$ is $R$-connected. 
    Let $\Gamma$ be the graph with vertex set $G$ and $R$-adjacency. 
    By the volume bounds for graphs of polynomial growth, the degree of $\Gamma$ is at most $\Delta = C_G R^d$. We bound the number of possible such blocks $A$ by counting the number of spanning trees of the induced subgraph of $\partial A$ in $\Gamma$ using a depth-first search algorithm.   
    Fix $y \in \partial A$ the boundary point nearest to $x$, chosen arbitrarily when $y$ is not unique. 
    A walk in $\Gamma$ starting at $y$ travels through each edge at most twice, and at each step there are at most $\Delta - 1$ non-backtracking possibilities. 
    It follows that there are at most $(\Delta-1)^{2m}$ walks starting at $y$ on the induced subgraph of $\partial A$ in $\Gamma$. 
    We now bound the number of potential starting points $y \in \partial A$ for the walk. 
    It follows from the volume bounds for transitive graphs of polynomial growth and the $d$-dimensional isoperimetric inequality that there exists $c_1 = c_1(G) > 0$ such that, independently of the choice of $A$, we have $d_G(x,y) \leq c_1 m^{1/{(d-1)}}$.
    Indeed, since $y$ is chosen to be the boundary point nearest to $x$ we have $B(x,d_G(x,y)) \subseteq A$ and hence
    \begin{equation}   
        \label{eq:ball_inequality}
        c_G (d_G(x,y))^d
        \leq 
        \# B(x,d_G(x,y)) 
        \leq 
        \# A
        \leq 
        \left( 
            m
            /
            \ciso
        \right)^{d/(d-1)}.
    \end{equation}
    Choosing $c_1 = c_1(G) > 0$ proves that $d_G(x,y) \leq c_1 m^{1/(d-1)}$ as required. 
    As a result, there are at most $\# B(x,c_1 m^{1/(d-1)}) \leq C_G (c_1 m^{1/(d-1)})^d$ possible boundary vertices $y \in \partial A$ to begin the walk. 
    It follows that
    \begin{equation}
        \label{eq:peierls_desired_result}
        \#
        \{
            A \in \AA 
            : 
            A \ni x,
            \# 
            \partial 
            A = m
        \}
        \leq 
        C_G (c_1 m^{1/(d-1)})^d
        (\Delta - 1)^{2m}
        \leq 
        \exp
        \left(
            \cpei
            m
        \right)
    \end{equation}
    for some $\cpei = \cpei(G) > 0$, concluding the proof. 
\end{proof}
\section{Decay of the one-arm event}
\label{sec:one_arm}
In this section we prove \Cref{thm:arm_exponent_1} on the decay of the one-arm event. 
In \Cref{lem:one_arm_decomp} we assumed the results on the decay of the distribution of finite clusters to reduce our analysis of the truncated one-arm event to the task of bounding 
$
    \P_{\beta}
    \left(
        o 
        \leftrightarrow 
        B(r)^c, 
        \# 
        \cluster 
        \leq 
        k(r) 
    \right),
$
where $k(r)$ is a polylogarithmic function of $r$.
The first step is to interpret this event in terms of the connectivity of the block graph decomposition of $\cluster$.

\subsection{One-arm event via isolated blocks}
\label{subsec:one_arm_block}
Let $G$ be a transitive graph of polynomial growth, and recall the definition of $\AA(b,\mm)$ the set of possible block decompositions in \eqref{eq:possible_block_decompositions}.
In terms of blocks, the event $\{o \leftrightarrow B(r)^c, \# \cluster \leq k(r)\}$ implies that, for some $b \in \N$ and $\mm \in \N^b$, the block decomposition $\BB(\cluster)$ contains a sequence of blocks $\ablocks = (A_0,\ldots,A_{b-1}) \in \AA(b,\mm)$ with $A_{b-1} \cap B(r)^c \neq \varnothing$, $\#A_{0} + \ldots + \# A_{b-1} \leq k(r)$, and such that the block graph $\HH (\BB(\cluster))$ contains the graph
\begin{equation}
    P'(\ablocks)
    =
    (
        \{0,\ldots,b-1\},
        \{(0,1),\ldots,(b-2,b-1)\}
    )
\end{equation}
as a subgraph.
For $b \in \N$ we let
\begin{equation}
\label{eq:candidate_conditions}
    \MM(b)
    =
    \{
        \mm 
        \in 
        \N^{b} 
        :
        m_0 + \ldots + m_{b-1} \leq k(r)
    \}.
\end{equation}
The maximal distance that can be covered inside each block is the size of its boundary, and the condition $A_{b} \cap B(r)^c \neq \varnothing$ implies that for some $\rr \in \N^{b-1}$ with $m_0 + r_0 + \ldots +  r_{b-2} + m_{b-1} > r$, the weighted block graph $\HH^*(\ablocks)$ contains the weighted graph
\begin{equation}
    \label{eq:candidate_path_graph}
    P(\ablocks,\rr)
    =
    \left(
        \{
            0,
            \ldots,
            b-1
        \},
        \{
            (0,1,r_0),
            \ldots,
            (b-2,b-1,r_{b-1})
        \}
    \right).
\end{equation}  
There may be also additional blocks $\bblocks \setminus \ablocks$ which do not participate in the path event. We may not require that edges between them and the blocks $\ablocks$ constituting the path are closed.
We let
\begin{equation}
    \AA'(\ablocks,\ell)
    =
    \left\{
        \bblocks
        \in
        \AA(\ell)
        :
        \ablocks
        \text{ and }
        \bblocks
        \text{ are }
        1 
        \text{-disconnected}
    \right\}
\end{equation}
be the set of possible completions of the block decomposition.
For our analysis, we will also require that the block $A_{b-1}$ in the path described by $P(\ablocks,\rr)$ is the first to \textit{potentially} reach $B(r)^c$, namely that $m_0 + r_0 + \ldots + r_{b-3} + m_{b-2} \leq r$.
For $\mm \in \MM(b)$ we let
\begin{multline}
    \RR(b,\mm,r)
    =
    \{
        \rr \in \N^{b-1}
        :
        m_0 + r_0 + \ldots + r_{b-3} + m_{b-2} \leq r,
        \\
        m_0 + r_0 + \ldots + r_{b-2} + m_{b-1} > r
    \}.
\end{multline}
We collect this argument into the following lemma, and we verify that the assumption that the final block $A_b$ is the first to potentially reach $B(r)^c$ gives an upper bound.
\begin{proposition}[Block graph event the for the one-arm event]
\label[proposition]{prop:block_one_arm}
    Let $G$ be a transitive graph of polynomial growth, and suppose that $J : V \times V \to \R_+$ is a transitive kernel.  
    Then
    \begin{multline}
        \P_{\beta}
        \left(
            o 
            \leftrightarrow 
            B(r)^c, 
            \# \cluster
            \leq 
            k(r)
        \right) 
        \\
        \leq 
        \sum_{\substack{
            b \in \N
            \\
            \mm \in \MM(b)
            \\
            \rr \in \RR(b,\mm,r)
            \\
            \ablocks \in \AA(b,\mm)
        }}
        \sum_{\substack{
            \ell \in \N
            \\
            \bblocks \in \AA'(\ell,\ablocks)
        }}
        \P_{\beta}
        \left(
            \BB(K)
            =
            \ablocks \cup \bblocks,
            P(\ablocks, \rr) 
            \subseteq 
            \HH^*(\ablocks),
            \ablocks 
            \text{ is } 
            \BB(\cluster) 
            \text{-isolated}
        \right).
    \end{multline}
\end{proposition}

\begin{proof}
    For $b \in \N$ and $\mm \in \MM(b)$, let 
    \begin{equation}
        \label{eq:alt_candidate_vectors}
        \RR'(b,\mm,r)
        =
        \{
            \rr \in \N^{b-1}
            :
            m_0 + r_0 + \ldots + r_{b-2} + m_{b-1} > r
        \}.
    \end{equation}
    Note that for $\rr \in \RR'(b,\mm,r)$, $P(\ablocks,\rr)$ is a candidate one-arm path: the endpoint of the path may be in $B(r)$, but the total length of the path is at least $r$.
    It follows from \Cref{lem:lattice_animal} and from the discussion above that 
	\begin{multline}
        \P_{\beta}
        \left(
            o 
            \leftrightarrow 
            B(r)^c, 
            \# \cluster
            \leq 
            k(r)
        \right) 
        \\
        \leq 
        \sum_{\substack{
            b \in \N
            \\
            \mm \in \MM(b)
            \\
            \rr \in \RR'(b,\mm,r)
            \\
            \ablocks \in \AA(b,\mm)
        }}
        \sum_{\substack{
            \ell \in \N
            \\
            \bblocks \in \AA'(\ell,\ablocks)
        }}
        \P_{\beta}
        \left(
            \BB(K)
            =
            \ablocks \cup \bblocks,
            P(\ablocks, \rr) 
            \subseteq 
            \HH^*(\ablocks),
            \ablocks 
            \text{ is } 
            \BB(\cluster) 
            \text{-isolated}
        \right).
    \end{multline}
    We also want to require that the block $A_{b-1}$ in the path described by $P(\ablocks,\rr)$ is the first to \emph{potentially} reach $B(r)^c$, namely that $m_0 + r_0 + \ldots + r_{b-3} + m_{b-2} \leq r$.
    Suppose that $\rr \in \RR'(b,\mm,r)$ satisfies $m_0 + r_0 + \ldots + r_{b'-2} + m_{b'-1} > r$ for some $b' < b$ taken to be minimal.
    This means that the blocks $A_0,\ldots,A_{b'-1}$ already constitute a candidate one-arm path, and we want to put lump the blocks $A_{b'},\ldots,A_{b-1}$ into the sequence of extraneous blocks $\bblocks$.
    If we let $\rr' = (r_0,\ldots,r_{b'-2})$, $\mm' = (m_0,\ldots,m_{b'-1})$, and $\ablocks' = (A_0,\ldots,A_{b'-1})$, we have
    \begin{multline}
    \label{eq:path_containment_1}
        \P_{\beta}
        \left(
            \BB(K)
            =
            \ablocks
            \cup
            \bblocks,
            P(\ablocks, \rr) 
            \subseteq 
            \HH^*(\ablocks),    
            \ablocks 
            \text{ is }
            \BB(\cluster)
            \text{-isolated}
        \right)
        \\
        \leq
        \P_{\beta}
        \left(
            \BB(\cluster)
            = 
            \ablocks'
            \cup 
            \left(
                \bblocks
                \cup 
                \left(
                    \ablocks 
                    \setminus
                    \ablocks'
                \right)
            \right),
            P(\ablocks',\rr;) 
            \subseteq 
            \HH^*(\ablocks'),
            \ablocks'
            \text{ is }
            \BB(\cluster)
            \text{-isolated}
        \right).
    \end{multline}
    This holds since $\{P(\ablocks',\rr') \subseteq \HH^*(\ablocks')\}$ requires fewer open edges than $\{P(\ablocks, \rr) \subseteq \HH^*(\ablocks)\}$, and $\{\ablocks' \text{ is } \BB(\cluster) \text{-isolated}\}$ requires fewer closed edges than $\{\ablocks \text{ is } \BB(\cluster) \text{-isolated}\}$. 
    The statement follows.
\end{proof}

By conditioning on the block decomposition, it follows from the independence of edges in \LRP{} that the one-arm candidate event can be separated into independent connection and isolation events.

\begin{lemma}
    \label[lemma]{lem:independence}
    Let $G$ be a transitive graph of polynomial growth with $d \geq 2$, and suppose that $J : V \times V \to \R_+$ is a transitive kernel with $J(x,y) = \Omega(d_G(x,y)^{-d \alpha})$ with $\alpha > 1$.  
    Then there exists $c = c(G,J) > 0$ such that
    \begin{multline}
        \label{eq:connection_isolation_independence}
        \P_{\beta}
        \left(
            o 
            \leftrightarrow 
            B(r)^c, 
            \# \cluster
            \leq 
            k(r)
        \right) 
        \\
        \leq 
        \sum_{\substack{
            b \in \N
            \\
            \mm \in \MM(b)
        }}  
        \prod_{i = 0}^{b-1}
        \exp
        \left(
            -c \beta m_i 
        \right)
        \sum_{\substack{
            \rr \in \RR(b,\mm,r)
            \\
            \ablocks \in \AA(b,\mm)
        }}
        \P_{\beta}
        \left(
            P(\ablocks,\rr)
            \subseteq
            \HH^*(\ablocks)
        \right).
    \end{multline}
\end{lemma}

\begin{proof}
    Consider the bound in \Cref{prop:block_one_arm} and condition on the event $\{\BB(\cluster) = \ablocks \cup \bblocks\}$.
    The connection event $\{P(\ablocks,\rr) \subseteq \HH^*(\ablocks)\}$ is independent of the conditioning, while the isolation event $\{\ablocks \text{ is } \BB(\cluster) \text{-isolated}\}$ depends on the conditioning but is conditionally independent of the connection event, so that
    \begin{multline}
            \P_{\beta}
            \left(
                o 
                \leftrightarrow 
                B(r)^c, 
                \# \cluster
                \leq 
                k(r)
            \right) 
            \leq 
            \label{eq:independence}
            \sum_{\substack{
                b \in \N,
                \mm \in \MM(b)
                \\
                \rr \in \RR(b,\mm,r),
                \ablocks \in \AA(b,\mm)
            }}
            \P_{\beta}
            \left(
                P(\ablocks,\rr)
                \subseteq
                \HH^*(\ablocks)
            \right)
            \\
            \times
            \sum_{\substack{
                \ell \in \N
                \\
                \bblocks \in \AA'(\ablocks,\ell)
            }}
            \P_{\beta}
            \left(
                \ablocks
                \text{ is }
                \BB(\cluster)
                \text{-isolated}
                \mid
                \BB(\cluster)
                =
                \ablocks
                \cup
                \bblocks
            \right)
            \P_{\beta}
            \left(
                \BB(\cluster)
                =
                \ablocks
                \cup
                \bblocks
            \right).
    \end{multline}
    Fix $b \in \N, \mm \in \MM(b)$, and $\ablocks \in \AA(b,\mm)$.
    \Cref{lem:one_arm_isolation} gives a uniform bound for the isolation of the blocks $\ablocks$ in the potential one-arm path in terms of their total boundary size, and hence
    \begin{multline}
        \sum_{\substack{
            \ell \in \N
            \\
            \bblocks \in \AA'(\ablocks,\ell)
        }}
        \P_{\beta}
        (
            \ablocks
            \text{ is }
            \BB(\cluster)
            \text{-isolated}
            \mid
            \BB(\cluster)
            =
            \ablocks
            \cup
            \bblocks
        ) 
        \P_{\beta}
        \left(
            \BB(\cluster)
            =
            \ablocks
            \cup
            \bblocks
        \right)
        \\
        \leq
        \prod_{i = 0}^{b-1}
        \exp
        \left(
            -c \beta m_i 
        \right)
        \sum_{\substack{
            \ell \in \N
            \\
            \bblocks \in \AA'(\ablocks,\ell)
        }}
        \P_{\beta}
        \left(
            \BB(\cluster)
            =
            \ablocks
            \cup
            \bblocks
        \right)
        \leq
        \prod_{i = 0}^{b-1}
        \exp
        \left(
            -c \beta m_i 
        \right)
    \end{multline}
    for some $c = c(G,J) > 0$, where we use that the second sum equals $\P_{\beta}(\BB(\cluster) \supseteq \ablocks)$, which is at most one.
    The result follows.
\end{proof}

\subsection{Uniform bounds on long-edge integrals}
We now consider the connectivity event in \eqref{eq:connection_isolation_independence}, namely the probability
$
    \P_{\beta}
    \left(
        P(\ablocks,\rr)
        \subseteq
        \HH^*(\ablocks)
    \right)
$
of a candidate one-arm path. 
The goal of this section is to prove the following proposition giving uniform bounds on the probability of long edges. 
Note that the main decay term $r^{d(1-\alpha)}$ already appears there.
Recall that $r$ is the length of the one-arm path, $b$ is the number of blocks participating in the candidate one-arm path, and $\mm \in \MM(b)$ is the vector fixing the boundary sizes of the blocks. 

\begin{lemma}
    \label[lemma]{prop:one_arm_connection}
    Let $G$ be a transitive graph of polynomial growth with $d \geq 2$, and suppose that $J : V \times V \to \R_+$ is a transitive kernel with $J(x,y) = O(d_G(x,y)^{-d \alpha})$ with $\alpha > 1$.  
    For $r \in \N$, $b \geq 2$, and $\mm \in \MM(b)$, there exist $c_1 = c_1(G,J) > 0$ and $c_2 = c_2(G) > 0$ such that 
    \begin{multline}        
        \sum_{\substack{
            \rr \in \RR(b,\mm,r)
            \\
            \ablocks \in \AA(b,\mm)
        }}
        \P_{\beta}
        \left(
            P(\ablocks,\rr)
            \subseteq
            \HH^*(\ablocks)
        \right)
        \\
        \leq 
        c_1 
        r^{d(1 - \alpha)}
        (b-2)^{d \alpha - d + 1}
        \expdeg{o}^{b-1}
        \prod_{i=0}^{b-1} 
        \exp
        \left(
            c_2 \beta m_i
        \right)
        .
    \end{multline}
\end{lemma}

We prove this result in two steps. 
First, we inductively apply the isoperimetric inequality and the Peierls argument from \Cref{lem:trivial_isoperimetry,lem:peierls_argument} respectively.
Then, we use the pigeonhole principle to obtain uniform bounds on the probability of long edges.

\begin{lemma}[Inductive step]
\label[lemma]{lem:sum_over_blocks}
    Let $G$ be a transitive graph of polynomial growth with $d \geq 2$, and suppose that $J : V \times V \to \R_+$ is a transitive kernel with $J(x,y) = O(d_G(x,y)^{-d \alpha})$ with $\alpha > 1$.  
    For $r \in \N$, $b \geq 2$, $\mm \in \MM(b)$, and $\rr \in \RR(b,\mm,r)$, there exists $c = c(G) > 0$ such that
    \begin{equation}
    \label{eq:sum_over_blocks}
        \sum_{\ablocks \in \AA(b,\mm)}
        \P_{\beta}
        \left(
            P(\ablocks,\rr)
            \subseteq
            \HH^*(\ablocks)
        \right)
        \leq
        \prod_{i = 0}^{b-1}
        \exp
        \left(
            c 
            m_i
        \right)
        \prod_{\substack{j = 0}}^{b-2}
        \expdegrad{o}{r_j}.
    \end{equation}
\end{lemma}

\begin{proof}
    We argue by induction on $b$. 
    Recall the definition of $P(\ablocks,\rr)$ in \eqref{eq:candidate_path_graph}.
    When $b = 2$, $\mm \in \MM(2)$, and $\rr \in \RR(2,\mm)$, the path consists of two blocks with a single long edge between them and
    \begin{equation}
        \sum_{\ablocks \in \AA(2,\mm)}
        \P_{\beta}
        \left( 
            P(\ablocks,\rr)
            \subseteq 
            \HH^*(\ablocks)
        \right)
        \leq
        \sum_{\substack{
            A_0 \in \AA, 
            A_0 \ni o
            \\
            \# \partial A_0 = m_0}
        }
        \sum_{y \not\in A_0}
        \sum_{\substack{
            A_1 \in \AA,
            A_1 \ni y
            \\
            \# \partial A_1 = m_1}
        }
        \P_{\beta}
        \left(
            A_0 
            \jump{0}
            y
        \right).
    \end{equation}
    By the Peierls argument in \Cref{lem:peierls_argument} we have
    \begin{align}
        \sum_{\substack{
            A_0 \in \AA,
            A_0 \ni o 
            \\
            \# \partial A_0 = m_0}
        }
        \sum_{y \not\in A_0}
        \sum_{\substack{
            A_1 \in \AA,
            A_1 \ni y
            \\
            \# \partial A_1 = m_1}
        }
        \P_{\beta}
        \left(
            A_0 
            \jump{0}
            y
        \right)
        & =
        \sum_{\substack{
            A_0 \in \AA,
            A_0 \ni o 
            \\
            \# \partial A_0 = m_0}
        }
        \sum_{y \not\in A_0}
        \P_{\beta}
        \left(
            A_0
            \jump{0}
            y
        \right)
        \sum_{\substack{
            A_1 \in \AA,
            A_1 \ni y
            \\
            \# \partial A_1 = m_1}
        }
        1
        \\
        & \leq 
        \exp
        \left(
            \cpei
            m_1
        \right)
        \sum_{\substack{
            A_0 \in \AA,
            A_0 \ni o 
            \\
            \# \partial A_0 = m_0}
        }
        \sum_{y \not\in A_0}
        \P_{\beta}
        \left(
            A_0
            \jump{0}
            y
        \right).
    \end{align}
    Recall that we defined the event $\{A_0 \jump{0} y\}$ to be $\{\exists x \in A_0, d_G(x,y)=r, x \sim y\}$.
    By the $d$-dimensional isoperimetric inequality \eqref{eq:isoperimetric_inequality} and by the transitivity of the kernel $J$,
    \begin{align}
        \sum_{\substack{
            A_0 \in \AA,
            A_0 \ni o 
            \\
            \# \partial A_0 = m_0}
        }
        \sum_{y \not\in A_0}
        \P_{\beta}
        \left(
            A_0
            \jump{0}
            y
        \right)
        & = 
        \sum_{\substack{
            A_0 \in \AA,
            A_0 \ni o 
            \\
            \# \partial A_0 = m_0}
        }
        \sum_{\substack{
            x \in A_0
            \\
            y \in S(x,r_0) \setminus A_0
        }}
        \P_{\beta}
        \left(
            x
            \sim
            y
        \right)
        \\
        & \leq 
        \ciso{}
        m_0^{d/(d-1)}
        \sum_{\substack{
            A_0 \in \AA,
            A_0 \ni o 
            \\
            \# \partial A_0 = m_0
        }}
        \sum_{y \in S(r_0)}
        \P_{\beta}
        \left(
            o
            \sim
            y
        \right)
        \\
        & \leq 
        \ciso{}
        \exp
        \left(
            \frac{
                m_0 d 
            }{
                d-1
            }
        \right)
        \sum_{\substack{
            A_0 \in \AA,
            A_0 \ni o 
            \\
            \# \partial A_0 = m_0
        }}
        \sum_{y \in S(r_0)}
        \P_{\beta}
        \left(
            o
            \sim
            y
        \right)
    \end{align}
    where in the last inequality we use the fact that $m^{c} \leq \exp(m c)$ for $m \in \N$ and $c > 1$.
    Recall the definition of $\expdegrad{o}{L}$ the expected degree of the origin with edges of length $L$ in \eqref{eq:def_expdegrad}. Note that $\sum_{y \in S(r_0)} \P_{\beta} \left(o \sim y \right) = \expdegrad{o}{r_0}$, and applying the Peierls argument a second time,
    \begin{align}
        \sum_{\substack{
            A_0 \in \AA,
            A_0 \ni o 
            \\
            \# \partial A_0 = m_0
        }}
        \sum_{y \in S(r_0)}
        \P_{\beta}
        \left(
            o
            \sim
            y
        \right)
        & =
        \expdegrad{o}{r_0}
        \sum_{\substack{
            A_0 \in \AA,
            A_0 \ni o 
            \\
            \# \partial A_0 = m_0
        }}
        1
        \\
        & \leq
        \exp
        \left(
            \cpei
            m_0
        \right)
        \expdegrad{o}{r_0}.
    \end{align}
    Letting $c = \cpei + d/(d-1) + \ciso{}$, we have shown that
    \begin{equation}
         \sum_{\ablocks \in \AA(2,\mm)}
        \P_{\beta}
        \left( 
            P(\ablocks,\rr)
            \subseteq 
            \HH^*(\ablocks)
        \right)
        \leq
        \exp
        \left(
            c
            \left(
                m_0
                + 
                m_1
            \right)
        \right)
        \expdegrad{o}{r_0},
    \end{equation}
    which concludes the base case of the induction.
    We may assume inductively that the statement holds for $k \leq b-1$. 
    Let $\mm \in \MM(b)$ and $\rr \in \RR(b,\mm,r)$.
    For a fixed block $A_{b-1} \in \AA$ with $\# \partial A_{b-1} = m_{b-1}$, the same calculations as in the base give 
    \begin{equation}
        \sum_{y \not \in A_{b-1}}
        \sum_{\substack{
            A_{b-1} \in \AA,
            A_{b-1} \ni y
            \\
            \# \partial A_{b-1} = m_{b-1}
            }
        }
        \P_{\beta}
        \left(
            A_{b-1} \jump{b-2} y
        \right)
        \leq
        \exp
        \left(
            c 
            m_{b-1}
        \right)
        \expdegrad{o}{r_{b-2}}.
    \end{equation}
    We write $\widehat{\mm} = (m_0,\ldots,m_{b-2})$ and $\widehat{\rr} = (r_0,\ldots,r_{b-3})$.
    By the independence of edges in \LRP{} and together with the induction hypothesis, we have shown that
    \begin{align}
        \begin{split}
            &
            \sum_{\ablocks \in \AA(b,\mm)}
            \P_{\beta}
            \left(
                P(\ablocks,\rr)
                \subseteq
                \HH^*(\ablocks)
            \right)
            \\
            & 
            \quad
            =
            \sum_{\ablocks \in \AA(b-1,\widehat{\mm})}
            \P_{\beta}
            \left(
                P(\ablocks,\widehat{\rr})
                \subseteq 
                \HH^*(\ablocks)
            \right)
            \sum_{y \not \in A_{b-2}}
            \sum_{\substack{
                A_{b-1} \in \AA,
                A_{b-1} \ni y
                \\
                \# \partial A_{b-1} = m_{b-1}}
            }
            \P_{\beta}
            \left(
                A_{b-1}
                \jump{b-2}
                y
            \right)
        \end{split}
        \\
        & 
        \quad
        \leq
        \label{eq:base_case_1}
        \prod_{i = 0}^{b-1}
        \exp
        \left(
            c 
            m_i
        \right)
        \prod_{\substack{j = 0}}^{b-2}
        \expdegrad{o}{r_j},
    \end{align}
    concluding the induction.
\end{proof}

Note that by the asymptotics on $J$ and the kernel bounds in \eqref{eq:kernel_bounds} we have
$
    \expdegrad{o}{r}
    =
    O
    \left(
        \beta
        \#       
        S(r)
        r^{-d \alpha}
    \right)
$,
where $S(r)$ denotes the sphere of radius $r$.
Finding sharp bounds for the volume of spheres in transitive graphs of polynomial growth, or even in Cayley graphs of groups of polynomial growth, remains an open problem (see the discussion in Section \ref{subsubsec:spheres_in_graphs}).
Instead of using the highly refined estimates in \cite{colding_liouville_1998,tessera_volume_2007} or \cite{breuillard_rate_2013} to immediately estimate the volume of spheres in \eqref{eq:base_case_1}, we carry forward the term $\expdegrad{o}{r}$ to the next sum and we use the ``spheres partition balls'' argument recorded in \Cref{lem:exp_deg_geq_r} to obtain uniform bounds on the probability of long edges.  

\begin{lemma}[Uniform bounds on long edges]
\label[lemma]{lem:unif_bounds_long_edges}
    Let $G$ be a transitive graph of polynomial growth with $d \geq 2$, and suppose that $J : V \times V \to \R_+$ is a transitive kernel with $J(x,y) = O(d_G(x,y)^{-d \alpha})$ with $\alpha > 1$.    
    For $r \in \N$, $b \geq 2$, and $\mm \in \MM(b)$, there exists $c = c(G,J) > 0$ such that for all $r$ sufficiently large,
    \begin{equation}
        \sum_{\rr \in \RR(b,\mm,r)}
        \prod_{i = 0}^{b-2}
        \expdegrad{o}{r_i}
        \leq 
        c
        r^{d(1-\alpha)}
        (b-2)^{d\alpha - d + 1}
        \expdeg{o}^{b-1}
        .
    \end{equation}
\end{lemma}

\begin{proof}
    Let $b \geq 2$ and $\mm \in \MM(b)$.
    Since $\rr \in \RR(b,\mm,r)$, we have $r_0 + \ldots + r_{b-2} > r - k(r)$ and hence $r_{b-2} \geq r - k(r) - r_0 - \ldots - r_{b-3} + 1$.
    We also have that $r_0 + \ldots + r_{b-3} < r - k(r) + 1$.
    Let us write $L = r - k(r) - r_0 - \ldots - r_{b-3} + 1 > 0$ for the minimum length of the final jump $r_{b-2}$.
    In the inner sum over $r_{b-2}$ below we use \Cref{lem:exp_deg_geq_r}, yielding
    \begin{align}
        \begin{split}
            \sum_{\rr \in \RR(b,\mm,r)}
            &
            \prod_{i = 0}^{b-2}
            \expdegrad{o}{r_i}
            \\
            & \leq
            \sum_{r_0 = 2}^{\infty}
            \ldots 
            \sum_{r_{b-3} = 2}^{\infty}
            \prod_{i = 0}^{b-3}
            \expdegrad{o}{r_i}
            \sum_{r_{b-2} = L}^{\infty}
            \expdegrad{o}{r_{b-2}}
        \end{split}
        \\
        & \leq
        c_1
        \beta
        \sum_{r_0 = 2}^{\infty}
        \ldots 
        \sum_{r_{b-3} = 2}^{\infty}
        L^{d(1 - \alpha)}
        \prod_{i = 0}^{b-3}
        \expdegrad{o}{r_i}
    \end{align}
    for some $c_1(G,J,\alpha) > 0$.
    Plotted over $b-2$ dimensions, the inequality $r_0 + \ldots + r_{b-3} \le  r - k(r) + 1$ draws a simplex.
    Fix some $0< \gamma < 1$ and consider two regions of this sum, the `triangle' region $T$ and the `slab' region $S$ where
    \begin{gather}
        T 
        =
        \left\{
            \widehat{\rr} \in \N^{b-2}
            :
            r_0 
            + 
            \ldots 
            + 
            r_{b-3} 
            \leq 
            \gamma (r - k(r) + 1)
        \right\},
        \\
        S 
        = 
        \left\{
            \widehat{\rr} \in \N^{b-2}
            :
            r_0 
            + 
            \ldots 
            + 
            r_{b-3}
            > 
            \gamma (r - k(r) + 1)
        \right\}.
    \end{gather}
    If we consider the sum over the triangle region, since $d(1-\alpha) < 0$ we have $L^{d(1-\alpha)} \leq ((1 - \gamma)(r - k(r) + 1))^{d(1-\alpha)}$ and again using \Cref{lem:exp_deg_geq_r} we have
    \begin{align}
        \sum_{\widehat{\rr} \in T}
        L^{d(1 - \alpha)}
        \prod_{i = 0}^{b-3}
        \expdegrad{o}{r_i}
        & \leq
        (
            (1 - \gamma)
            (r - k(r) + 1)
        )^{d(1-\alpha)}
        \sum_{\widehat{\rr} \in T}
        \prod_{i = 0}^{b-3}
        \expdegrad{o}{r_i}
        \\
        & \leq 
        (
            (1 - \gamma)
            (r - k(r) + 1)
        )^{d(1-\alpha)}
        \expdeg{o}^{b-2}.
    \end{align}
    We now bound the sum over the slab region $S$.
    By the pigeonhole principle at least one of the $r_i$ satisfies $r_i > \gamma(r - k(r) + 1)/(b-2)$, and we upper-bound the sum over $S$ by summing over the rectangular regions where one variable $r_i$ satisfies this lower bound and the remaining variables are unbounded. 
    Note that in the slab region $S$ we have $L^{d(1-\alpha)} \leq 1$, and again using \Cref{lem:exp_deg_geq_r} we have
    \begin{align}
        \begin{split}
            \sum_{\widehat{\rr} \in S}
            L^{d(1 - \alpha)}
            &
            \prod_{i = 0}^{b-3}
            \expdegrad{o}{r_i}
            \\
            & \leq 
            \sum_{i = 0}^{b - 3}
            \sum_{r_i > \gamma(r - k(r) + 1)/(b-2)}
            \expdegrad{o}{r_i}
            \prod_{j \neq i}
            \sum_{r_j = 2}^{\infty}
            \expdegrad{o}{r_j}
        \end{split}
        \\
        & \leq 
        c_2 
        \gamma^{d(1-\alpha)}
        (r-k(r)+1)^{d(1-\alpha)}
        \beta
        (b-2)^{d \alpha -d + 1}
        \expdeg{o}^{b-3}
        \\
        & \leq
        c_3 
        (r-k(r)+1)^{d(1-\alpha)}
        (b-2)^{d \alpha -d + 1}
        \expdeg{o}^{b-2}
    \end{align}
    for some $c_2(G,J) > 0$ and $c_3(G,J) > 0$.
    Piecing together the bounds for the triangle and slab regions, we have shown that
    \begin{equation}
        \sum_{\rr \in \RR(b,\mm,r)}
        \prod_{i = 0}^{b-2}
        \expdegrad{o}{r_i}
        \leq 
        c_4
        (r-k(r)+1)^{d(1-\alpha)}
        (b-2)^{d\alpha - d + 1}
        \expdeg{o}^{b-1}
    \end{equation}
    for some $c_4(G,J,a) > 0$.
    Choosing $k(r)$ as in \Cref{lem:one_arm_decomp} we have that $(r - k(r) + 1)^{d(1-\alpha)} \leq r^{d(1-\alpha)}(1-\epsilon)$ for all $r$ sufficiently large, and hence
    \begin{equation}
        \sum_{\rr \in \RR(b,\mm,r)}
        \prod_{i = 0}^{b-2}
        \expdegrad{o}{r_i}
        \leq 
        c_5
        r^{d(1-\alpha)}
        (b-2)^{d\alpha - d + 1}
        \expdeg{o}^{b-1}
    \end{equation}
    for some $c_5(G,J) > 0$ as required.
\end{proof}

\begin{proof}[Proof of \Cref{prop:one_arm_connection}]
    The statement follows directly from \Cref{lem:sum_over_blocks,lem:unif_bounds_long_edges}.
\end{proof}

\subsection{Truncated one-arm event:  upper bound}
\begin{proof}[Proof of \Cref{thm:arm_exponent_1}]
    By \Cref{lem:independence,prop:one_arm_connection}, there exist $c_1 = c_1(G,J),c_2 = c_2(G,J),c_3 = c_3(G) > 0$ such that
    \begin{multline}   
    \label{eq:all_together_now}
        \P_{\beta}
        \left(
            o \leftrightarrow B(r)^c, 
            \# \cluster \leq k(r)
        \right)
        \leq
        c_1
        r^{d(1-\alpha)}
        \sum_{b \in \N}
        \Bigg(
            \expdeg{o}^{b-1}
            (b-2)^{d \alpha - d + 1}
            \\
            \times
            \sum_{\mm \in \MM(b)}
            \prod_{i = 0}^{b-1}
            \exp
            \left(
                m_i
                \left(
                    - c_2 
                    \beta
                    + 
                    c_3
                \right)
            \right)
        \Bigg)
    \end{multline} 
    for all $r$ sufficiently large.
    We bound the two sums in turns. 
    Fixing $b \in \N$, we have
    \begin{equation}
        \label{eq:inner_sum}
        \sum_{\mm \in \MM(b)}
        \prod_{i = 0}^{b-1}
        \exp
        \left(
            m_i
            \left(
                - c_2 
                \beta
                + 
                c_3
            \right)
        \right)
        \leq
        \left(
            \sum_{m = 1}^{\infty} 
            \exp
            \left(
                - 
                c_2
                \beta
                +
                c_3
            \right)^m
        \right)^b.
    \end{equation}
    The sum over $m$ is a geometric series with ratio $q_1(\beta) = \exp(-c_2 \beta + c_3)$, and we can choose $\beta$ sufficiently large such that $0 < q_1(\beta) < 1/2 < 1$.
    For this choice of $\beta$, the sum over $b \in \N$ in \eqref{eq:all_together_now} is bounded above by
    \begin{equation}
    \label{eq:arithmetico_geometric_series}
        \sum_{b \in \N}
        (b-2)^{d \alpha - d + 1}
        \left(
            \frac{
                \expdeg{o}
                q_1(\beta)
            }{
                1 - q_1(\beta)
            }
        \right)^{b-1}.
    \end{equation}
    By \Cref{lem:finite_expdeg} we have $\E_{\beta} \left[ \deg(o) \right] \leq c_4 \beta$ for some $c_4 = c_4(G,J,\alpha) > 0$, and the geometric part of the series has ratio $q_2(\beta) = c_4 \beta q_1(\beta) / (1 - q_1(\beta))$.
    Since $q_1(\beta)$ decays exponentially in $\beta$, we may once again choose $\beta$ sufficiently large such that $0 < q_2(\beta) < 1$. 
    For this choice of $\beta$ the series over $m \in \N$ converges and hence 
    \begin{equation}
        \label{eq:desired_bound}
        \P_{\beta}
        \left(
            o 
            \leftrightarrow 
            B(r)^c, 
            \# \cluster
            \leq 
            k(r)
        \right)
        \leq
        c_5
        r^{d(1-\alpha)}
    \end{equation}
    for some $c_5 = c_5(\beta,G,J) > 0$.
    Together with \Cref{lem:one_arm_decomp} we have shown that there exists $c_6 = c_6(\beta,G,J) > 0$ such that 
    \begin{equation}
        \label{eq:one_arm_desired_bound}
        \P_{\beta}
        \left(
            o 
            \leftrightarrow 
            B(r)^c, 
            \# 
            \cluster 
            < 
            \infty 
        \right) 
        \leq
        c_6
        r^{d(1-\alpha)}
    \end{equation}
    for all $r \in \N$, concluding the proof.
\end{proof} 
\section{Cluster-size decay for \texorpdfstring{$\alpha > 1 + 1/d$}{alpha > 1 + 1/d}}
\label{sec:cluster_decay_above}
In this section we prove \Cref{thm:cluster_size_strong_decay}.
The layout of the proof is similar to that of the proof of \Cref{thm:arm_exponent_1}, we inductively apply an isoperimetric inequality and the Peierl's argument, but with a different defining connectivity event.
For the distribution of finite clusters, we will require that the block graph is connected, and we bound this probability by counting spanning trees.
It will also be necessary to use the isoperimetric inequality given by \Cref{lem:want_to_show} instead of \Cref{lem:trivial_isoperimetry}.
This will explain the transition in the driving exponent of the decay of the distribution of finite clusters, discussed in detail in \Cref{rem:crucial_alpha_assumption}.
As before, the first step is to express $\{k \leq \# \cluster < \infty\}$ in terms of the block graph decomposition of $K$.

\subsection{Cluster-size decay event via isolated blocks}
In terms of blocks, the event $\{k \leq \# \cluster \leq \infty\}$ implies that the block decomposition $\BB(\cluster)$ consists of a sequence of blocks with total size greater than $k$ and such that the block graph $\HH(\BB(\cluster))$ is connected.
For $b \in \N$ and $\mm \in \N^b$, recall the definition of $\AA(b,\mm)$ the set of possible block decompositions in \eqref{eq:possible_block_decompositions}.
Suppose that $\BB(\cluster) = \bblocks = (B_0,\ldots,B_{b-1}) \in \AA(b,\mm)$ with $\# \cup_i B_i \geq k$.
It follows from the $d$-dimensional isoperimetric inequality that
\begin{equation}
    \label{eq:concavity_argument}
    m_0 
    + 
    \ldots 
    + 
    m_{b-1}
    \geq 
    \ciso
    \left(
        \# B_0^{(d-1)/d} 
        + 
        \ldots
        +
        \# B_{b-1}^{(d-1)/d} 
    \right)
    \geq 
    \ciso k^{(d-1)/d},
\end{equation}
where we use that the function $f(x) = x^{(d-1)/d}$ is concave and increasing.
We let
\begin{equation}
    \label{eq:cluster_size_blocks}
    \MM(b,k)
    =
    \{ 
        \mm 
        \in 
        \N^{b}
        : 
        m_0 
        + 
        \ldots
        +
        m_{b-1}
        \geq
        \ciso k^{(d-1)/d}
    \}.
\end{equation}  

\begin{proposition}[Block graph event for the cluster-size decay event]
\label[proposition]{prop:block_cluster}
    Let $G$ be a transitive graph of polynomial growth, and suppose that $J : V \times V \to \R_+$ is a transitive kernel. 
    Then for all $k \in \N$,
    \begin{multline}    
        \label{eq:cluster_size_upper_bound}
        \P_{\beta}
        \left(
            k
            \leq 
            \#
            \cluster
            < \infty
        \right)
        \leq
        \sum_{\substack{
            b \in \N
            \\
            \mm \in \MM(b,k)
        }}
        \frac{1}{(b-1)!}
        \\
        \times
        \sum_{\bblocks \in \AA(b,\mm)}
        \P_{\beta}
        \left(
            \BB(\cluster)
            = 
            \bblocks,
            \HH(\bblocks)
            \text{ is connected},
            \bblocks
            \text{ is }
            \BB(\cluster)
            \text{-isolated}
        \right).
    \end{multline}
\end{proposition}

\begin{proof}
    This follows from Lemma \ref{lem:lattice_animal} with the additional requirement that the total size of the blocks is at least $k$.
    When summing over $\bblocks \in \AA(b,\mm)$, we sum over the location of the blocks in the block decomposition and also over the labelling of the blocks. 
    The label of the block $B_0$ is fixed since $o \in B_0$, and hence a given block decomposition is counted $(b-1)!$ times.
    The factor $1/(b-1)!$ accounts for the permutations of the labels of the blocks. 
\end{proof}

\begin{lemma}
    \label[lemma]{prop:cluster_size_independence}
    Let $G$ be a transitive graph of polynomial growth with $d \geq 2$, and suppose that $J : V \times V \to \R_+$ is a transitive kernel with $J(x,y) = \Omega(d_G(x,y)^{-d \alpha})$ with $\alpha > 1$. 
    Then there exists $c = c(G,J) > 0$ such that for all $k \in \N$
    \begin{multline}
        \label{eq:cluster_size_independence}
        \P_{\beta}
        \left(
            k
            \leq 
            \#
            \cluster
            < \infty
        \right)
        \\
        \leq
        \sum_{\substack{
            b \in \N
            \\
            \mm \in \MM(b,k)
        }}
        \prod_{i = 0}^{b-1}
        \exp
        \left(
            -c \beta m_i 
        \right)
        \frac{1}{(b-1)!}
        \sum_{\bblocks \in \AA(b,\mm)}
        \P_{\beta}
        \left(
            \HH(\bblocks) \text{ is connected}
        \right).
    \end{multline}
\end{lemma}

\begin{proof}
    Consider the bound in \Cref{prop:block_cluster} and condition on the event $\{\BB(\cluster) = \bblocks\}$. 
    The connection event $\{\HH(\bblocks) \text{ is connected}\}$ is independent of the conditioning, while the isolation event $\{\bblocks \text{ is } \BB(\cluster) \text{-isolated}\}$ depends on the conditioning but is conditionally independent of the connection event.
    Just as in Lemma \ref{lem:independence}, the result follows from \Cref{lem:one_arm_isolation}.
\end{proof}

\subsection{Connectivity event: spanning trees}
In this section we bound the connectivity event in \eqref{eq:cluster_size_independence}.
For an arbitrary transitive kernel $J$, recall the rotationally symmetric kernel $\jmax$ associated to $J$ introduced in Section \ref{subsec:rotation_inv_ker}.

\begin{lemma}
    \label[lemma]{prop:cluster_size_connectivity}
    Let $G$ be a transitive graph of polynomial growth with $d \geq 2$, and suppose that $J : V \times V \to \R_+$ is a transitive kernel with $J(x,y) = O(d_G(x,y)^{-d \alpha})$ for $\alpha > 1$. 
    For $k \in \N$, $b \in \N$, and $\mm \in \MM(b,k)$, there exists $c = c(G) > 0$ such that
    \begin{equation}
        \frac{1}{(b-1)!}
        \sum_{\bblocks \in \AA(b,\mm)}
        \P_{\beta}
        \left(
            \HH(\bblocks) \text{ is connected}
        \right)
        \leq
        \prod_{i = 0}^{b-1}
        \left(
            c 
            \E_{\beta,\jmax}
            \left[
                \mathrm{length}(o)
            \right]
        \right)^{m_i}.
    \end{equation}
\end{lemma}

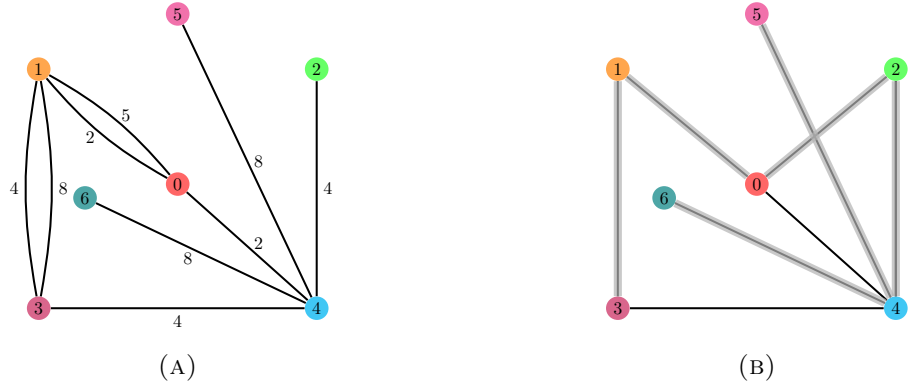
\begin{figure}
    \centering
    \hfill
    \begin{subfigure}[b]{0.32\textwidth}
        \centering
        \begin{tikzpicture}[
    block/.style={
        circle,
        minimum size=14pt,
        inner sep=2pt,
        draw=none
    },
    Acol/.style={fill=red!60},
    Bcol/.style={fill=orange!70},
    Ccol/.style={fill=green!60},
    Dcol/.style={fill=purple!60},
    EFcol/.style={fill=cyan!60},
    Gcol/.style={fill=teal!70},
    Hcol/.style={fill=magenta!70},
    edge/.style={line width=1.25pt}
]
\useasboundingbox (-4,-3) rectangle (4,4.5);

% ---------------------------------------------------------
% Nodes with labels
% ---------------------------------------------------------
\node[block, Acol]   (A)  at (0,0.3) {0};
\node[block, Bcol]   (B)  at (-3,2.8) {1};
\node[block, Ccol]   (C)  at (3,2.8) {2};
\node[block, Dcol]   (D)  at (-3,-2.4) {3};
\node[block, EFcol]  (EF) at (3,-2.4) {4};
\node[block, Hcol]   (H)  at (0,4) {5};
\node[block, Gcol]   (G)  at (-2,0) {6};

% ---------------------------------------------------------
% Edges with weights
% ---------------------------------------------------------
% Edges from 0 to 1 (weights 5 and 2)
\draw[edge, bend right=10] (A) to node[midway, right, yshift = 2pt] {5} (B);
\draw[edge, bend left=10]  (A) to node[midway, left, yshift = -2pt]  {2} (B);

% Edges from 1 to 3 (weights 4 and 8)
\draw[edge, bend right=10] (B) to node[midway, left] {4} (D);
\draw[edge, bend left=10]  (B) to node[midway, right]{8} (D);

% Edge from 3 to 4 (weight 4)
\draw[edge] (D) -- node[midway, below] {4} (EF);

% Edge from 4 to 2 (weight 4)
\draw[edge] (EF) -- node[midway, right] {4} (C);

% Edge from 6 to 4 (weight 8)
\draw[edge] (G) -- node[midway, left, xshift=-1pt, yshift = -3pt] {8} (EF);

% Edge from 0 to 4 (weight 2)
\draw[edge] (A) -- node[midway, right, yshift = 2pt] {2} (EF);

% Edge from 4 to 5 (weight 8)
\draw[edge] (EF) -- node[midway, right] {8} (H);

\end{tikzpicture}
        \caption{}
        \label{fig:block_and_spanning_graph}
    \end{subfigure}
    \hfill
    \begin{subfigure}[b]{0.32\textwidth}
        \centering
        \begin{tikzpicture}[
    block/.style={
        circle,
        minimum size=14pt,
        inner sep=2pt,
        draw=none
    },
    Acol/.style={fill=red!60},
    Bcol/.style={fill=orange!70},
    Ccol/.style={fill=green!60},
    Dcol/.style={fill=purple!60},
    EFcol/.style={fill=cyan!60},
    Gcol/.style={fill=teal!70},
    Hcol/.style={fill=magenta!70},
    edge/.style={line width=1.25pt},
    highlight/.style={line width=5pt, gray!60, opacity=0.7}
]
\useasboundingbox (-4,-3) rectangle (4,4.5);

% ---------------------------------------------------------
% Nodes with labels inside the nodes
% ---------------------------------------------------------
\node[block, Acol]   (A)  at (0,0.3) {0};
\node[block, Bcol]   (B)  at (-3,2.8) {1};
\node[block, Ccol]   (C)  at (3,2.8) {2};
\node[block, Dcol]   (D)  at (-3,-2.4) {3};
\node[block, EFcol]  (EF) at (3,-2.4) {4};
\node[block, Hcol]   (H)  at (0,4) {5};
\node[block, Gcol]   (G)  at (-2,0) {6};

% ---------------------------------------------------------
% Edges
% ---------------------------------------------------------
\draw[edge] (A) -- (C);
\draw[edge] (A) -- (EF);
\draw[edge] (B) -- (D);
\draw[edge] (G) -- (EF);
\draw[edge] (H) -- (EF);
\draw[edge] (B) -- (A);
\draw[edge] (EF) -- (C);
\draw[edge] (EF) -- (C);
\draw[edge] (EF) -- (D);

% ---------------------------------------------------------
% Highlighted edges
% ---------------------------------------------------------
\draw[highlight] (A) -- (B); % 0 to 1
\draw[highlight] (B) -- (D); % 1 to 3
\draw[highlight] (EF) -- (G); % 4 to 6
\draw[highlight] (EF) -- (C); % 4 to 2
\draw[highlight] (EF) -- (H); % 4 to 5
\draw[highlight] (A) -- (C); % 0 to 2

\end{tikzpicture}
        \caption{}
        \label{fig:weighted_block_graph}
    \end{subfigure}
    \hfill
    \caption{
        Figures \ref{fig:block_and_spanning_graph} and \ref{fig:weighted_block_graph} represent the block graph and the weighted block graph associated to the percolation configuration in Figure \ref{fig:lrp_configuration} with an arbitrary choice of labeling.
        In Figure \ref{fig:block_and_spanning_graph}, the highlighted subgraph is a spanning tree of the block graph, and the vector of forward degrees $\ff = (2,1,1,0,2,0,0)$ is an $\ff$-tree.
    }
    \label{fig:spanning_tree_example}
\end{figure}

We prove this lemma by considering rooted labelled spanning trees.
Let $b \in \N$ and consider some block decomposition $\BB(\cluster) \in \AA(b)$ with $b$ blocks.
We begin by choosing a rooted labelled spanning tree of the block graph $\HH(\BB(\cluster))$, so that the root is at the origin and the labels of the vertices are the labels of the block decomposition.
We let
\begin{equation}
    \label{eq:f_trees}
    \FF(b)
    = 
    \left\{
        (f_0,\ldots,f_{b-1})
        \in
        \N^b
        :
        \sum_{i = 0}^{b-1}
        f_i 
        = 
        b - 1,
        \sum_{i = 0}^{j}
        f_i 
        \geq j
        \text{ for all }
        j 
        \leq 
        b - 1
    \right\},
\end{equation}
we call $\ff \in \FF(b)$ the \textbf{vector of forward degrees}, and we say that a rooted labelled spanning tree of $\HH(\BB(\cluster))$ is an \textbf{$\ff$-tree} if its root has label $0$, the vertex with label $0$ has an edge to the vertices with labels $1,2,\ldots,f_0 - 1,f_0$, the vertex with label $1$ has an edge to the vertices with labels $f_0 + 1,f_0 +2 \ldots,f_0 + f_1 - 1, f_0 + f_1$, and inductively the vertex with label $j$ has an edge to each of the vertices with labels $\sum_{i = 0}^{j-1} f_i + 1, \sum_{i = 0}^{j-1} f_i + 2, \ldots, \sum_{i=1}^j f_i$.
In words, the integer $f_i$ corresponds to the number of outgoing edges from the vertex $i$ with respect to the spanning tree, and a rooted labelled spanning tree of $\HH(\BB(\cluster))$ is an $\ff$-tree if its labelling agrees with the vector of forward degrees $\ff$. 
Figure \ref{fig:spanning_tree_example} illustrates an $\ff$-tree on the block decomposition $\BB(\cluster)$ previously sketched in Figure \ref{fig:lrp_configuration}.
Finally, we say that the block graph $\HH(\BB(\cluster))$ is \textbf{$\ff$-connected} if it contains an $\ff$-tree on its vertices rooted at 0.

In the following lemma we reduce the question of bounding the probability that the block graph $\BB(\cluster)$ is connected to the question of bounding the probability that the block graph $\BB(\cluster)$ is $\ff$-connected for some vector of forward degrees.

\begin{lemma}
    \label[lemma]{lem:block_permutations}
    Let $G$ be a transitive graph of polynomial growth, and suppose that $J : V \times V \to \R_+$ is a transitive kernel with $J(x,y) = O(d_G(x,y)^{-d \alpha})$ with $\alpha > 1$.
    For $k \in \N$ and $\mm \in \MM(b,k)$,
    \begin{multline}
        \frac{1}{(b-1)!}
        \sum_{\bblocks \in \AA(b,\mm)}
        \P_{\beta}
        \left(
            \HH(\bblocks) \text{ is connected}
        \right)
        \\
        \leq 
        \sum_{\ff \in \FF(b)}
        \prod_{i = 1}^{b-1}
        \frac{1}{f_i !}
        \sum_{\bblocks \in \AA(b,\mm)}
        \P_{\beta}
        \left(
            \HH(\bblocks) \text{ is } \ff \text{-connected}
        \right).
    \end{multline}
\end{lemma}

\begin{proof}
    Suppose that $\BB(\cluster) = \bblocks = (B_0,\ldots,B_{b-1}) \in \AA(b,\mm)$, and write $S_{b-1}$ for the group of permutations on $(1,\ldots,b-1)$. 
    For $\ff \in \FF(b)$ and $\sigma \in S_{b-1}$, suppose that the block graph $\HH((B_0,B_{\sigma(1)},\ldots,B_{\sigma(b-1)}))$ is $\ff$-connected.
    By counting rooted isomorphisms, it follows that there are $\prod_{i=1}^{b-1} f_i!$ distinct pairs $(\ff',\sigma')$ with $\ff' \in \FF(b)$ and $\sigma' \in S_{b-1}$ such that $\HH((B_0,B_{\sigma'(1)},\ldots,B_{\sigma'(b-1)}))$ is $\ff'$-connected.
    Recall that the labelling of the block $B_0$ is fixed since $o \in B_0$.
    We have 
    \begin{align}
        \begin{split}
            &
            \frac{1}{(b-1)!}
            \sum_{\bblocks \in \AA(b,\mm)}
            \P_{\beta}
            \left(
                \HH(\bblocks) \text{ is connected}
            \right)
            \\
            &
            \leq 
            \frac{1}{(b-1)!}
            \sum_{\ff \in \FF(b)}
            \prod_{i = 1}^{b-1}
            \frac{1}{f_i !}
            \sum_{\bblocks \in \AA(b,\mm)}
            \sum_{\sigma \in S_{b-1}}
            \P_{\beta}
            \left(
                \HH((B_0,B_{\sigma'(1)},\ldots,B_{\sigma'(b-1)})) 
                \text{ is } 
                \ff \text{-connected} 
            \right)
        \end{split}
        \\
        &
        \leq 
        \frac{1}{(b-1)!}
        \sum_{\ff \in \FF(b)}
        \prod_{i = 1}^{b-1}
        \frac{1}{f_i !}
        \sum_{\bblocks \in \AA(b,\mm)}
        \P_{\beta}
        \left(
            \HH(\bblocks) \text{ is } \ff \text{-connected}
        \right)
        \sum_{\sigma \in S_{b-1}}
        1
    \end{align}
    where we use that \LRP{} is independent of the labelling.
    The permutation group $S_{b-1}$ has $\# S_{b-1} = (b-1)!$, so that the sum over $S_{b-1}$ and the factor $1/(b-1)!$ cancel to give the desired result.
\end{proof}

\begin{lemma}[Inductive step]
    \label[lemma]{lem:f_connectivity}
    Let $G$ be a transitive graph of polynomial growth with $d \geq 2$, and suppose that $J : V \times V \to \R_+$ is a transitive kernel with $J(x,y) = O(d_G(x,y)^{-d \alpha})$ with $\alpha > 1$.
    For $k \in \N$, $b \in \N, \mm \in \MM(b,k)$, and $\ff \in \FF(b)$,
    \begin{equation}
        \sum_{\bblocks \in \AA(b,\mm)}
        \P_{\beta}
        \left(
            \HH(\bblocks)
            \text{ is }
            \ff \text{-connected}
        \right)
        \leq 
        \E_{\beta,\jmax}
        \left[
            \mathrm{length}(o)
        \right]^{b-1}
        \prod_{i = 0}^{b-1}
        m_i^{f_{i}}
        \exp
        \left(
            \cpei m_i
        \right).
    \end{equation}
\end{lemma}

\begin{proof}
    We argue by induction on $b$.
    When $b = 1$, $m_0 \in \MM(1,k)$, and $\ff \in \FF(1)$, then $f_0 = 0$ and an $\ff$-tree is a tree on a single vertex. 
    By \Cref{lem:peierls_argument} it follows that
    \begin{align}
        \sum_{\bblocks \in \AA(1,m_0)}
        \P_{\beta}
        \left(
            \HH(\bblocks) \text{ is } \ff \text{-connected}
        \right)
        & =
        \# 
        \{
            B_0 \in \AA 
            :
            B_0 \ni o,
            \# \partial B_0
            =
            m_0
        \}
        \\
        & \leq
        \exp(\cpei m_0)
    \end{align}
    for $\cpei >0$, as required.
    We may assume inductively that the statement holds for $k \leq b-1$.
    Let $\mm \in \MM(b,k)$ and $\ff \in \FF(b)$.
    By construction, the vertex $b - 1$ is a leaf of the tree, so $f_{b-1} = 0$. 
    The parent of the vertex $b-1$ is the vertex $\ell$ where $\ell$ is the maximal index such that $f_{\ell} \neq 0$.
    The remaining labelled graph when we remove the vertex $b - 1$ is a tree with forward degree given by $\widehat{\ff} = (f_0, \ldots, f_{\ell - 1}, f_{\ell} - 1, f_{\ell + 1}, \ldots, f_{b-2}) \in \FF(b-1)$.
    For a fixed block $B_{\ell} \in \AA$ with $\# \partial B_{\ell} = m_{\ell}$, the Peierls argument in \Cref{lem:peierls_argument} yields
    \begin{align}
        \label{eq:spanning_tree_connectivity}
        \sum_{y \not \in B_{\ell}}
        \sum_{\substack{
            B_{b-1} \in \AA,
            B_{b-1} \ni y
            \\
            \# \partial B_{b-1} = m_{b-1}
        }}
        \P_{\beta}
        \left(
            B_{\ell} 
            \sim
            y
        \right)
        & =
        \sum_{\substack{
            x \in B_{\ell}
            \\
            y \not \in B_{\ell}
        }}
        \P_{\beta}
        \left(
            x
            \sim
            y
        \right)
        \sum_{\substack{
            B_{b-1} \in \AA,
            B_{b-1} \ni y
            \\
            \# \partial B_{b-1} = m_{b-1}
        }}
        1
        \\
        & \leq   
        \exp(\cpei m_{b-1})
        \sum_{\substack{
            x \in B_{\ell}
            \\
            y \not \in B_{\ell}
        }}
        \P_{\beta}
        \left(
            x
            \sim
            y
        \right).
    \end{align}
    Rewriting and using the isoperimetric inequality in \Cref{lem:want_to_show},
    \begin{align}
        \sum_{\substack{
            x \in B_{\ell}
            \\
            y \not \in B_{\ell}
        }}
        \P_{\beta}
        \left(
            x
            \sim
            y
        \right)
        =
        \sum_{r \in \N}
        \sum_{\substack{
            x \in B_{\ell}
            \\
            y \not \in B_{\ell}
        }}
        \P_{\beta}
        \left(
            x
            \jump{}
            y
        \right)
        & \leq
        \sum_{r \in \N}
        \max_{\substack{
            x,y \in G
            \\
            d_G(x,y) = r
        }}
        \P_{\beta}
        \left(
            x
            \jump{}
            y
        \right)
        \sum_{\substack{
            x \in B_{\ell}
            \\
            y \not \in B_{\ell}
        }}
        1
        \\
        & \leq 
        m_{\ell}
        \sum_{r \in \N}
        \max_{\substack{
            x,y \in G
            \\
            d_G(x,y) = r
        }}
        \P_{\beta}
        \left(
            x
            \jump{}
            y
        \right)
        r 
        \# S(r)
    \end{align}
  We use the rotationally symmetric kernel $\jmax$ associated to $J$, defined in Section \ref{subsec:rotation_inv_ker}, yielding
    \begin{equation}
        \sum_{r \in \N}
        \max_{\substack{
            x,y \in G
            \\
            d_G(x,y) = r
        }}
        \P_{\beta}
        \left(
            x
            \jump{}
            y
        \right)
        r 
        \# S(r)
        =
        \E_{\beta,\jmax}
        \left[
            \mathrm{length}(o)
        \right].
    \end{equation}
    We have shown that
    \begin{equation}
        \sum_{y \not \in B_{\ell}}
        \sum_{\substack{
            B_{b-1} \in \AA,
            B_{b-1} \ni y
            \\
            \# \partial B_{b-1} = m_{b-1}
        }}
        \P_{\beta}
        \left(
            B_{\ell}
            \sim
            y
        \right)
        \leq
        m_{\ell}
        \exp
        \left(
            \cpei{}
            m_{b-1}
        \right)
        \E_{\beta,\jmax}
        \left[
            \mathrm{length}(o)
        \right].
    \end{equation}
    Let $\widehat{\mm} = (m_0,\ldots,m_{b-2})$.
    By the independence of edges in \LRP{} and together with the induction hypothesis, we have shown that
    \begin{align}
        \begin{split}
            \sum_{\bblocks \in \AA(b,\mm)}
            & 
            \P_{\beta}
            \left(
                \HH(\bblocks) 
                \text{ is } 
                \ff 
                \text{-connected}
            \right)
            \\
            & \leq 
            \sum_{\bblocks \in \AA(b-1,\widehat{\mm})}
            \P_{\beta}
            \left(
                \HH(\bblocks) \text{ is } \widehat{\ff} \text{-connected}
            \right)
            \sum_{y \not\in B_{\ell}}
            \sum_{\substack{
                B_{b-1} \in \AA, B_{b-1} \ni y
                \\
                \# \partial B_{b-1} = m_{b-1}
            }}
            \P_{\beta}
            \left(
                B_{\ell} 
                \sim
                y
            \right)
        \end{split}
        \\
        \begin{split}
            & \leq
            \E_{\beta,\jmax}
            \left[
                \mathrm{length}(o)
            \right]^{b-2}
            \left(
                \prod_{i = 0}^{b-2}
                m_i^{\widehat f_{i}}
                \exp
                \left(
                    \cpei m_i
                \right)
            \right)
            \\
            & \qquad \qquad \qquad \qquad \qquad \qquad
            \times 
            m_{\ell}
            \exp
            \left(
                \cpei m_{b-1}
            \right)
            \E_{\beta,\jmax}
            \left[
                \mathrm{length}(o)
            \right]
        \end{split}
        \\
        & \leq 
        \E_{\beta,\jmax}
        \left[
            \mathrm{length}(o)
        \right]^{b-1}
        \prod_{i = 0}^{b-1}
        m_i^{f_{i}}
        \exp
        \left(
            \cpei m_i
        \right),
    \end{align}
    where we use that $\widehat f_{\ell} = f_{\ell} - 1$ and $f_{b-1} = 0$.
    This concludes the induction.
\end{proof}

\begin{proof}[Proof of \Cref{prop:cluster_size_connectivity}]
    Let $b \in \N$ and $\mm \in \MM(b,k)$.
    By \Cref{lem:block_permutations,lem:f_connectivity}, we have shown
    \begin{multline}
        \frac{1}{(b-1)!}
        \sum_{\bblocks \in \AA(b,\mm)}
        \P_{\beta}
        \left(
            \HH(\bblocks) \text{ is connected}
        \right)
        \\
        \leq
        \E_{\beta,\jmax}
        \left[
            \text{length}(o)
        \right]^{b-1}
        \sum_{\ff \in \FF(b)}
        \prod_{i = 0}^{b-1}
        \frac{
            m_i^{f_{i}}
            \exp
            \left(
                \cpei m_i
            \right)
        }{
            f_i!
        }.
    \end{multline}
    It follows from considering the series expansion of $\exp(m)$ that $m^{f}/f! \leq \exp(m)$ for $m,f \in \N$, and hence 
    \begin{equation}
        \prod_{i = 0}^{b-1}
        \frac{
            m_i^{f_{i}}
            \exp
            \left(
                \cpei m_i
            \right)
        }{
            f_i!
        }
        \leq 
        \prod_{i = 0}^{b-1}
        \exp
        \left(
            \left(
                \cpei
                + 
                1
            \right)
            m_i
        \right).
    \end{equation}
    Recall the definition of $\FF(b)$ in \eqref{eq:f_trees}.
    There are at most $\binom{2b-2}{b-1} \leq 2^{2b} \leq \exp(2b)$ ways to write the integer $b-1$ as a sum of $b$ non-negative integers (the sequence of forward degrees in $\ff$), and since $b\le \sum_{i=0}^{b-1} m_i$, we have $\# \FF(b) \leq \exp(2 \sum_{i = 0}^{b-1} m_i)$.
    We have shown that
    \begin{align}
        \begin{split}
            \frac{1}{(b-1)!}
            \sum_{\bblocks \in \AA(b,\mm)}
            \P_{\beta}
            (
                &
                \HH(\bblocks) \text{ is connected}
            )
            \\
            & \leq 
            \E_{\beta,\jmax}
            \left[
                \mathrm{length}(o)
            \right]^{b-1}
            \prod_{i = 0}^{b-1}
            \exp
            \left(
                \left(
                    \cpei
                    +
                    3
                \right)
                m_i
            \right)
        \end{split}
        \\
        & \leq
        \prod_{i = 0}^{b-1}
        \left(
            c
            \E_{\beta,\jmax}
            \left[
                \mathrm{length}(o)
            \right]
        \right)^{m_i}
    \end{align}
    for some $c = c(G) > 0$, concluding the proof.
\end{proof}

\subsection{Cluster-size decay: upper bound Theorem~\ref{thm:cluster_size_strong_decay}}
\begin{proof}[Proof of \Cref{thm:cluster_size_strong_decay}]
    By \Cref{prop:cluster_size_independence,prop:cluster_size_connectivity} we have shown that
        \begin{equation}
        \label{eq:cluster_size_series_bound}
        \P_{\beta}
        \left(
            k
            \leq 
            \# 
            K
            < 
            \infty
        \right)
        \leq
        \sum_{b \in \N}
        \sum_{\mm \in \MM(b,k)}
        \prod_{i = 0}^{b-1}
        \left(
            c_1
            \E_{\beta,\jmax}
            \left[
                \mathrm{length}(o)
            \right]
            \exp
            \left(
                - 
                c_2 
                \beta
            \right)
        \right)^{m_i}
    \end{equation} 
    for some $c_1 = c_1(G),c_2 = c_2(G,J) > 0$.
    By \Cref{lem:finite_exp_leng}, we have that $\E_{\beta,\jmax}\left[ \mathrm{length}(o)\right] < \infty$ if and only if $\alpha > 1 + 1/d$. 
    For these values of $\alpha$ we have
    $
        \E_{\beta,\jmax}
        \left[
            \mathrm{length}(o)
        \right]
        \leq 
        c_3
        \beta
    $
    for some $c_3 =c_3(G,J) > 0$.
    Recall the definition of $\MM(b,k)$ in \eqref{eq:cluster_size_blocks}.
    If we write $\ell = m_0 + \ldots + m_{b-1}$, then $\ell \geq b$ and $\ell \geq \ciso k^{(d-1)/d}$. 
    There are at most $\binom{\ell - 1}{b - 1} \leq 2^{\ell}$ ways to write an integer $\ell$ as a sum of $b$ strictly positive integers and hence
    \begin{align}
        \P_{\beta}
        \left(
            k
            \leq 
            \# 
            K
            < 
            \infty
        \right)
        & \leq
        \sum_{\ell = \lfloor \ciso k^{(d-1)/d} \rfloor}^{\infty}
        \sum_{b = 1}^\ell
        \sum_{\substack{
            \mm \in \N^b
            \\
            \sum_{i \leq b} m_i = \ell
        }}
        \prod_{i = 0}^{b-1}
        \left(
            c_4
            \beta
            \exp
            \left(
                - 
                c_2 
                \beta
            \right)
        \right)^{m_i}
        \\
        & \leq 
        \sum_{\ell = \lfloor \ciso k^{(d-1)/d} \rfloor}^{\infty}
        \left(
            c_4
            \beta
            \exp
            \left(
                - 
                c_2 
                \beta
            \right)
        \right)^\ell
        \sum_{b = 1}^\ell
        \sum_{\substack{
            \mm \in \N^b
            \\
            \sum_{i = 0}^b m_i = \ell
        }}
        1
        \\
        & \leq
        \sum_{\ell = \lfloor \ciso k^{(d-1)/d} \rfloor}^{\infty}
        \left(
            2
            c_4
            \beta
            \exp
            \left(
                - 
                c_2 
                \beta
            \right)
        \right)^{\ell}
    \end{align}  
    for some $c_4 = c_4(G,J) > 0$.
    Let $c_5 = c_5(G,J) > 0$ be such that 
    $2 c_4 \beta \exp \left(- c_2 \beta \right) \leq \exp(-c_5 \beta)$ for all $\beta$ sufficiently large (the threshold depending only on $c_2,c_4, c_5$). 
    The sum over $\ell$ can be bounded from above by a geometric series with ratio $q(\beta) = \exp\left(- c_5 \beta \right)$. 
    We can choose $\beta$ sufficiently large such that $0 < q(\beta) < 1/2 < 1$, and hence
    \begin{equation}
        \P_{\beta}
        \left(
            k
            \leq 
            \# 
            K
            < 
            \infty
        \right)
        \leq
        2
            \exp
            \left(
                - 
                c_5 
                \beta
        \right)^{\lfloor \ciso k^{(d-1)/d} \rfloor}
        \leq
        \exp
        \left(
            -
            \beta k^{(d-1)/d}
            /
            A
        \right)
    \end{equation}
    for some $A = A(G,J) > 0$. This concludes the proof.
\end{proof}

\begin{remark}
    \label[remark]{rem:crucial_alpha_assumption}
    The assumption that $\alpha > 1 + 1/d$ was necessary in the proof of \Cref{thm:cluster_size_strong_decay} for the expected total edge-length $\E_{\beta,\jmax} \left[ \mathrm{length}(o) \right]$ to be finite.
    The difference with the argument for the truncated one-arm, where the result holds for all $\alpha > 1$, occurs at the induction step, \Cref{lem:sum_over_blocks,lem:f_connectivity} respectively.
    In the truncated one-arm, the path we are considering has fixed outgoing degree $1$: in \Cref{lem:sum_over_blocks} we used the  $d$-dimensional isoperimetric inequality, and the term $m_i^{d/(d-1)}$ can be absorbed by the term coming from the Peierls argument.
    In the cluster-size, the outgoing degree $f_i$ is unbounded, and if we were to use the $d$-dimensional isoperimetric inequality the additional term $m_i^{d/(d-1)}$ cannot be absorbed. 
    Intuitively, the isoperimetric inequality in Proposition \ref{lem:want_to_show} is stronger than that in Lemma \ref{lem:trivial_isoperimetry} for short edges, and in the cluster-size we are considering unboundedly many short edges.
\end{remark}

\subsection{Finite set connects to infinity}
A modification of the proof of \Cref{thm:cluster_size_strong_decay} gives the following result.
For $a > 0$, let 
\begin{equation}
    \SS(a) 
    = 
    \{
        S \subset V 
        : 
        S 
        \text{ is finite with } 
        \ell
        \text{ components where } 
        \ell
        \log 
        \ell
        \leq 
        a 
        (\#S)^{(d-1)/d} 
        \}.
\end{equation}

\begin{theorem}
    \label[theorem]{thm:s_few_components}
    Let $G$ be a transitive graph of polynomial growth with $d \geq 2$, and suppose that $J: V \times V \to \R_+$ is a transitive kernel with $J(x,y) = O(d_G(x,y)^{-d\alpha})$ with $\alpha > 1 + 1/d$. 
    For $a > 0$ and $\beta > \beta_c$ sufficiently large, there exists $A = A(G,J) > 0$ such that for every $S \in \SS(a)$
    \begin{equation}
        \P_{\beta}
        \left(
            S 
            \not\leftrightarrow
            \infty
        \right) 
        \leq 
        \exp
        \left(
            - 
            \beta (\# 
            S)^{(d-1)/d}
            /
            A
        \right).
    \end{equation}
\end{theorem}

\begin{proof}
    Let $a > 0$, let $S \in \SS(a)$, and let $S_1,\ldots,S_\ell$ be the components of $S$. 
    By the translation invariance of \LRP{}, we may assume without loss of generality that $o \in S$. 
    Slightly abusing notation, we write $K(S)=\cup_{v\in S}K(v)$ for the union of the clusters of $S$, which we consider as a set of vertices of $G$.
    Recall the definitions of the block decomposition and the associated block graph in Section \ref{subsec:block_decomp}.
    The event $\{S \not\leftrightarrow \infty\}$ is equivalent to the event $\{ \# K(S) < \infty\}$, and in terms of blocks the event $\{ \# K(S) < \infty\}$ implies that the block decomposition $\BB(K(S))$ consists of a sequence of blocks containing $S$ such that the block graph $\HH(\BB(K))$ consists of $q$ components for some $q \leq \ell$ and with each component containing at least one of the $S_i$.
    Let $\PP(\ell,q)$ be the set of partitions of $\{1,\ldots,\ell\}$ into $q$ many sets, so that $P \in \PP(\ell,q)$ is given by $P = \{P_1,\ldots,P_q\}$ with $P_i \subseteq \{1,\ldots,\ell\}$ for all $i$, $P_i \cap P_j = \varnothing$ for all $i \neq j$, and $\cup_i P_i = \{1, \ldots, \ell\}$.
    For $b,k \in \N$, recall the definition of $\MM(b,k)$ the set of possible boundary sizes of the blocks in \eqref{eq:cluster_size_blocks} and $\AA(b,\mm)$ the set of possible block decompositions with given boundary size in \eqref{eq:possible_block_decompositions}.
    For each $i \leq q$, let $T^{(i)} = \sum_{j \in P_i} \# S_j$. 
    Then 
    \begin{multline}    
        \P_{\beta}
        \left(
            S 
            \not\leftrightarrow
            \infty
        \right)
        \leq
        \sum_{\substack{
            q \leq \ell
            \\
            P \in \PP(\ell,q)
        }}
        \sum_{\substack{
            i \leq q,
            b^{(i)} \in \N
            \\
            \mm^{(i)} \in \MM(b^{(i)},T^{(i)})
        }}
        \left(
            \prod_{i \leq q}
            \frac{1}{(b^{(i)} - 1)!}
        \right)
        \sum_{\bblocks^{(i)} \in \AA(b^{(i)},\mm^{(i)})}
        \\
        \P_{\beta}
        \big(
            \BB(\cluster(S))
            = 
            \cup_{i \leq q} \bblocks^{(i)},
            \HH(
                \bblocks^{(i)}   
            )
            \text{ is connected for each } i \leq q,
            \\
            \bblocks^{(i)} 
            \text{ is }
            \BB(\cluster(S))
            \text{-isolated for each } i \leq q
        \big).
    \end{multline}
    This bound is similar to that in \Cref{prop:block_cluster}, except that in this setting we are requiring that the block graph consists of $q$ components rather than being a connected graph. 
    We sum over partitions of $S_1,\ldots,S_\ell$ so that each of the $q$ components of the block graph corresponds to one of these partitions, and for each component we sum over all possible connected block graphs respecting that the block graph must contain the required sets $S_j$. 
    In particular, the indexing over $i$ above corresponds to the individual components of the block graph. 
    As the blocks $ \bblocks^{(i)}$ contain the sets $S_j$ for $j \in P_i$, there are in total at least $T^{(i)}=\sum_{j\in P_j}\#S_j$ vertices in the blocks $\bblocks^{(i)}$. As a result of the isoperimetric inequality, the total boundary length of the blocks $\sum_{j=0}^{b^{(i)}-1} m_j^{(i)} \ge  c_{\mathrm{iso}} (T^{(i)})^{(d-1)/d}$. This means that the boundary lengths must be in the set $\mathcal M(b^{(i)}, T^{(i)})$ in this case, see more on this around \eqref{eq:concavity_argument}.
    The remainder of the proof follows the proof of \Cref{thm:cluster_size_strong_decay} with some minor differences which we now point out. 
    The independence of events is similar to as in \Cref{prop:cluster_size_independence}, except that the events that the blocks in the different components are isolated are not independent, but the upper bound from \Cref{lem:one_arm_isolation} factorises.
    This yields
    \begin{multline}    
        \P_{\beta}
        \left(
            S 
            \not\leftrightarrow
            \infty
        \right)
        \leq
        \sum_{\substack{
            q \leq \ell
            \\
            P \in \PP(\ell,q)
        }}
        \prod_{i \leq q}
        \sum_{\substack{
            b^{(i)} \in \N
            \\
            \mm^{(i)} \in \MM(b^{(i)},T^{(i)})
        }}
        \frac{
            1
        }{
            (b^{(i)} - 1)!
        }
        \\
        \times
        \exp\left(-c \beta \sum_{j=0}^{b^{(i)}-1}m_j^{(i)} \right)
        \sum_{\bblocks^{(i)} \in \AA(b^{(i)},\mm^{(i)})}
        \P_{\beta}
        \left(
            \HH(
                \bblocks^{(i)}   
            )
            \text{ is connected}
        \right)
    \end{multline}
    for $c = c(G, J) > 0$ as in \Cref{lem:one_arm_isolation}.
    From here, \Cref{prop:cluster_size_connectivity,lem:block_permutations,lem:f_connectivity} and the arguments in the proof of \Cref{thm:cluster_size_strong_decay} apply nearly verbatim to give
    \begin{align}
        \P_{\beta}
        \left(
            S 
            \not\leftrightarrow
            \infty
        \right)
        & \leq
        \sum_{\substack{
            q \leq \ell
            \\
            P \in \PP(\ell,q)
        }}
        \prod_{i \leq q}
        \exp
        \left(
            -\beta (T^{(i)})^{d/(d-1)} 
            /
            A_1
        \right)
        \\
        & \leq
        \sum_{\substack{
            q \leq \ell
            \\
            P \in \PP(\ell,q) 
        }}
        \exp
        \left(
            -\beta \# S
            ^{d/(d-1)} 
            /
            A_1
        \right) 
       \end{align}
    for some $A_1 = A_1(G,J) > 0$, where in the second inequality we use that $\sum_{i \leq q} T^{(i)} = \# S$ and that the function $f(x) = x^{(d-1)/d}$ is concave and increasing.
    The number of partitions of $\{1,\ldots,\ell\}$ into $q$ many sets is at most $q^\ell \leq \ell^\ell$.
    Provided that $\ell \log \ell \leq a \# S^{(d-1)/d}$, as is the case for sets $S \in \SS(a)$, we can choose $\beta$ sufficiently large such that $a<\beta/(2A_1)$ and hence
    \begin{equation}
        \P_{\beta}
        \left(
            S 
            \not\leftrightarrow
            \infty
        \right)
        \leq
        \exp 
        \left(
            -\beta  
            \# 
            S^{(d-1)/d}
            /
            A_2
        \right)
    \end{equation}
    for some $A_2 = A_2(G,J) > 0$, as required. 
\end{proof}

\begin{proof}[Proof of \Cref{thm:s_connects_infty}]
    This follows immediately from \Cref{thm:s_few_components}.
\end{proof}

\bibliographystyle{alpha}
\bibliography{references.bib}

\end{document}